\documentclass[12pt]{article}
\usepackage{amsmath}
\usepackage{amssymb}
\usepackage{amsfonts}
\usepackage{upref}
\usepackage{theorem}

\usepackage{psfig}

\usepackage{mathrsfs}      
\usepackage{wasysym}       

\title
{\huge\bf 
$\,$\\[-5ex]
Optimal Tristance Anticodes in~Certain~Graphs
}

\author{\\ 
        {\bf Tuvi Etzion}\\ 
        \small Technion --- Israel Institute of Technology \vspace*{-1.2ex}\\
        \small Department of Computer Science \vspace*{-1.2ex}\\
        \small Technion City, Haifa 32000, Israel \vspace*{-0.75ex}\\
        \small\tt etzion@cs.technion.ac.il\\
   \and \\ 
        {\bf Moshe Schwartz}\\ 
        \small Technion --- Israel Institute of Technology \vspace*{-1.2ex}\\
        \small Department of Computer Science \vspace*{-1.2ex}\\
        \small Technion City, Haifa 32000, Israel \vspace*{-0.75ex}\\
        \small\tt moosh@cs.technion.ac.il\\
   \and \\ \\[-1.25ex]
     {\bf Alexander Vardy}${}^{*}$\\
     \small University of California San Diego\vspace*{-1.0ex}\\
     \small Department of Computer Science\vspace*{-1.0ex}\\
     \small Department of Electrical Engineering\vspace*{-1.0ex}\\
     \small Department of Mathematics\vspace*{-1.0ex}\\
     \small 9500 Gilman Drive, 
            La Jolla, CA 92093, U.S.A.\vspace*{-.75ex}\\
        \small\tt vardy@kilimanjaro.ucsd.edu \\ \\[.5ex]  
       }

\date{June 5, 2004}

\theoremstyle{plain} 
\theorembodyfont{\normalfont\slshape}

\newtheorem{thm}{Theorem}  
\newenvironment{theorem}{\begin{thm}\hspace*{-1.25ex}{\bf.}}{\end{thm}}
\newcommand{\Tref}[1]{Theo\-rem\,\ref{#1}}

\newtheorem{lem}[thm]{Lemma} 
\newenvironment{lemma}{\begin{lem}\hspace*{-1.25ex}{\bf.}}{\end{lem}}
\newcommand{\Lref}[1]{Lemma\,\ref{#1}}

\newtheorem{prop}[thm]{Proposition}

\newtheorem{cor}[thm]{Corollary} 
\newenvironment{corollary}{\begin{cor}\hspace*{-1.25ex}{\bf.}}{\end{cor}}
\newcommand{\Cref}[1]{Corollary\,\ref{#1}}

\newtheorem{defi}{Definition}
\newenvironment{definition}{\begin{defi}\hspace*{-1ex}{\bf.}}{\end{defi}}
\newcommand{\Dref}[1]{Definition\,\ref{#1}}

\newtheorem{conj}[thm]{Conjecture} 
\newenvironment{conjecture}{\begin{conj}\hspace*{-1.25ex}{\bf.}}{\end{conj}}

\theorembodyfont{\normalfont}

\newtheorem{exam}{Example}

\newenvironment{proof}  
{$\,$\vspace{-4ex}\\ \hspace*{15pt}{\it Proof.}}{}

\newcommand{\dfn}{\sffamily\slshape}
\newcommand{\deff}{\mbox{$\stackrel{\rm def}{=}$}}

\renewcommand{\em}{\sl}
\renewcommand{\emph}[1]{{\em {#1}\/}}

\renewcommand{\le}{\leqslant}
\renewcommand{\leq}{\leqslant}
\renewcommand{\ge}{\geqslant}

\newcommand{\bc}{\begin{center}}
\newcommand{\ec}{\end{center}}

\newcommand{\be}{\begin{equation}}
\newcommand{\ee}{\end{equation}}

\newcommand{\Eq}[1]{\begin{equation} #1 \end{equation}}
\newcommand{\eq}[1]{{\rm (\ref{#1})}}

\newcommand{\abs}[1]{\left|#1\right|}

\newcommand{\set}[1]{\left\{#1\right\}}

\newcommand{\ceilenv}[1]{\left\lceil #1 \right\rceil}
\newcommand{\floorenv}[1]{\left\lfloor #1 \right\rfloor}

\newcommand{\shalf}[1]
{\mbox{\raisebox{.8mm}{\footnotesize $\scriptstyle #1$} 
       \footnotesize$\!\!\! / \!\!\!$
       \raisebox{-.8mm}{\footnotesize $\scriptstyle 2$}\kern 1pt}
}

\newcommand{\sthird}[1]
{\mbox{\raisebox{.8mm}{\footnotesize $\scriptstyle #1$} 
       \footnotesize$\!\!\! / \!\!\!$
       \raisebox{-.8mm}{\footnotesize $\scriptstyle 3$}\kern 1pt}
}

\newcommand{\sqrthalf}[1]
{\mbox{\raisebox{.8mm}{\footnotesize $\scriptstyle \sqrt{#1}$} 
       \footnotesize$\!\!\! / \!\!\!$
       \raisebox{-.8mm}{\footnotesize $\scriptstyle 2$}\kern 1pt}
}

\newcommand{\Z}{\mathbb{Z}}
\newcommand{\R}{\mathbb{R}}
\newcommand{\C}{\mathbb{C}}

\newcommand{\Strut}[2]{\rule[-#2]{0cm}{#1}}
\newcommand{\inv}{{}^{\raisebox{-1ex}{$\scriptstyle -1$}}}

\newcommand{\al}{\alpha}

\newcommand{\icos}{ico\-si\-hexa\-hedron}
\newcommand{\icoses}{ico\-si\-hexa\-hedra}

\DeclareFontFamily{U}{ding}{}
\DeclareFontShape{U}{ding}{m}{n}{  <-> bbding10 }{}

\newcommand{\plus}{\mbox{\tiny$\rm \boxplus$}}
\renewcommand{\star}{\raisebox{.25ex}{\tiny\rm $\infty$}}
\newcommand{\hex}{\raisebox{0.20ex}{\tiny\rm $\hexagon$}}

\newcommand{\tplus}{\raisebox{-.30ex}{\large $\rm \boxplus$}}
\newcommand{\tstar}{\raisebox{+.00ex}{\large $\infty$}}
\newcommand{\thex}{\raisebox{-.30ex}{\Large $\rm \hexagon$}}

\newcommand{\Gplus}{\cG^{\plus}_2}
\newcommand{\Gplusthree}{\cG^{\plus}_3}

\newcommand{\Gstar}{\cG^{\star}_2}
\newcommand{\Ghex}{\cG^{\hex}_2}

\newcommand{\dstar}{d^{\star}_3}
\newcommand{\dhex}{d^{\kern1pt\hex}_3}
\newcommand{\dtwohex}{d^{\kern1pt\hex\kern-1pt}}

\newcommand{\xmax}{x_{\max}}
\newcommand{\ymax}{y_{\max}}
\newcommand{\zmax}{z_{\max}}
\newcommand{\xmin}{x_{\min}}
\newcommand{\ymin}{y_{\min}}
\newcommand{\zmin}{z_{\min}}
\newcommand{\xmid}{x_{\rm mid}}
\newcommand{\ymid}{y_{\rm mid}}

\DeclareMathOperator{\hmax}{\overline{\max}}
\DeclareMathOperator{\hmin}{\underline{\min}}

\newcommand{\x}{{\sf x}}
\newcommand{\y}{{\sf y}}
\newcommand{\z}{{\sf z}}
\newcommand{\bx}{\bar{\sf x}}
\newcommand{\by}{\bar{\sf y}}
\newcommand{\bz}{\bar{\sf z}}

\newcommand{\p}{{\sf \chi}}
\newcommand{\q}{{\sf \psi}}
\renewcommand{\r}{{\sf \omega}}
\newcommand{\bp}{\bar{\sf \chi}}
\newcommand{\bq}{\bar{\sf \psi}}
\newcommand{\br}{\bar{\sf \omega}}

\newcommand{\cA}{{\cal A}}

\newcommand{\cC}{{\cal C}}

\newcommand{\cE}{{\cal E}}

\newcommand{\cG}{{\cal G}}

\newcommand{\cI}{{\cal I}}

\newcommand{\cP}{{\cal P}}

\newcommand{\cR}{{\cal R}}
\newcommand{\cS}{{\cal S}}

\newcommand{\cV}{{\cal V}}

\newcommand{\sD}{{\mathscr D}}

\newcommand{\sH}{{\mathscr H}}
\newcommand{\sI}{{\mathscr I}}

\newcommand{\sO}{{\mathscr O}}

\newcommand{\sS}{{\mathscr S\kern -4pt}}

\outer\def\proclaim #1. #2\par{\medbreak
 \noindent{\bf#1.\enspace}{\sl#2\par}%
 \ifdim\lastskip<\medskipamount \removelastskip\penalty55\medskip\fi}

\def\qed{\hskip 3pt \hbox{\vrule width4pt depth2pt height6pt}}

\mathchardef\inn="3232

\makeatletter

\gdef\@punct{.\ \ }  
\def\@sect#1#2#3#4#5#6[#7]#8{%
  \ifnum #2>\c@secnumdepth
     \def\@svsec{}
  \else
     \refstepcounter{#1}\edef\@svsec{%
     \ifnum #2>0{{\csname the#1\endcsname}}.\fi%
    \hskip .5em}
  \fi
  \@tempskipa #5\relax
  \ifdim \@tempskipa>\z@
     \begingroup #6\relax
       \@hangfrom{\hskip #3\relax\@svsec}{\interlinepenalty \@M #8\par}
     \endgroup
     \csname #1mark\endcsname{#7}
     \addcontentsline{toc}{#1}{\ifnum #2>\c@secnumdepth\else
          \protect\numberline{\csname the#1\endcsname}\fi#7}
  \else
     \def\@svsechd{#6\hskip #3\@svsec #8\@punct\csname #1mark\endcsname{#7}
     \addcontentsline{toc}{#1}{\ifnum #2>\c@secnumdepth \else
          \protect\numberline{\csname the#1\endcsname}\fi#7}}
  \fi
  \@xsect{#5}}

\def\@ssect#1#2#3#4#5{\@tempskipa #3\relax
  \ifdim \@tempskipa>\z@
     \begingroup #4\@hangfrom{\hskip #1}{\interlinepenalty \@M #5\par}\endgroup
  \else \def\@svsechd{#4\hskip #1\relax #5\@punct}\fi
  \@xsect{#3}}

\makeatother

\begin{document}


\maketitle
\vspace{-1.25ex}

\begin{center}
{\bf Abstract} \vspace{-2.5ex}
\end{center}
\begin{list}{}
{
\addtolength{\leftmargin}{-1.5ex}
\setlength{\rightmargin}{\leftmargin}
}
\item\noindent
\baselineskip2.6ex
\looseness=-1
For $z_1,z_2,z_3 \in \Z^n$, the {\dfn tristance\/} $d_3(z_1,z_2,z_3)$ is
a generalization of the $L_1$-distance on $\Z^n$ to a quantity 
that reflects the relative dispersion of three points rather than two.
A {\dfn tristance anticode\/} $\cA_d$ of diameter $d$ is a subset 
of $\Z^n$ with the property that $d_3(z_1,z_2,z_3) \le d$ for all
$z_1,z_2,z_3 \in \cA_d$. 
An anticode is {\dfn optimal\/} if it has the largest possible
cardinality for its diameter $d$. 
We determine the cardinality and completely classify the optimal tristance 
anticodes in $\Z^2$ for all diameters $d \ge 1$. 
We then generalize this result to two related distance
models: a different distance structure on $\Z^2$ where $d(z_1,z_2) = 1$
if $z_1,z_2$ are adjacent either horizontally, vertically, 
or {\sl diagonally}, and the distance structure obtained when
$\Z^2$ is replaced~by the hexagonal lattice $A_2$.
We also investigate optimal tristance anticodes~in~$\Z^3$
and optimal quadristance anticodes in $\Z^2$, and provide bounds
on their~cardinal\-ity.
We conclude 
with a brief discussion of the applications
of our results to multi-dimensional interleaving schemes
and to connectivity loci in the game of Go. 
\end{list}



\footnotetext[1]
{\\
\hspace*{-1ex}${}^*$Research supported 
by the David and Lucile Packard Foundation
and by the National Science Foundation.\vspace*{-9ex}
}

 
\thispagestyle{empty}
\addtocounter{page}{-1}
\newpage

\section{Introduction}  
\vspace{-1.5ex} 
\label{sec1}

\noindent
Given two points $z = (z_1,z_2,\ldots,z_n)$ and 
$z' = (z'_1,z'_1,\ldots,z'_n)$ in $\Z^n$, the 
{$L_1$-distance\/}~between $z$ and $z'$ 
is defined as 
$
d(z,z') = |z_1 - z'_1| + |z_2 - z'_2| + \cdots + |z_n - z'_n|
$.
Alternatively, let~$\cG^{\plus}_n = (V,E)$ denote the 
{\dfn grid graph\/} of $\Z^n$ whose vertex set is $V\! = \Z^n$
and whose edges are defined as follows: $\set{z,z'} \in E$
if and only if $d(z,z') = 1$. Then the $L_1$-distance between
$z$ and $z'$ in $\Z^n$ is the number of edges in the shortest
path joining $z$ and $z'$ in~$\cG^{\plus}_n$. 
The latter point of view leads to a natural generalization of 
the $L_1$-distance on $\Z^n$ to a~quantity that reflects the 
relative dispersion of three points rather than two.

\begin{definition}
\label{tristance-def}
Let $z_1,z_2,z_3 \in \Z^n$. Then the {\dfn tristance\/}
$d_3(z_1,z_2,z_3)$ is defined as the number of edges in 
a minimal spanning tree for $z_1,z_2,z_3$ in the grid
graph $\cG^{\plus}_n$ of $\Z^n$. 
\end{definition}

\looseness=-1
The notion of tristance defined above can be further generalized 
in two different ways. First, the {\dfn quadristance} $d_4(z_1,z_2,z_3,z_4)$,
the {\dfn quintistance} $d_5(z_1,z_2,z_3,z_4,z_5)$, and more generally
the {\dfn $r$-dispersion\/} $d_r(z_1,z_2,\ldots,z_r)$ may be 
defined~\cite{EV,ME,SE2} as the
number of edges in a~minimal spanning tree for $z_1,z_2,\ldots,z_r$ 
in the grid graph $\cG^{\plus}_n$. In this paper, we consider only
the tristance $d_3(z_1,z_2,z_3)$ and, briefly, the quadristance 
$d_4(z_1,z_2,z_3,z_4)$ in \S4.2.

Another way to generalize \Dref{tristance-def} is to replace the 
grid graph $\cG^{\plus}_n$ by a different~graph. We will consider two 
alternative graphs $\Gstar$ and $\Ghex$ that are useful in applications 
to two-dimensional interleaving~\cite{BBF,BBV,EV,SE2}.
The graph $\Gstar$ has $\Z^2$ as its vertex set, with $z = (x,y)$ 
and $z' = (x',y')$ in $\Z^2$ being adjacent if and only if
$$
d^{\star}(z,z') \ \ \deff\ \ \max\set{|x-x'|,|y-y'|} \ = \ 1
$$
Thus tristance in $\Gstar$ may be thought of as a generalization 
of the $L_{\infty}$-distance on $\Z^2$.
The vertex set of the graph $\Ghex$ is the hexagonal lattice 
$A_2 = \set{ (\shalf{1} v, u \,{+}\, \sqrthalf{3} v) : \,u,v \in \Z}$,
with two points $z = (x,y)$ and $z' = (x',y')$ in $A_2$ being adjacent 
iff 
$$
d_{\rm E}(z,z') \ \ \deff\ \ \sqrt{(x-x')^2 + (y-y')^2} \ = \ 1
$$
The three graphs $\Gplus$, $\Gstar$, and $\Ghex$ are illustrated in Figure\,1.
We will \mbox{sometimes~refer~to} 
the graphs $\Gplus$, $\Gstar$, and $\Ghex$ 
as the {grid graph}, the {\dfn infinity graph}, 
and the {\dfn hexagonal graph}, respectively 
(or the $\tplus$ model, the $\tstar$ model, 
and the $\thex$ model, for short).

\looseness=-1
Given a set $\cS$ and a definition of distance between points of $\cS$,
a {\dfn code\/} $\cC \subseteq \cS$ of minimum distance $d$ is characterized
by the property that the distance between any two distinct points of $\cC$
is \emph{at least} $d$. Similarly, an {\dfn anticode\/} $\cA \subseteq \cS$
of diameter $d$ is characterized by the property that the distance between 
any two distinct points of $\cA$ is \emph{at most} $d$. One is usually 
interested in codes and anticodes of the largest possible cardinality for 
a~given minimum distance or diameter --- such codes/anticodes are said
to be {\dfn optimal}. An encyclopedic survey of optimal codes in the
Hamming graph may be found in~\cite{MWS,Handbook}; for codes in other
graphs, see~\cite{Bannai,Biggs,Delsarte,Kratochvil,SE,SE2}. Optimal anticodes
in the Hamming metric and related distance-regular graphs have been 
studied in~\cite{AK1,AK2,Delsarte,MZ,SE} and other papers.

\looseness=-1
The concepts of a code and an anticode can be generalized using the 
notion of tristance in \Dref{tristance-def}. Thus a tristance code 
$\cC \subseteq \Z^n$ of minimum tristance $d$ 
is a subset of $\Z^n$ such 

\hspace*{0.10in}%
\psfig{figure=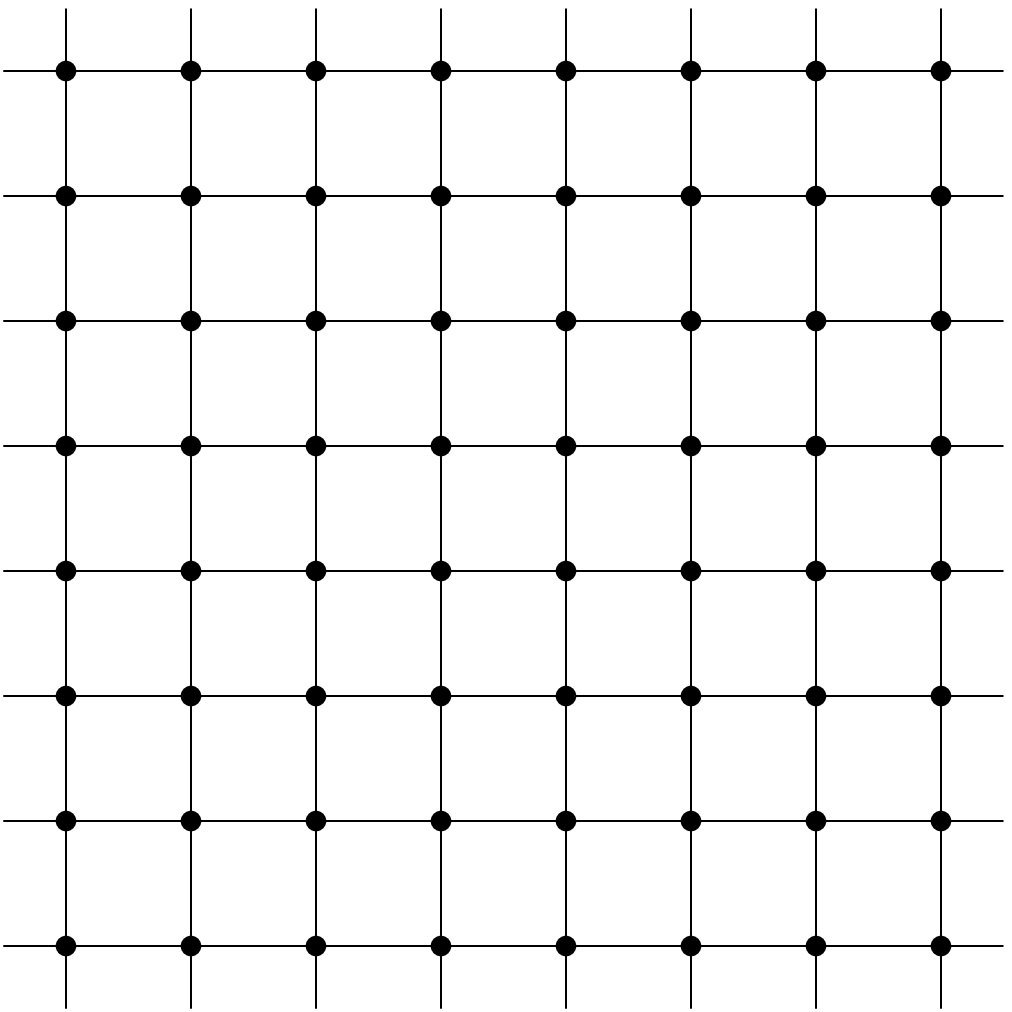,width=1.850in,silent=}%
\hspace{.25in}%
\psfig{figure=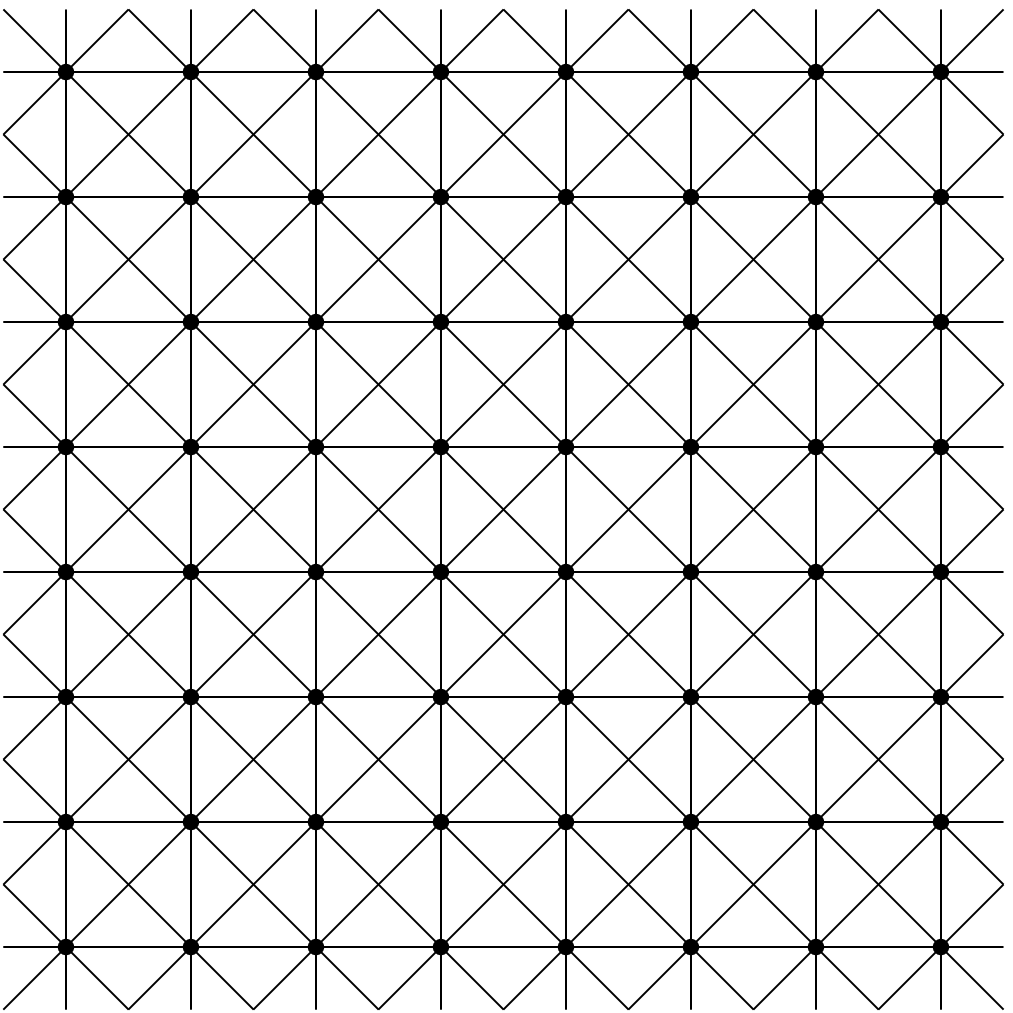,width=1.850in,silent=}%
\hspace{3.45in}%
\raisebox{-7pt}%
{\psfig{figure=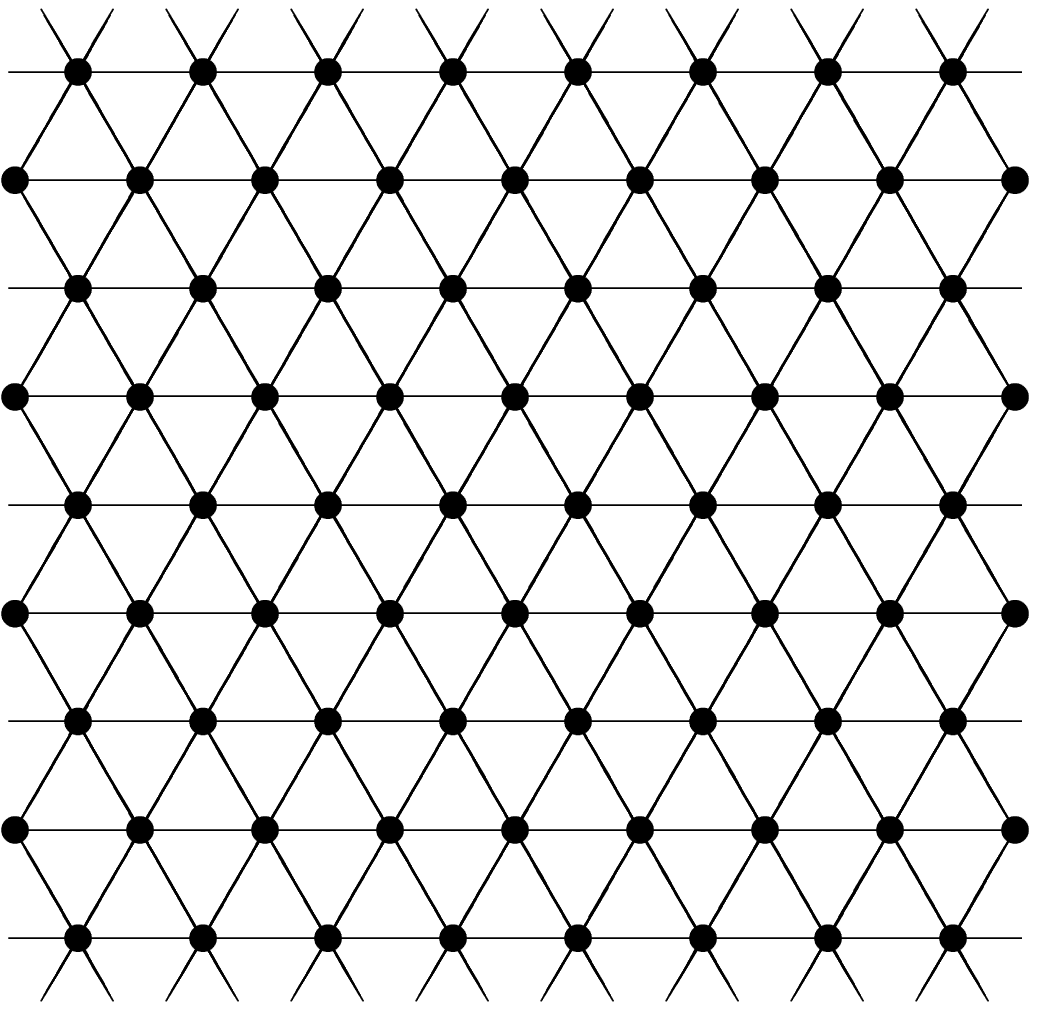,width=1.850in,angle=90,silent=}}%
\hspace*{-220pt}%
\vspace*{-.75ex}

\begin{center}
{\bf Figure\,1.}
{\sl The grid graph $\Gplus$, the infinity graph $\Gstar$, 
and the hexagonal graph $\Ghex$}\vspace{1.0ex}
\end{center}

\looseness=-1
that $d_3(z_1,z_2,z_3) \ge d$ for all $z_1,z_2,z_3 \in \cC$.
Numerous results on optimal tristance codes in~$\Z^2$ can 
be found in~\cite{EV,ME,SE2}. Optimal tristance \emph{anticodes} are the
subject of this paper.

\begin{definition}
\label{anticode-def}
A set $\cA_d \subset \Z^n$ is a {\dfn tristance anticode}
of diameter $d$ if $d_3(z_1,z_2,z_3) \le d$ for all $z_1,z_2,z_3 \in \cA_d$;
it is {\dfn optimal} if it has the largest possible
cardinality for its diameter $d$.\vspace*{-.60ex}
\end{definition}

Observe that \Dref{anticode-def} can be extended in the obvious way to 
other graphical models (such as the $\tstar$ and the $\thex$ models)
as well as to higher dispersions (such as quadristance).

\looseness=-1
One can also define tristance anticodes \emph{centered} about a given
point or a pair of points. Given $z_0 \,{\in}\, \Z^n$, a set 
$\cA_d(z_0) \subset \Z^n$ 
is said to be a tristance anticode of diameter~$d$ {\dfn centered about $z_0$}
if $d_3(z_0,z_1,z_2) \le d$ for all $z_1,z_2 \,{\in}\, \cA_d(z_0)$.
Given distinct $z_1,z_2 \in \Z^n$,~a~set 
$\cA_d(z_1,z_2) \subset \Z^n$ 
is said to be a tristance anticode of diameter $d$ {\dfn centered about 
$z_1$~and~$z_2$} if $d_3(z_1,z_2,z) \le d$ for all $z \,{\in}\,\cA_d(z_1,z_2)$.
Once again, we are interested in {\dfn optimal\/} centered tristance 
anticodes that have the largest possible cardinality for their diameter
and center(s).
We note that the corresponding problem for conventional (distance) 
anticodes is trivial. The unique optimal anticode of diameter $d$
centered about a given point $z_0 \in \Z^n$ is 
$
\sS_d(z_0) = \set{ z \in \Z^n : d(z,z_0) \le d}
$,
which is just a {\dfn sphere of radius $d$ about $z_0$}. 

\looseness=-1
The rest of this paper is organized as follows.
The next section is concerned with optimal tristance
anticodes in $\Z^2$. We determine the cardinality and classify the optimal
tristance anticodes $\cA_d(z_1,z_2)$, $\cA_d(z_0)$, and $\cA_d$ 
in the grid graph $\Gplus$, for all $d \ge 1$. We also introduce 
in \S2 certain methods and techniques that will be useful 
throughout this paper.
In~\S3,~we extend the results of \S2 to the $\tstar$
model and the $\thex$ model.
%
In~\S4,~we pursue generalizations\linebreak[3.99]
to higher dimensions and to
higher dispersions: we investigate optimal tristance
anticodes in $\Z^3$ and optimal quadristance anticodes in $\Z^2$.
In \S5, we discuss some of the applications of 
our results to multi-dimensional interleaving schemes with 
repetitions~\cite{BBV,BBF,EV,SE2} and to the study of connectivity loci
in the game of Go.  
\pagebreak[3.99]

\newpage
\vspace{3ex}
\section{Optimal tristance anticodes in the grid graph}
\vspace{-1.5ex} 
\label{sec2}

We will first need some auxiliary results. Trivially, the $L_1$-distance
between $(x_1,y_1)$ and $(x_2,y_2)$ in $\Z^2$ can
be written as 
$
(\max_{1 \leq i \leq 2} x_i - \min_{1 \leq i \leq 2} x_i)
 + 
(\max_{1 \leq i \leq 2} y_i - \min_{1 \leq i \leq 2} y_i)
$.
The~following theorem of~\cite{EV} shows that a similar
expression holds for tristance in $\Z^2$.

\begin{theorem}
\label{tristance-formula}
Let  $z_1 = (x_1,y_1)$, $z_2 = (x_2,y_2)$,
$z_3 = (x_3 ,y_3 )$ be distinct 
points in $\Z^2$. Then\vspace{-.75ex}
\Eq{
\label{d3}
d_3(z_1,z_2,z_3)
\ = \ 
\left(\max_{1 \leq i \leq 3} x_i - \min_{1 \leq i \leq 3} x_i \right)
\ + \ 
\left(\max_{1 \leq i \leq 3} y_i - \min_{1 \leq i \leq 3} y_i \right)
\vspace{1.5ex}
}
\end{theorem}

Next, we recall some results on optimal \underline{distance}
anticodes in the grid graph $\Gplus$. Let $d$ be even, let $z_0 = (x_0,y_0)$
be an arbitrary point in $\Z^2$, and consider the set
\Eq{
\label{L_1-sphere}
\sS_{d/2}(z_0)
\ = \ 
\set{ (x,y) \in \Z^2 ~:~ |x - x_0| + |y - y_0| \le \frac{d}{2}\, }
}
which is the $L_1$-sphere of radius $d/2$ about $z_0$.
By triangle inequality, for all $z_1,z_2\,{\in}\,\sS_{d/2}(z_0)$ we have 
$
d(z_1,z_2) \le d(z_1,z_0) + d(z_0,z_2) \le d
$, 
so $\sS_{d/2}(z_0)$ is an anticode of diameter~$d$.
It is easy to see that $\abs{\sS_{d/2}(z_0)} = d^2\!/2 + d + 1$.
On the other hand, it is shown in~\cite{BBV} that~$\Z^2$ can be
partitioned into $d^2\!/2 + d + 1$ codes 
such that the minimum $L_1$-distance of each code is $d+1$. Obviously,
any anticode of diameter $d$ can contain at most one point from 
each such code. 
It follows that the anticode $\sS_{d/2}(z_0)$ is optimal for all even $d$.
For odd $d$, we let $z_0$ be an arbitrary point in $(\shalf{1},0) + \Z^2$
or in $(0,\shalf{1}) + \Z^2$. Then $\sS_{d/2}(z_0)$ defined in~\eq{L_1-sphere}
is again an anticode of diameter $d$ (by triangle inequality), and 
$\abs{\sS_{d/2}(z_0)} = d^2\!/2 + d + \shalf{1}$.
Once again, it is shown in~\cite{BBV} that~$\Z^2$ can be
partitioned into $d^2\!/2 + d + \shalf{1}$ codes 
with even minimum distance $d+1$, so the anticode $\sS_{d/2}(z_0)$ 
is optimal for all odd $d$.

\subsection{Uniqueness of optimal anticodes in the grid graph} 
\vspace{-1ex}
\label{sec2.1}

\looseness=-1
We now show that optimal distance anticodes in $\Gplus$ 
are unique: if $\cA$ is an optimal~anti\-code of diameter $d$, then
$\cA = \sS_{d/2}(z_0)$, where $z_0\,{\in}\,\,\Z^2$ if $d$ is even, whereas 
if $d$ is odd then $z_0 \in (\shalf{1},0) + \Z^2$ or 
$z_0 \in (0,\shalf{1}) + \Z^2$. 
This result is established in a series of lemmas below. 

\looseness=-1
A set $\cS\,{\subseteq}\,\Z^2$ is said to be 
{\dfn vertically contiguous\,}
if it has the following property:~if~$(x,y_1)\,{\in}\,\cS$ 
and $(x,y_2)\,{\in}\,\cS$, then $(x,y)\,{\in}\,\cS$ for all $y$ in the 
range $\min\set{y_1,y_2} \le y \le \max\set{y_1,y_2}$.
Similarly, we say that $\cS \subseteq \Z^2$ is 
{\dfn horizontally contiguous\,} if 
the fact that $(x_1,y) \in \cS$ and $(x_2,y) \in \cS$ implies
that $(x,y) \in \cS$ for all $\min\set{x_1,x_2} \le x \le \max\set{x_1,x_2}$.

\begin{lemma}
\label{contiguous}
Let $\cA$ be an optimal anticode of diameter $d$ in the grid graph $\Gplus$.
Then $\cA$ is both vertically contiguous and horizontally contiguous.
\end{lemma}

\begin{proof}
\looseness=-1
Suppose that $z_1 = (x_0,y_1)$ and $z_2 = (x_0,y_2)$ are points in $\cA$.
Assume w.l.o.g.\ that~$y_2 > y_1$, and consider a point 
$z_3 = (x_0,y_3)$ with $y_1 \le y_3 \le y_2$. If $z = (x,y) \in \cA$,
then
$$
d(z,z_3)
\: = \:
|x - x_0| \,+\, |y - y_3|
\ \le \
|x - x_0| \,+\, \max\set{|y - y_1|,|y - y_2|}
\: = \: 
\max\set{d(z,z_1),d(z,z_2)}
$$
Thus $d(z,z_3) \le d$ for all $z \in \cA$, 
and if $\cA$ is optimal, it must 
contain~the~point $z_3 = (x_0,y_3)$.
Hence $\cA$ is vertically contiguous.
By a similar argument, $\cA$ is horizontally 
contiguous.~\qed\pagebreak[3.99]
\end{proof}

Given $z = (x,y) \in \Z^2$, we say that the points $(x-1,y)$ and 
$(x+1,y)$ are the {\dfn horizontal neighbors\/} of $z$, while the 
points $(x,y-1)$ and $(x,y+1)$ are the {\dfn vertical neighbors\/} of $z$.

\begin{lemma}
\label{neighbors}
Let $\cA$ be an optimal anticode of diameter $d$ in the grid graph $\Gplus$.
If~$\cA$~contains the two horizontal neighbors of a point $z \in \Z^2$,
or if $\cA$ contains the two vertical neighbors of $z$, then $\cA$ 
necessarily contains $z$ itself and \underline{all the four} neighbors of~$z$.
\end{lemma}

\begin{proof} 
Suppose $\cA$ contains the points 
$z_1 = (x_0-1,y_0)$ and $z_2 = (x_0+1,y_0)$.
Since $\cA$ is horizontally contiguous by \Lref{contiguous},
it also contains the point $z_0 = (x_0,y_0)$. Moreover, 
if $z = (x,y)$ 
is any point in $\cA$, 
then 
$
d(z,z_0) = \max\set{d(z,z_1),d(z,z_2)} - 1 \le d-1
$. 
Now, let $z_3 = (x_0,y_0+1)$. Then we have
$$
d(z,z_3)
\ = \ 
|x - x_0| + |y - y_0 - 1|
\ \le \
|x - x_0| + |y - y_0| + 1
\ = \
d(z,z_0) + 1
\ \le \
d
$$
Hence, if $\cA$ is optimal, it must contain $z_3 = (x_0,y_0+1)$.
By a similar argument, $\cA$ also contains the point $(x_0,y_0-1)$.
The claim for vertical neighbors follows by symmetry.\,~\qed\vspace{1.25ex}
\end{proof}

Given a set $\cS \subseteq \Z^2$, 
let $\Gplus(\cS)$ denote the induced subgraph of $\Gplus$,
consisting of $\cS$ and the edges of $\Gplus$ with both 
endpoints in $\cS$. 
We say that $z \in \cS$ is an {\dfn internal point\/} of $\cS$
if $z$ has degree $4$ in $\Gplus(\cS)$; otherwise we say that
$z$ is a {\dfn boundary point\/} of $\cS$.

\begin{lemma}
\label{boundary}
An optimal anticode of diameter $d$ in $\Gplus$ 
has at most $2d$ boundary points. 
\end{lemma}

\begin{proof}
Let $\cA$ be an optimal anticode of diameter $d$ in $\Gplus$,
and define the integers $\xmin$, $\xmax$ as follows:
$\xmin = \min\set{ x : (x,y) \in \cA}$
and
$\xmax = \max\set{ x : (x,y) \in \cA}$.
Clearly $\Delta = \xmax - \xmin \le d$.
Let us partition $\cA$ into $\Delta+1$ vertical segments
\Eq{
\label{V_i}
\cV_i 
\ \ \deff \ \
\Bigl\{ (x,y) \in \cA ~:~ x = \xmin + i\, \Bigr\}
\hspace{6ex}
\text{for $i = 0,1,\ldots,\Delta$}
}
Since $\cA$ is vertically contiguous by \Lref{contiguous}, 
for each $i = 0,1,\ldots,\Delta$, the vertical segment $\cV_i$ 
in~\eq{V_i} can be written as
\Eq{
\label{V_i2}
\cV_i 
\ = \
\Bigl\{ 
(\xmin{+}\,i,y_{\min,i}),
(\xmin{+}\,i,y_{\min,i}+1),
\ldots,
(\xmin{+}\,i,y_{\max,i})
\Bigr\}
}
for some integers $y_{\min,i} \le y_{\max,i}$.
Notice that for each point $z \in \cV_i$, except 
$(\xmin{+}\,i,y_{\min,i})$ and $(\xmin{+}\,i,y_{\max,i})$,
both vertical neighbors of $z$ are in $\cV_i$, and hence
also in $\cA$. \Lref{neighbors} thus implies that $z$
has degree $4$ in $\Gplus(\cA)$. It follows that each $\cV_i$
contains at most two boundary points of $\cA$, so that $\cA$ has
at most $2(\Delta+1)$ boundary points altogether. If $\Delta \le d-1$
we are done, so it remains to consider the case $\Delta = d$.
But then $|\cV_0| = |\cV_{\Delta}| = 1$
\looseness=-1
(if $y_{\min,0} = y_{\max,0} = y_{\min,\Delta} = y_{\max,\Delta}$
does not hold, there are points in $\cV_0 \cup \cV_{\Delta}$
at distance $\ge d+1$ from each other).
Thus $\cA$ has at most $2d$ boundary points in this 
case~as~well.~\qed\vspace{1.25ex}
\end{proof}

\begin{lemma}
\label{internal}
Let $\cA$ be an anticode of diameter $d \ge 2$ in the grid graph $\Gplus$.
Then the set of internal points of $\cA$, if nonempty, forms an anticode 
of diameter $d-2$.
\end{lemma}

\begin{proof}
\looseness=-1
By convention, a set of size $\le 1$ has diameter zero. 
Otherwise, if $z_1,z_2$~are~distinct internal
points of $\cA$, then all of their neighbors
are also in $\cA$.
Observe that the set of $4$ neighbors of $z_1$ always contains
at least one point $z$ such that $d(z,z_2) = d(z_1,z_2) +1$.
It follows that among the 
neighbors of $z_1$ and $z_2$, there are (at least) two points at distance
$d(z_1,z_2) + 2$ from each other.
Hence if $\cA$ has diameter $d$, then 
$d(z_1,z_2) \le d-2$.\,~\qed\pagebreak[3.99]
\end{proof}

\begin{theorem}
\label{distance-uniqueness}
Let $\cA$ be an optimal anticode of diameter $d$ in the grid graph $\Gplus$.
Then 
$\cA = \sS_{d/2}(z_0)$, where $z_0\,{\in}\,\,\Z^2$ if $d$ is even
and $z_0 \in \{ (\shalf{1},0) + \Z^2 \} \cup \{ (0,\shalf{1}) + \Z^2\}$
otherwise.
\end{theorem}

\begin{proof}
We will only prove the theorem for even $d$; the proof for odd $d$
is similar. We proceed by induction on $d$. For $d=2$, it can be 
readily verified that an anticode~of~diameter $2$ and size $5$ is
necessarily the $L_1$-sphere $\sS_1(z_0)$ for some $z_0 \in \Z^2$.
Now,~let~$\cA$~be an anticode of diameter $d$ and cardinality
$\abs{\cA} = d^2\!/2 + d + 1$. Let $\sD(\cA)$ denote the set 
of internal points of $\cA$. Then $\sD(\cA)$ is an anticode 
of diameter $d-2$ by \Lref{internal}, and
$$
\abs{\sD(A)} 
\ \ge \
\abs{\cA} - 2d
\ = \
\frac{(d\,{-}\,2)^2}{2} \, + \, (d\,{-}\,2) \, + \, 1
$$
by \Lref{boundary}. It follows that $\cA$ has \emph{exactly} $2d$
boundary points and $\sD(\cA)$ is an \emph{optimal} anticode of 
diameter $d-2$. Hence, by induction hypothesis, 
$\sD(\cA) = \sS_{(d/2)-1}(z_0)$~for~some $z_0 \in \Z^2$.
Referring to~\eq{V_i} and~\eq{V_i2}, we see that each 
vertical segment $\cV_1,\cV_2,\ldots,\cV_{\Delta-1}$
in $\cA$ must have exactly two boundary points of $\cA$, and
$|\cV_0| = |\cV_{\Delta}| = 1$.
In other words, there is a unique way to adjoin $2d$ boundary 
points to $\sD(\cA) = \sS_{(d/2)-1}(z_0)$ to obtain an optimal
anticode, and it is easy to 
\looseness=-1
see that the result is precisely the $L_1$-sphere 
$\sS_{d/2}(z_0)$.~~\qed\vspace{.75ex}
\end{proof}

The foregoing theorem, which is the main result of this subsection,
establishes the uniqueness of optimal distance anticodes in $\Gplus$.
All such anticodes are $L_1$-spheres of radius $d/2$.
We are now in a position to begin the classification of optimal
\underline{tristance} anticodes in $\Gplus$.

\subsection{Centered tristance anticodes in the grid graph} 
\vspace{-1ex}
\label{sec2.2}

Recall that $\cA_d(z_0) \subset \Z^2$ 
is a tristance anticode of diameter~$d$ 
{\em centered about $z_0 \in \Z^2$},
if $d_3(z_0,z_1,z_2) \le d$ for all $z_1,z_2 \,{\in}\, \cA_d(z_0)$.
First assume that $d$ is even, and consider~the~$L_1$-sphere 
$\sS_{d/2}(z_0)$. For all $z_1,z_2 \in \sS_{d/2}(z_0)$, we have
$
d_3(z_0,z_1,z_2) 
\le 
d(z_0,z_1)  +  d(z_0,z_2)
\le 
d
$.
It follows that $\sS_{d/2}(z_0)$ is a centered tristance anticode
of diameter $d = 2t$ and cardinality $2t^2 + 2t + 1$. 
Now suppose that $d = 2t+1$ is odd, and consider the $L_1$-sphere 
$\sS_{d/2}(z_0+\xi)$, where $\xi = (\shalf{1},0)$. Let $z_0 = (x_0,y_0)$,
and let $z_1 = (x_1,y_1)$ be any point in $\sS_{d/2}(z_0+\xi)$.
Then $d(z_0,z_1) \le t+1$ and, moreover, $d(z_0,z_1) \le t$ unless $x_1 > x_0$.
It follows that $\sS_{d/2}(z_0+\xi)$ is a tristance anticode 
centered about $z_0$, of diameter $d=2t+1$ and cardinality 
$2t^2 + 4t + 2$. 
The following theorem shows that the anticodes
constructed above are the {\sl unique\/} optimal 
tristance anticodes centered about $z_0$ in the 
grid graph $\Gplus$.

\begin{theorem}
\label{one-center}
Let $\cA_d(z_0)$ be an optimal tristance anticode of diameter $d$ 
in $\Gplus$ centered about $z_0 \in \Z^2$. 
If $d$ is even, then 
$\cA_d(z_0) = \sS_{d/2}(z_0)$. If $d$ is odd, then 
$\cA_d(z_0) = \sS_{d/2}(z_0+\xi)$ for some 
$\xi \in \{(\shalf{1},0), (0,\shalf{1}), (-\shalf{1},0), (0,-\shalf{1})\}$.
\vspace{-0.5ex}
\end{theorem}

\begin{proof}
We will prove the theorem for even $d$ only; the proof for odd $d$
is similar. 
For all $z_1,z_2 \in \cA_d(z_0)$, we have 
$d(z_1,z_2) \le d_3(z_0,z_1,z_2) \le d$
by definition. Hence $\cA_d(z_0)$~is~also
a~distance anticode of diameter $d$ and
\Eq{
\label{centered-aux}
|\cA_d(z_0)|
\ \le \
d^2\!/2 + d + 1
}
We have shown that $\sS_{d/2}(z_0)$ is a tristance anticode
of diameter $d$ centered about $z_0$.~~Since 
$|\sS_{d/2}(z_0)| = d^2\!/2 + d + 1$,
this anticode is optimal in view of~\eq{centered-aux}.
Moreover, by \Tref{distance-uniqueness}, equality 
in~\eq{centered-aux} is possible only if 
$\cA_d(z_0) = \sS_{d/2}(z_0)$.~~\qed\vspace*{-1.00ex}\pagebreak[3.99]
\end{proof}

\looseness=-1
Next, we consider the anticodes $\cA_d(z_1,z_2) \subset \Z^2$ centered
about a pair of points $z_1$ and $z_2$ and defined by the property
that $d_3(z_1,z_2,z) \le d$ for all $z\,{\in}\,\cA_d(z_1,z_2)$. Such
anticodes~can\-not have an arbitrary diameter: if $d(z_1,z_2) = \Delta$,
then $\cA_d(z_1,z_2) = \varnothing$ unless $d \ge \Delta$.
For $d=\Delta$, it turns out that the optimal anticode 
$\cA_{\Delta}(z_1,z_2)$ is the bounding rectangle~of~$z_1,\,z_2$.

\begin{definition}
\label{rectangle-def}
\looseness=-1
Let $z_1 = (x_1,y_1), z_2 = (x_2,y_2), \ldots, z_n = (x_n,y_n)$ 
be $n$ distinct points~in~$\Z^2$. The {\dfn bounding rectangle}
of $z_1,z_2,\ldots,z_n$ is the smallest rectangle 
$\cR(z_1,z_2,\ldots,z_n)$ with edges parallel to the axes that 
contains all the $n$ points. Explicitly, let
$
x_{\max} = \max\{x_1,\ldots,x_n\}
$, 
$
x_{\min} = \min\{x_1,\ldots,x_n\}
$,
$
y_{\max} = \max\{y_1,\ldots,y_n\}
$, 
and
$
y_{\min} = \min\{y_1,\ldots,y_n\} 
$.
Then
\Eq{
\label{rectangle}
\cR(z_1,z_2,\ldots,z_n)
\ \ \deff \ \
\left\{ \Strut{2.5ex}{0ex}
(x,y) \in \Z^2 ~:~ 
x_{\min} \le x \le x_{\max} \ \text{ and }\ y_{\min} \le y \le y_{\max} \,
\right\} 
}
\end{definition}

\looseness=-1
By \Tref{tristance-formula}, 
the tristance of {\sl any\/} three points in $\cR(z_1,z_2)$
is at most $d(z_1,z_2)$. 
Moreover, if $z \not\in \cR(z_1,z_2)$, then 
$d_3(z_1,z_2,z) > d(z_1,z_2)$. This implies
that if $d = d(z_1,z_2)$, then $\cR(z_1,z_2)$ is
the optimal anticode $\cA_d(z_1,z_2)$.
For general $d$, we have the following theorem.

\begin{theorem}
\label{centeroct}
Let $z_1 = (x_1,y_1)$, $z_2 = (x_2,y_2)$ be distinct 
points in $\Z^2$ and assume,~w.l.o.g., 
that $x_2 \ge x_1$ and\, $y_2 \ge y_1$ 
so that $d(z_1,z_2) = (x_2 - x_1) + (y_2 - y_1)$.
Let $\cA_d(z_1,z_2)$~be~the optimal tristance anticode in $\Gplus$
of diameter $d \ge d(z_1,z_2)$ centered about $z_1$ and $z_2$. 
Write $c = d - d(z_1,z_2)$.
Then $\cA_d(z_1,z_2)$ consists of all $z = (x,y)$ in $\Z^2$ 
such that \vspace{-.50ex}
\begin{eqnarray}
\label{centeroct1}
x_1 - c \ \le\ x \ \le\ x_2 + c\, ,
&\hspace{1ex}& 
x_1 + y_1 - c \ \le\ x+y \ \le\ x_2 + y_2 + c
\\[1.0ex]
\label{centeroct2}
y_1 - c \ \le\ y \ \le\ y_2 + c\, ,
&\hspace{1ex}& 
x_1 - y_2 - c \ \le\ x-y \ \le\ x_2 - y_1 + c
\end{eqnarray}
\end{theorem}

\begin{proof}
The given points $z_1,z_2$ completely determine all the other points 
in $\cA_d(z_1,z_2)$, as follows:
$
\cA_d(z_1,z_2)
 = 
\left\{ z \in \Z^2 : d_3(z_1,z_2,z) \le d \right\}
$.
It is now easy to see that 
$z \in \cA_d(z_1,z_2)$
if and only if the $L_1$-distance from $z$ to (the closest point of)
the {bounding rectangle}~$\cR(z_1,z_2)$ is at most $c = d - d(z_1,z_2)$.
This is precisely the property expressed 
by \eq{centeroct1} and \eq{centeroct2}.~~\qed\vspace{1.00ex}
\end{proof}

\subsection{General tristance anticodes in the grid graph} 
\vspace{-1ex}
\label{sec2.3}

We will use the results of \S\ref{sec2.2}, especially \Tref{centeroct},
to classify unrestricted (non-centered) optimal tristance anticodes in 
$\Gplus$. The subset of $\Z^2$ defined by 
equations \eq{centeroct1}, \eq{centeroct2}~in~\Tref{centeroct}
is an example of a set we call an {\sl octagon}. 
We formalize this as follows.

\begin{definition}
Let $\alpha_1,\alpha_2,\ldots,\alpha_8$ be arbitrary real
constants.
An {\dfn octagon $\sO(\alpha_1,\alpha_2,\ldots,\alpha_8)$} 
is a subset of\, $\Z^2$ defined by the inequalities
\begin{eqnarray}
\label{oct-def1}
\alpha_1 \le x \le \alpha_5\, ,
&\hspace{4ex}& 
\alpha_3 \le x+y \le \alpha_7\\[0.750ex]
\label{oct-def2}
\alpha_2 \le y \le \alpha_6\, ,
&& 
\alpha_4 \le x-y \le \alpha_8
\end{eqnarray}
\end{definition}

A generic octagon $\sO(\alpha_1,\alpha_2,\ldots,\alpha_8)$ is 
illustrated in Figure\,2.
Note that $\sO(\alpha_1,\alpha_2,\ldots,\alpha_8)$
may have fewer than eight sides (say, if $\al_3 \le \al_1+\al_2$),
or may be empty altogether. 
Octagons will
play an important role in this paper. Observe that the 
$L_1$-spheres $\sS_{d/2}(z_0)$ and $\sS_{d/2}(z_0+\xi)$ in 
\Tref{one-center} are octagons. 
By~\Tref{centeroct}, the optimal anticode 
$\cA_d(z_1,z_2)$ is also an octagon. The~following lemma establishes 
another useful property of octagons.\pagebreak[3.99]

\centerline{\psfig{figure=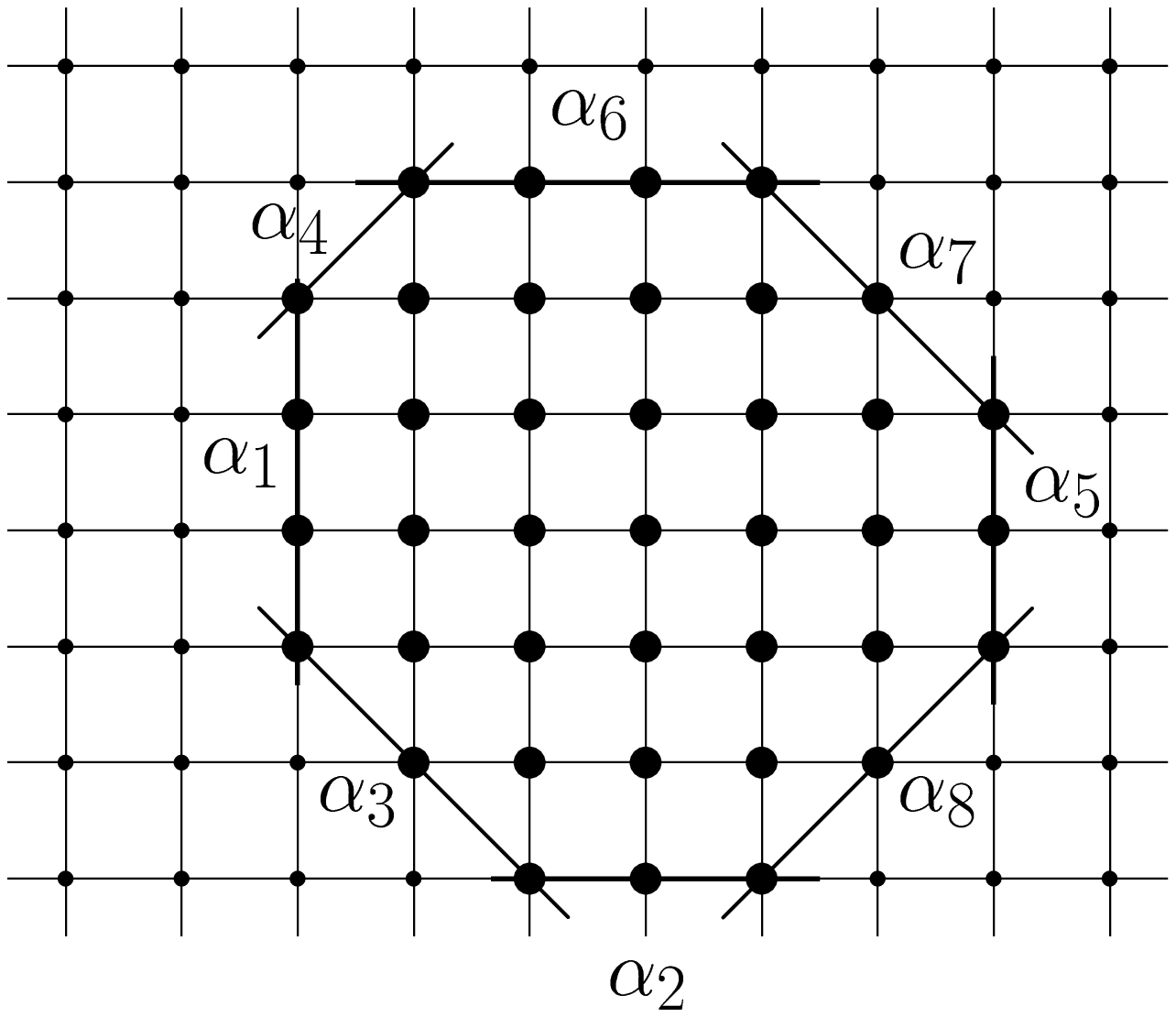,width=2.25in,silent=}}
\vspace{-2.5ex}

\begin{center}
{\bf Figure\,2.}
{\sl A generic octagon $\sO(\alpha_1,\alpha_2,\ldots,\alpha_8)$
in the grid graph $\Gplus$}
\vspace{1.5ex}
\end{center}

\begin{lemma}
\label{octint}
The intersection of any two octagons is an octagon.
\vspace{-.5ex}
\end{lemma}

\begin{proof}
It is obvious that 
$
\sO(\alpha_1,\alpha_2,\ldots,\alpha_8)
\,\cap\,
\sO(\beta_1,\beta_2,\ldots,\beta_8)
 = 
\sO(\gamma_1,\gamma_2,\ldots,\gamma_8)
$,~where 
$\gamma_i = \max\{\alpha_i,\beta_i\}$ for $i = 1,2,3,4$
and $\gamma_i = \min\{\alpha_i,\beta_i\}$ for 
$i = 5,6,7,8$.~~\qed\vspace{1.00ex}
\end{proof}

\looseness=-1
We are now in a position to begin the classification of unrestricted
optimal tristance~anticodes in $\Gplus$. The next lemma establishes
a certain closure property of such anticodes.

\begin{lemma}
\label{intersection}
Let $\cA_d$ be an optimal tristance anticode of diameter $d$ in $\Gplus$. 
Then $\cA_d$ is closed under intersection
with anticodes centered about pairs of its own points, namely
\Eq{
\cA_d
\ = \!\! 
\bigcap_{z_1,z_2 \in \cA_d} \!\!\! \cA_d(z_1,z_2)
\vspace{.75ex}
}
\end{lemma}

\begin{proof}
If $z\in \cA_d$ then, by the definition of an anticode,
we have $d_3(z_1,z_2,z)\leq d$ for all $z_1,z_2\in \cA_d$.
Thus $z\in \cA_d(z_1,z_2)$\, for all $z_1,z_2\in \cA_d$, and therefore
$z\in \bigcap_{z_1,z_2\in \cA_d}\cA_d(z_1,z_2)$. Hence, for 
{\sl any\/} tristance anticode $\cA$ of diameter $d$, we have
$$
\cA
\ \subseteq \!\! 
\bigcap_{z_1,z_2 \in \cA} \!\! \cA_d(z_1,z_2)
\vspace{-1ex}
$$ 
Now, if $\cA_d$ is optimal and $z\not\in \cA_d$, then 
there exist $z_1,z_2 \in \cA_d$ such that $d_3(z_1,z_2,z) > d$;
otherwise, we could adjoin $z$ to $\cA_d$ to obtain a larger anticode. 
For these $z_1,z_2 \in \cA_d$, we have
$z\not\in \cA_d(z_1,z_2)$. Hence 
$z\not\in \bigcap_{z_1,z_2\in \cA_d}\cA_d(z_1,z_2)$,
and the lemma follows.~~\qed\vspace{1.00ex}
\end{proof}

Combining \Tref{centeroct}, \Lref{octint}, and \Lref{intersection}
makes it possible to determine the {\sl shape}
of optimal tristance anticodes in the grid graph.

\begin{lemma}
\label{shape}
Let $\cA_d$ be an optimal tristance anticode of diameter $d$ in $\Gplus$. 
Then $\cA_d$
is an octagon $\sO(\alpha_1,\alpha_2,\ldots,\alpha_8)$
for some $\alpha_1,\alpha_2,\ldots,\alpha_8 \in \Z$.
\vspace{-.5ex}
\end{lemma}

\begin{proof}
By \Lref{intersection}, we have
$\cA_d = \bigcap_{z_1,z_2\in \cA_d}\cA_d(z_1,z_2)$.
By \Tref{centeroct}, each of the sets $\cA_d(z_1,z_2)$ is
an octagon. By \Lref{octint}, an intersection of octagons
is also an octagon.~\qed\vspace{1.00ex}
\end{proof}

Using translations in $\Z^2$, we may always assume 
w.l.o.g. that $\al_1 = \al_2 = 0$ in~\eq{oct-def1} and~\eq{oct-def2}.
Thus, in view of \Lref{shape}, we have
$
\cA_d
 = 
\sO(0,0,\alpha_3,\alpha_4,\ldots,\alpha_8)
$,
and it remains~to~determine the six integer parameters
$\alpha_3,\alpha_4,\alpha_5,\alpha_6,\alpha_7,\alpha_8$
as a~function of the diameter~$d$.

To this end, we first rewrite the definition of an 
octagon  $\sO(0,0,\alpha_3,\alpha_4,\ldots,\alpha_8)$
in a~different form.
This octagon can be defined 
as the set of all $(x,y) \in \Z^2$ such that
\begin{eqnarray}
\label{new-oct1}
0\kern-1pt \ \le\ x \ \le\ \kern-1pt a\, ,
&\hspace{1ex}&  
c_0 \ \le\ x+y \ \le\  a + b - c_2\\[0.750ex]
\label{new-oct2}
0 \ \le\ y \ \le\ b\, ,
&& 
c_3 - b \ \le\ x-y \ \le \ a - c_1
\end{eqnarray}
where 
$a = \alpha_5$, 
$b = \alpha_6$, 
$c_0 \,{=}\, \alpha_3$,
$c_1 \,{=}\, \alpha_5 - \alpha_8$, 
$c_2 \,{=}\, \alpha_5 + \alpha_6 - \alpha_7$, 
and
$c_3 \,{=}\, \alpha_4 + \alpha_6$
(cf.\,Figure\,3).
We omit the tedious, but easy, proof of the transformation 
from~\eq{oct-def1}--\eq{oct-def2}
to~\eq{new-oct1}--\eq{new-oct2}.\vspace{2.5ex}

\centerline{\psfig{figure=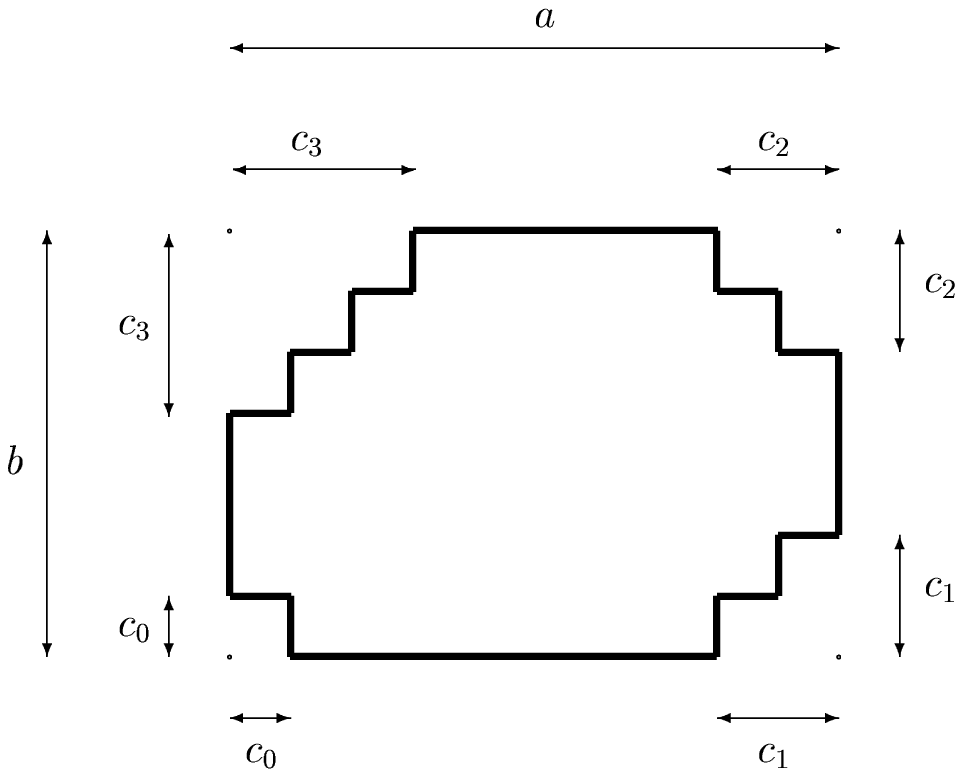,width=3.35in,silent=}}
\vspace{-2.5ex}

\begin{center}
{\bf Figure\,3.}
{\sl An octagon $\sO(0,0,\alpha_3,\alpha_4,\ldots,\alpha_8)$
defined in terms of $a$, $b$, $c_0$, $c_1$, $c_2$, and $c_3$}\vspace{.5ex}
\end{center}

We will use $\sO(a,b,c_0,c_1,c_2,c_3)$ to denote an 
octagon $\sO(0,0,\alpha_3,\alpha_4,\ldots,\alpha_8)$
specified in terms of the alternative parameters 
$a$, $b$, $c_1$, $c_2$, $c_3$, $c_4$ of~\eq{new-oct1}
and \eq{new-oct2}. It is easy to see from Figure\,3 that 
the size of $\sO(a,b,c_0,c_1,c_2,c_3)$ is given by
\Eq{
\label{A_d-size}
\left|\sO(a,b,c_0,c_1,c_2,c_3)\right|
\ = \
(a+1)(b+1) \ - \ \sum_{i=0}^3 \frac{c_i(c_i+1)}{2}
}
The next step is to determine the maximum tristance $d$ of an optimal
tristance anticode $\cA_d = \sO(a,b,c_0,c_1,c_2,c_3)$ as a function 
of its parameters. 
We say that points $z_1,z_2,z_3 \in \cA_d$ are 
{\dfn diametric\/} if they attain the maximum tristance in $\cA_d$,
that is, if $d_3(z_1,z_2,z_3) = d$.

\begin{lemma}
\label{diametric}
An optimal tristance anticode 
$\cA_d = \sO(a,b,c_0,c_1,c_2,c_3)$ 
always contains~diametric points $z^*_1,z^*_2,z^*_3$
such that $z^*_1=(x^*_1,0)$ and $z^*_3=(x^*_3,b)$ 
for some $x^*_1,x^*_3 \in \Z$.\vspace{-.50ex}
\end{lemma}

\begin{proof}
Let $z_1 = (x_1,y_1)$, $z_2 = (x_2,y_2)$, and $z_3 = (x_3,y_3)$
be three diametric points~in $\cA_d$. W.l.o.g.\
assume that $\min\{y_1,y_2,y_3\} = y_1$ and $\max\{y_1,y_2,y_3\} = y_3$.
Since \mbox{$d_3(z_1,z_2,z_3) = d$}, \Tref{tristance-formula} implies
that the point $z_1 - (0,1) = (x_1,y_1{-}1)$ is not in $\cA_d$.
Referring~to~\mbox{\eq{oct-def1}--\eq{new-oct2}} along with Figure\,2 
and Figure\,3,
this in turn implies that either $y_1 = 0$ or $x_1 + y_1 = c_0$
or $x_1 - y_1 = a - c_1$. If $y_1 = 0$, we take $z^*_1 = z_1$.
Otherwise, if $x_1 + y_1 = c_0$, let
$$
z'_1 \ \ \deff \ \ z_1 \ - \ (x_1-c_0,y_1) \ = \ (c_0,0)
$$
(the point $z'_1$ is obtained from $z_1$ by moving down along the
South-West edge of the~octagon until reaching the South edge).
Clearly $z'_1 \,{\in}\, \cA_d$.
Moreover \mbox{$d_3(z'_1,z_2,z_3) \ge d_3(z_1,z_2,z_3) = d$},\pagebreak[3.99]
since replacing $y_1 = \min\{y_1,y_2,y_3\}$ by $0$ increases the 
tristance by $y_1$ (cf.\,\Tref{tristance-formula}) while replacing
$x_1$ by $c_0 = x_1 + y_1$ decreases the tristance by at most $y_1$.
Hence the points $z'_1,z_2,z_3$ are diametric, and we take $z^*_1 = z'_1$.
If $x_1 - y_1 = a - c_1$, we replace $z_1$ by
$
z'_1 
= (a{-}c_1,0)
$.
Again, it is easy to see that $z'_1,z_2,z_3$ are diametric, 
and we take $z^*_1 = z'_1$.
We can now proceed with the diametric points 
$z^*_1,z_2,z_3$ in a similar fashion to obtain 
$z^*_3$.~~\qed\vspace{6.00ex}
\end{proof}

\centerline{%
\psfig{figure=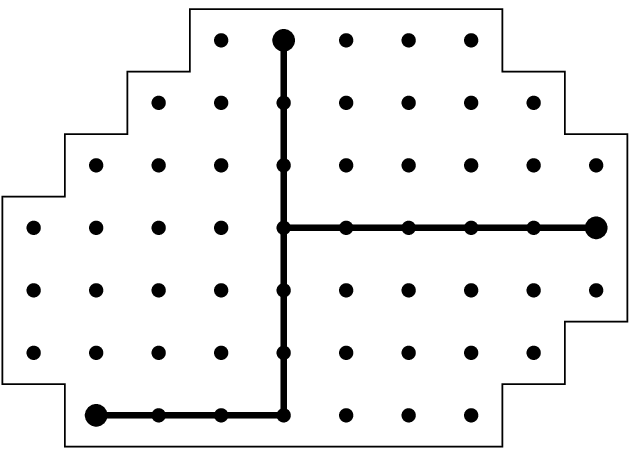,width=2.40in,silent=}}
\vspace{-.5ex}

\begin{center}
{\bf Figure\,4.}
{\sl A diametric configuration in $\cA_d$}
\vspace{1.5ex}
\end{center}

\begin{corollary}
\label{diameter}
The diameter of an optimal tristance 
anticode $\cA_d = \sO(a,b,c_0,c_1,c_2,c_3)$~is 
given by\vspace{-.50ex}
\Eq{
\label{Cor13}
d \ = \  a + b - \min\{c_0,c_1,c_2,c_3\}
}
\end{corollary}

\begin{proof}
W.l.o.g.\ assume that $\min\{c_0,c_1,c_2,c_3\} = c_0$.
Let $z_1 = (x_1,0)$, $z_2 = (x_2,y_2)$, and $z_3 = (x_3,b)$
be the diametric triple exhibited in \Lref{diametric}.
Then, in view of \Tref{tristance-formula}, 
\Eq{
d 
\ = \
d_3(z_1,z_2,z_3) 
\ = \
b \ + \ \max\{x_1,x_2,x_3\} \ - \ \min\{x_1,x_2,x_3\} 
}
It is easy to see that
$
\max\{x_1,x_2,x_3\} - \min\{x_1,x_2,x_3\} \le a - \min\{c_0,c_1,c_2,c_3\}
= a - c_0
$.
Equality in this bound is achieved for $x_2 = a$ and $x_1 = c_0$,
as illustrated in Figure\,4.~~\qed\vspace{.75ex}
\end{proof}

\begin{corollary}
\label{c-equality}
If $\cA_d = \sO(a,b,c_0,c_1,c_2,c_3)$ 
is an optimal tristance anticode, then
\Eq{
\label{c-equal}
c_0 \ = \ c_1 \ = \ c_2 \ = \ c_3
}
\end{corollary}

\begin{proof}
Obviously, \eq{c-equal}
maximizes the size of $\cA_d = \sO(a,b,c_0,c_1,c_2,c_3)$ 
in~\eq{A_d-size} for the given diameter
$
d  =   a + b - \min\{c_0,c_1,c_2,c_3\}
$.~~\qed\vspace{.75ex}
\end{proof}

\begin{theorem}
\label{optimal-tristance}
Let $\cA_d$ be an optimal tristance anticode of diameter~$d$
in the grid graph~$\Gplus$. Then
$$
|\cA_d|
\ = \
\ceilenv{\frac{2d^2 + 6d + 4}{7}}
\ = \
\ceilenv{\frac{2(d+1)(d+2)}{7}}
\vspace{.5ex}
$$
Moreover, up to rotation by an angle of $\pi/2$ and translation, 
$\cA_d = \sO(a,b,c,c,c,c)$
where the parameters 
$a$, $b$, and $c$ are given as a function of $d$ in Table\,1.
\end{theorem}

\begin{proof}
It follows from \Lref{shape} in conjunction with \Cref{c-equality} that 
$\cA_d$ is an octagon of the form $\sO(a,b,c,c,c,c)$,
for some $a,b,c \in \Z$. The size of $\cA_d$ is\pagebreak[3.99]
$(a+1)(b+1) - 2c(c+1)$

$$
\begin{array}{|c|c|c|c|c|}
\hline \hline
& & & & \\[-1.50ex]
\,d \hspace{-1.50ex}\pmod{7} & a & b & c & |\cA_d| 
\\[0.75ex]
\hline\hline
& & & & \\[-1.50ex]
0  & \frac{4d}{7} & \frac{4d}{7} & \frac{d}{7} & \frac{2d^2+6d+7}{7} 
\\[0.55ex]
\hline
& & & & \\[-1.95ex]
1  & \frac{4d+3}{7} & \frac{4d-4}{7} & \frac{d-1}{7} & \frac{2d^2+6d+6}{7} 
\\[0.55ex]
\hline
& & & & \\[-1.95ex]
2  & \frac{4d-1}{7} & \frac{4d-1}{7} & \frac{d-2}{7} & \frac{2d^2+6d+8}{7} 
\\[0.55ex]
\hline
& & & & \\[-1.95ex]
3  & \frac{4d+2}{7} & \frac{4d-5}{7} & \frac{d-3}{7} & \frac{2d^2+6d+6}{7} 
\\[0.55ex]
\hline
& & & & \\[-1.95ex]
4  & \frac{4d-2}{7} & \frac{4d-2}{7} & \frac{d-4}{7} & \frac{2d^2+6d+7}{7} 
\\[0.55ex]
\hline
& & & & \\[-1.50ex]
5  & 
\begin{array}{c}
\frac{4d+1}{7}\\[.75ex]
\frac{4d+1}{7}
\end{array}
& 
\begin{array}{c}
\frac{4d+1}{7}\\[.75ex]
\frac{4d-6}{7}
\end{array}
&
\begin{array}{c}
\frac{d+2}{7}\\[.75ex]
\frac{d-5}{7}
\end{array}
& 
\frac{2d^2+6d+4}{7}
\\[2.25ex]
\hline
& & & & \\[-1.50ex]
6  & 
\begin{array}{c}
\frac{4d-3}{7}\\[.75ex]
\frac{4d+4}{7}
\end{array}
& 
\begin{array}{c}
\frac{4d-3}{7}\\[.75ex]
\frac{4d-3}{7}
\end{array}
&
\begin{array}{c}
\frac{d-6}{7}\\[.75ex]
\frac{d+1}{7}
\end{array}
& 
\frac{2d^2+6d+4}{7}
\\[3.15ex]
\hline\hline
\end{array}
$$

\begin{center}
{\bf Table\,1.}
{\sl Parameters of optimal tristance anticodes in the grid graph $\Gplus$}
\vspace{.5ex}
\end{center}

by~\eq{A_d-size} and its diameter is 
$
d  = a + b - c
$
by \Cref{diameter}.
To complete the proof, it remains to maximize
$(a+1)(b+1) - 2c(c+1)$ subject to the
constraint\, $a + b - c = d$.
The solution to this simple optimization problem is 
given in Table\,1.~~\qed\vspace{6ex}
\end{proof}

\centerline{%
\psfig{figure=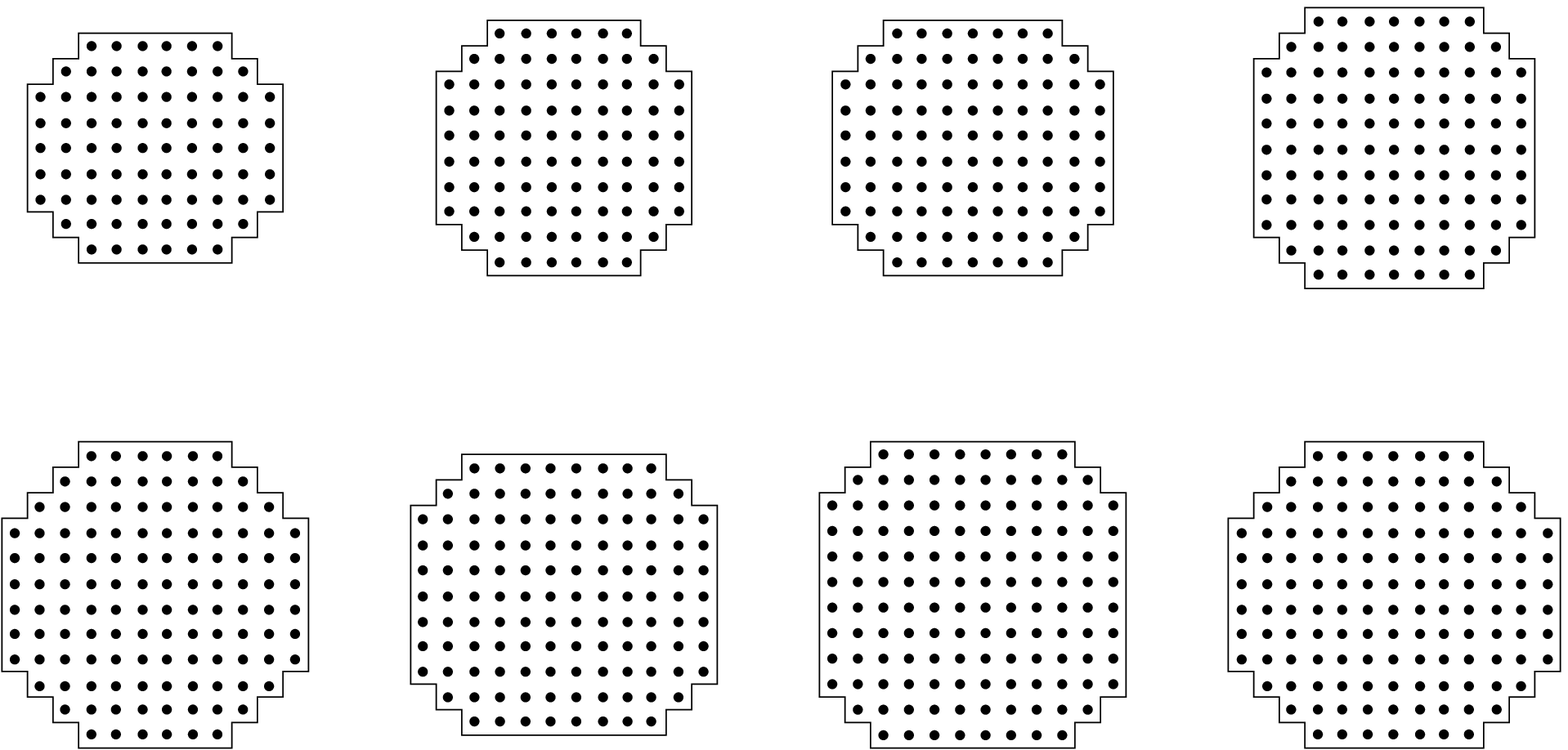,width=6.00in,silent=}}
\vspace{1.00ex}

\begin{center}
{\bf Figure\,5.}
{\sl Optimal tristance anticodes in $\Gplus$ of diameter $d = 15,16,\ldots,20$}
\vspace{.5ex}
\end{center}

\Tref{optimal-tristance} completely characterizes the optimal 
tristance anticodes of a given diameter in the grid graph~$\Gplus$.
Some examples of such anticodes are illustrated in\vspace*{-1.00ex} 
Figure\,5.\pagebreak[3.99]

\vspace{3ex}
\section{Optimal tristance anticodes in the {\Large $\infty$}
         and \protect\raisebox{-0.10ex}{\LARGE\rm $\hexagon$} models}
\vspace{-1.5ex} 
\label{sec3}

We now classify the optimal tristance anticodes in two related
graphical models: the infinity graph $\Gstar$ and the hexagonal
graph $\Ghex$ (defined in \S1). 
To obtain the classification~for~$\Gstar$, we make use of a mapping
$\varphi:\R^2 \to \R^2$ which takes $\Gstar$ into a power
graph of $\Gplus$. 
For the hexagonal graph $\Ghex$, we first derive an expression
for the corresponding tristance $\dhex(\cdot)$, and then follow
the same line of argument as in the previous section.

\subsection{Optimal tristance anticodes in the infinity graph}
\vspace{-1ex}
\label{sec3.1}

\looseness=-1
It is easy to see that the unique, up to translation, optimal 
anticode of diameter $d$~in~$\Gstar$ is the~square 
$
\cS_d
= 
\left\{\,
(x,y) \in \Z^2 \,:\, 0 \le x \le d \, \text{ and } \, 0 \le y \le d
\,\right\}
$.
To deal with tristance anticodes, 
we~first need an expression for tristance in $\Gstar$. It turns out that
tristance in $\Gstar$ is related~to~tri\-stance in $\Gplus$ via the 
mappings $\varphi:\R^2 \to \R^2$ and $\varphi\inv\!:\R^2 \to \R^2$
defined by
\Eq{
\label{varphi}
\varphi(x,y) \ = \ (x-y,x+y)
\hspace{3ex}
\text{and}
\hspace{3ex}
\varphi\inv(x,y) 
\ = \ 
\left(\frac{x+y}{2}\,,\: \frac{y-x}{2}\right)
}
Geometrically, the mapping $\varphi$ is simply a rotation 
by an angle of $\pi/4$ followed by scaling 
by a factor of $\sqrt{2}$. 
Note that $\varphi(\Z^2) = D_2$, where 
$D_2 = \left\{(x,y) \in \Z^2 : x + y \equiv 0\!\! \mod 2 \right\}$~is
the two-dimensional checkerboard lattice. 
Also note that for all $z_1,z_2 \in \Z^2$, we have
\Eq{
\label{dinfty}
d^{\star}(z_1,z_2) 
\ = \ 
\frac{|(x_1 - y_1) - (x_2 - y_2)|}{2} 
\ + \
\frac{|(x_1 + y_1) - (x_2 + y_2)|}{2} 
\ = \
\frac{d\bigl(\varphi(z_1),\varphi(z_2)\bigr)}{2}
}
\looseness=-1
Now define the
graph $\varphi(\Gstar) = (V,E)$~as~follows: $V = \varphi(\Z^2) = D_2$
and $\{z_1,z_2\} \in E$~if~and only if
$\{\varphi\inv(z_1), \varphi\inv(z_2)\}$ is an edge in~$\Gstar$.
Clearly, the graphs $\Gstar$ and $\varphi(\Gstar)$ are isomorphic.
It follows from~\eq{dinfty} that 
the edges of $\varphi(\Gstar)$ are precisely the 
paths of length~$2$~in~the grid graph $\Gplus$.
This fact was used in~\cite{EV} to prove the following theorem.

\begin{theorem}
\label{star-tristance}
Let  $z_1 = (x_1,y_1)$, $z_2 = (x_2,y_2)$,
$z_3 = (x_3 ,y_3 )$ be three distinct 
points in $\Z^2$,
and let $z'_1,z'_2,z'_3 \in D_2$
denote their images under $\varphi$.
Then
$
\dstar(z_1,z_2,z_3)
= 
\ceilenv{d_3(z'_1,z'_2,z'_3)/2}
$.
\end{theorem}

\Tref{star-tristance} and the mappings in~\eq{varphi} make it 
possible to classify the optimal tristance anticodes in $\Gstar$
using the classification of tristance anticodes in $\Gplus$
carried out in \S2.

\begin{theorem}
\label{star-centeroct}
Let $z_1 = (x_1,y_1)$, $z_2 = (x_2,y_2)$ be distinct 
points in $\Z^2$ and assume,~w.l.o.g., 
that $x_2 - x_1 \ge |y_2 - y_1|$ 
so that $d^{\star}(z_1,z_2) = x_2 - x_1$.
Let $\cA^{\star}_d(z_1,z_2)$~be~the optimal tristance anticode in $\Gstar$
of diameter $d \ge d^{\star}(z_1,z_2)$ centered about $z_1$ and $z_2$. 
Let $c = d - d^{\star}(z_1,z_2)$.
Then $\cA^{\star}_d(z_1,z_2)$ consists of all $z = (x,y)$ in $\Z^2$ 
such that 
\begin{eqnarray}
\label{star-centeroct1}
x_1 - c \ \le\hspace{-1.25ex} &x& \hspace{-1.25ex}\le\ x_2 + c\\
x_1 + y_1 - 2c\ \le\hspace{-1.25ex} &x+y& \hspace{-1.25ex}\le\ x_2 + y_2 + 2c\\
x_1 - y_1 - 2c\ \le\hspace{-1.25ex} &x-y& \hspace{-1.25ex}\le\ x_2 - y_2 + 2c\\
\label{star-centeroct4}
(x_1 + y_1) - (x_2 - y_2) - 2c\ \le\hspace{-1.25ex} &2y&  
\hspace{-1.25ex}\le\ (x_2 + y_2) - (x_1 - y_1) + 2c
\end{eqnarray}
\end{theorem}

\begin{proof}
Let $z'_1,z'_2 \in D_2$ be the images of $z_1$ and $z_2$ 
under $\varphi$. Let $z \in \cA^{\star}_d(z_1,z_2)$.
Then $\dstar(z_1,z_2,z) \le d$ 
and $d_3(z'_1,z'_2,\varphi(z)) \le 2d$ by \Tref{star-tristance},
so that $\varphi(z) \in \cA_{2d}(z'_1,z'_2)$.~~Since\pagebreak[3.99]
$\varphi(z) \in D_2$ for all $z \in \Z^2$, it 
follows that 
$
\varphi\bigl(\cA^{\star}_d(z_1,z_2)\bigr) 
\subseteq 
\cA_{2d}(z'_1,z'_2) \cap D_2
$. 
Conversely, let 
$z' \in \cA_{2d}(z'_1,z'_2) \cap D_2$.
Then $d_3(z'_1,z'_2,z') \le 2d$ and 
$\dstar(z_1,z_2,\varphi\inv(z')) \le d$
by \Tref{star-tristance}.
Hence 
$
\varphi\inv
\bigl(\cA_{2d}(z'_1,z'_2) \cap D_2\bigr) \subseteq \cA^{\star}_d(z_1,z_2)
$
and therefore
\Eq{
\label{oct-mapping}
\varphi\bigl(\cA^{\star}_d(z_1,z_2)\bigr) 
\ = \
\cA_{2d}(z'_1,z'_2) \cap D_2
}
In view of~\eq{oct-mapping} and \Tref{centeroct}, 
a point $z \in \Z^2$ belongs to $\cA^{\star}_d(z_1,z_2)$
if and only if~$\varphi(z)$~satisfies 
conditions~\eq{centeroct1}--\eq{centeroct2} for
$z'_1 = \varphi(z_1)$ and $z'_2 = \varphi(z_2)$,
with $c$ replaced by $2d - 2\dstar(z_1,z_2)$ in view of~\eq{dinfty}.
This is precisely the property expressed 
by \eq{star-centeroct1}--\eq{star-centeroct4}.~~\qed\vspace{1.00ex}
\end{proof}

\begin{corollary}
\label{star-shape}
Let $\cA^{\star}_d$ be an optimal tristance anticode of diameter $d$ 
in $\Gstar$. 
Then $\cA^{\star}_d$
is an octagon $\sO(\alpha_1,\alpha_2,\ldots,\alpha_8)$
for some $\alpha_1,\alpha_2,\ldots,\alpha_8 \in \Z$.
\vspace{-.5ex}
\end{corollary}

\begin{proof}
It is obvious from equations~\eq{star-centeroct1}--\eq{star-centeroct4}
in \Tref{star-centeroct} that the set $\cA^{\star}_d(z_1,z_2)$~is 
an octagon. Since the closure property of \Lref{intersection}
holds regardless of a particular~distance model, the corollary 
now follows in exactly the same way as \Lref{shape}.~~\qed\vspace{1.00ex}
\end{proof}

As in~\eq{new-oct1}--\eq{new-oct2}, we can use translations
in $\Z^2$ to write the octagon 
$\cA^{\star}_d$ 
as 
$
\sO(a,\kern-1pt b,\kern-1pt c_0, \kern-1pt c_1,\kern-1pt c_2,\kern-1pt c_3)$ 
for some $a,b,c_i \in \Z$.
Then the cardinality of $\cA^{\star}_d$ is given by~\eq{A_d-size},
and the next step is to determine its diameter $d$ as a function
of the parameters $a,b,c_0,c_1,c_2,c_3$.

\begin{lemma}
\label{star-diameter}
The diameter of an optimal tristance 
anticode $\cA^{\star}_d = \sO(a,b,c_0,c_1,c_2,c_3)$
in the infinity graph $\Gstar$ is given by
\Eq{
\label{star-d}
d 
\ = \
a+b \ - \ \floorenv{\, 
\frac{\min\{a+c_0+c_1,a+c_2+c_3,b+c_0+c_3,b+c_1+c_2\} }{2}\,}
\vspace{2ex}
}
\end{lemma}

\begin{proof}
We again make use of the mapping in~\eq{varphi}, in conjunction
with \Cref{diameter}. First consider the set $\varphi(\cA^{\star}_d)$.
Even though $\cA^{\star}_d = \sO(a,b,c_0,c_1,c_2,c_3)$, the set
$\varphi(\cA^{\star}_d)$ is, in general, not an octagon, since
$\varphi(\cA^{\star}_d) \subset D_2$. However, we can convert 
this into~an~octagon by adjoining the ``missing'' points as
follows:
\Eq{
\label{A'-def}
\cA'
\ \ \deff \ \
\varphi(\cA^{\star}_d) \cup
\bigl\{\, 
z \in \Z^2 ~:~
\text{at least $3$ of the $4$ neighbors of $z$ in $\Gplus$
are in $\varphi(\cA^{\star}_d)$}
\,\bigr\}
}
A straightforward analysis of the effect of the mapping $\varphi$
on~\eq{new-oct1}--\eq{new-oct2} now shows that $\cA'$ is an 
octagon $(c_3-b,c_0) + \sO(a',b',c'_0,c'_1,c'_2,c'_3)$, where\vspace{.75ex}
\begin{eqnarray}
\label{oct'-1}
a' \ = \ a + b - (c_1{+}c_3) \, ,
&\hspace{1ex}&
c'_0 \ = \ b - (c_0{+}c_3) \, ,
\hspace{5ex}
c'_1 \ = \ a - (c_0{+}c_1) 
\\[1ex]
\label{oct'-2}
b' \ = \ a + b - (c_0{+}c_2) \, ,
&&
c'_2 \ = \ b - (c_1{+}c_2) \, ,
\hspace{5ex}
c'_3 \ = \ a - (c_2{+}c_3) 
\end{eqnarray}
Let $d'$ denote the diameter of $\cA'$, and let $z_1,z_2,z_3$
be diametric points in $\cA^{\star}_d$. Then 
$\varphi(z_1)$, $\varphi(z_2)$, $\varphi(z_3)$ are in $\cA'$
and their tristance in $\Gplus$ is at least $2d-1$ by
\Tref{star-tristance}. Hence $d' \ge 2d-1$. 
Now let $z'_1,z'_2,z'_3$ be diametric points in $\cA'$.
If $z'_1 \not\in \varphi(\cA^{\star}_d)$, then it has~at~least 
three neighbors in $\varphi(\cA^{\star}_d)$ by~\eq{A'-def}.
By \Tref{tristance-formula}, this means that 
we can replace $z'_1,z'_2,z'_3$ by another diametric 
configuration $z''_1,z'_2,z'_3 \in \cA'$, where $z''_1$ 
is a neighbor of $z'_1$ such that $z''_1 \in \varphi(\cA^{\star}_d)$.
Repeating the argument for $z'_2$ and $z'_3$, we can find points 
$z'_1,z'_2,z'_3 \in \varphi(\cA^{\star}_d)$ 
such that $d_3(z'_1,z'_2,z'_3) = d'$.
We now have 
$
\dstar(\varphi\inv(z'_1),\varphi\inv(z'_2),\varphi\inv(z'_2)) 
= \ceilenv{d'/2} \le d
$,
in view of \Tref{star-tristance}.
Hence $d = \ceilenv{d'/2}$.
But $d' = a' + b' - \min\{c'_0,c'_1,c'_2,c'_3\}$ by
\Cref{diameter}. The lemma now follows straightforwardly 
from~\eq{oct'-1}--\eq{oct'-2}.~~\qed\pagebreak[3.99]
\end{proof}

\begin{theorem}
\label{optimal-star-tristance}
Let $\cA^{\star}_d$ be an optimal tristance anticode of diameter~$d$
in the infinity~graph. Then\vspace{-.25ex}
$$
|\cA^{\star}_d|
\ = \
\ceilenv{\frac{4d^2 + 8d + 2}{7}}
\vspace{.5ex}
$$
Moreover, up to rotation by an angle of $\pi/2$ and translation, 
$\cA_d^{\star} = \sO(a,b,c_0,c_1,c_2,c_3)$~whe\-re the parameters 
$a,b,c_0,c_1,c_2,c_3$ are given as a function of $d$ in 
Table\,2.\vspace{1.50ex}
$$
\begin{array}{|c|c|c|c|c|}
\hline \hline
 & & & & \\[-1.50ex]
\,d \hspace{-1.50ex}\pmod{7} & a, b & c_0,c_2 & c_1,c_3 & |\cA^{\star}_d| 
\\[0.75ex]
\hline\hline
 & & & & \\[-1.50ex]
0  & \frac{6d}{7} & \frac{2d}{7} & \frac{2d}{7} &  
\frac{4d^2+8d+7}{7} 
\\[0.55ex]
\hline
 & & & & \\[-1.95ex]
1  & \frac{6d+1}{7} & \frac{2d-2}{7} & \frac{2d+5}{7} &  
\frac{4d^2+8d+2}{7} 
\\[0.55ex]
\hline
 & & & & \\[-1.95ex]
2  & \frac{6d+2}{7} & \frac{2d+3}{7} & \frac{2d+3}{7} &  
\frac{4d^2+8d+3}{7} 
\\[0.55ex]
\hline
 & & & & \\[-1.95ex]
3  & \frac{6d-4}{7} & \frac{2d-6}{7} & \frac{2d-6}{7} &  
\frac{4d^2+8d+3}{7} 
\\[0.55ex]
\hline
 & & & & \\[-1.95ex]
4  & \frac{6d-3}{7} & \frac{2d-8}{7} & \frac{2d-1}{7} &  
\frac{4d^2+8d+2}{7} 
\\[0.55ex]
\hline
 & & & & \\[-1.95ex]
5  & \frac{6d-2}{7} & \frac{2d-3}{7} & \frac{2d-3}{7} &  
\frac{4d^2+8d+7}{7} 
\\[0.55ex]
\hline
 & & & & \\[-1.95ex]
6  & \frac{6d-1}{7} & \frac{2d-5}{7} & \frac{2d+2}{7} &  
\frac{4d^2+8d+4}{7} 
\\[1.15ex]
\hline\hline
\end{array}
$$

\begin{center}
{\bf Table\,2.}
{\sl Parameters of optimal tristance anticodes in the infinity graph $\Gstar$}
\vspace{.5ex}
\end{center}
\end{theorem}

\begin{proof}
In view of \Cref{star-shape} and \Lref{star-diameter},
we need to maximize the cardinality of $\cA^{\star}_d$ 
given by~\eq{A_d-size} subject to the constraint~\eq{star-d}.
Let $t = a + b - 2d + \max\{a,b\}$. If $t$ is even, then 
choosing $c_0 = c_1 = c_2 = c_3 = t/2$ satisfies~\eq{star-d}
and maximizes~\eq{A_d-size}. If $t$ is odd, then the 
corresponding extremal values are $c_0 = c_2 = (t-1)/2$
and $c_1 = c_3 = (t+1)/2$ (or vice versa). 
If we now assume w.l.o.g.\ that $a \ge b$, then
the cardinality of $\cA^{\star}_d$ is given by
$
(a+1)(b+1) - \floorenv{ (2a + b - 2d + 1)^2/{2} }
$.
It remains to maximize this expression, 
subject to $a \ge b$.
The solution to this
optimization problem is given in Table\,2.~~\qed\vspace{1.00ex}
\vspace{5ex}
\end{proof}

\centerline{%
\psfig{figure=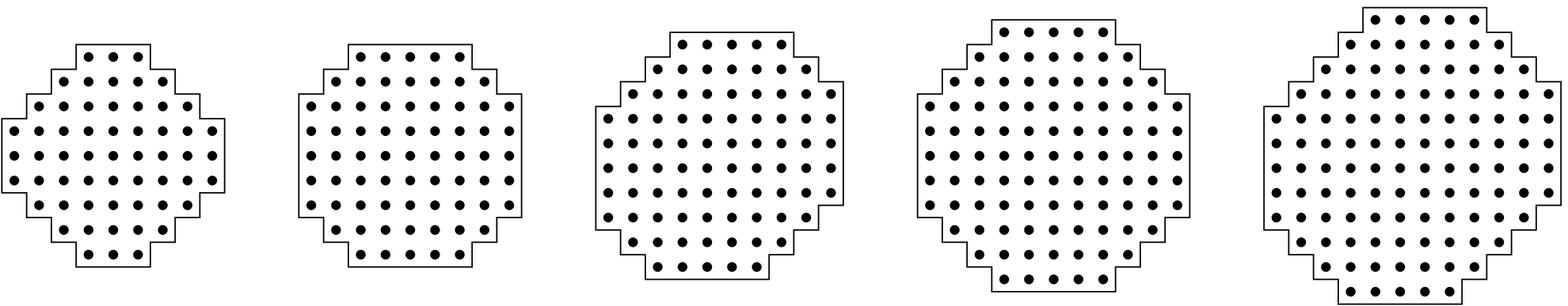,width=6.20in,silent=}}

\begin{center}
{\bf Figure\,6.}
{\sl Optimal tristance anticodes in $\Gstar$ of diameter $d = 9,10,\ldots,13$}
\vspace{.25ex}
\end{center}

\Tref{optimal-star-tristance} completes our classification of 
optimal tristance anticodes in the infinity graph~$\Gstar$. 
Some examples of such anticodes are illustrated in\vspace*{-1.00ex} 
Figure\,6.\pagebreak[3.99]

\newpage

\subsection{Optimal tristance anticodes in the hexagonal graph}
\vspace{-1ex}
\label{sec3.2}

\looseness=-1
Many different coordinate systems for the hexagonal lattice $A_2$
are known~\cite{SPLAG}. For our purposes, it would be most convenient
to identify $A_2$ with the Eisenstein integers~\cite[p.\,53]{SPLAG}.
That is, we write $A_2 = \{ x + \omega y\, : \,x,y \in \Z \}$, where
$\omega = -\shalf{1} + \sqrthalf{3}\, i\,$ is a complex cube root of 
unity.$^*$\footnote
{\\
\hspace*{-1.25ex}${}^*$\,Note that 
the numbers $1$, $\omega$, and $\omega^2$ represent
the three different edge orientations in the hexagonal graph
(cf.~Figure\,7).
However, two rather than three coordinates suffice to describe $\Ghex$,
since $\omega^2 = \bar{\omega} = -1 - \omega$.
\vspace*{-3.5ex}
}
Thus a generic vertex $v$ of $\Ghex$ will be written as $v = (x,y)$,
with~the~understanding that $v = x + \omega y$. The resulting labeling
of the hexagonal graph is shown
in Figure\,7.\vspace{3.5ex}

\centerline{%
\psfig{figure=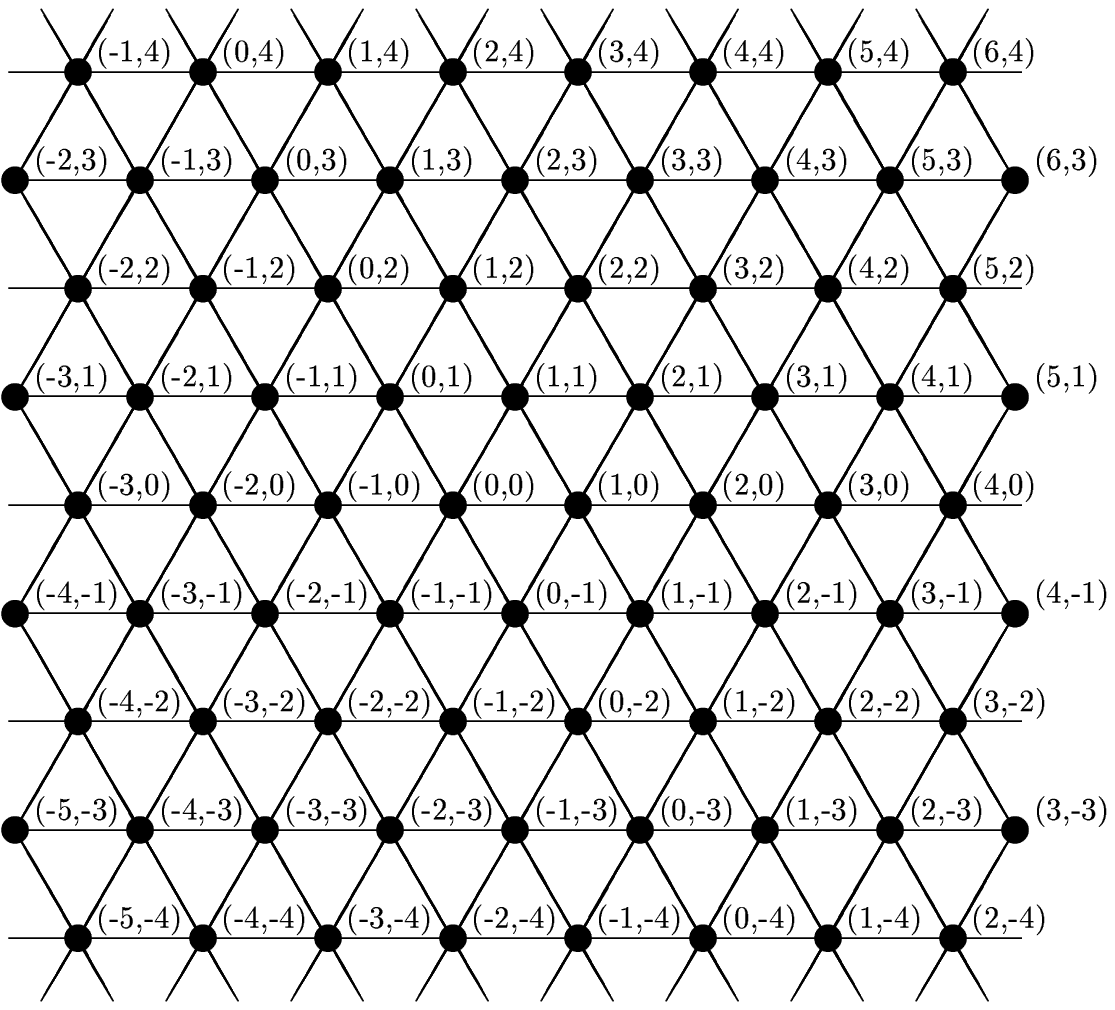,width=2.75in,angle=90,silent=}}
\vspace{-1.50ex}

\begin{center}
{\bf Figure\,7.}
{\sl The graph $\Ghex$ with vertices labeled by the Eisenstein integers}
\vspace{.25ex}
\end{center}

Our first task is to find expressions for distance and tristance in 
$\Ghex$. To this end, let us introduce the following notation:
given $a,b \,{\in}\, \Z$, we shall write
$\,\hmax\{a,b\} = \max\{a,b,0\}$~and 
$\hmin\{a,b\} = \min\{a,b,0\}$.
Now let $v_1 = (x_1,y_1)$ and $v_2 = (x_2,y_2)$ be two 
arbitrary points of~$A_2$. It is easy to see that the distance between 
$v_1$ and $v_2$ in $\Ghex$ is given by
\Eq{
\label{hex-distance}
\dtwohex(v_1,v_2)
\ = \
\hmax\,\{x_1\,{-}\,x_2,y_1\,{-}\,y_2\}
\ - \ 
\hmin\,\{x_1\,{-}\,x_2,y_1\,{-}\,y_2\} 
}
Note that if $(x_1{-}x_2)(y_1{-}\,y_2) \le 0$, then 
\eq{hex-distance} reduces to 
$
\dtwohex(v_1,v_2)
 = 
|x_1{-}x_2|\, + \,|y_1{-}\,y_2|
$~while 
if $(x_1\,{-}\,x_2)(y_1\,{-}\,y_2) \ge 0$
then
$
\dtwohex(v_1,v_2)
 = 
\max\{|x_1\,{-}\,x_2|,|y_1\,{-}\,y_2|\}
$.
Thus distance in $\Ghex$~is, in a sense, ``half-way'' between 
the $L_1$-distance of $\Gplus$ and the $L_{\infty}$-distance
of $\Gstar$. Deriving an expression for tristance in $\Ghex$ 
is a bit more involved. First, we need a lemma.

\begin{lemma}
\label{hextristance-lemma}
Let  $v_1$, $v_2$, and $v_3$ be distinct 
points in $A_2$. Then there exists
a point $v \,{\in}\, A_2$, such that
$
\dhex(v_1,v_2,v_3)
 = 
\dtwohex(v_1,v) + \dtwohex(v_2,v) + \dtwohex(v_3,v)
$.
\end{lemma}

\begin{proof}
By definition, $\dhex(v_1,v_2,v_3)$ is the number of edges in
a minimal spanning tree for $v_1$, $v_2$, $v_3$ in the hexagonal
graph $\Ghex$. Let $T$ be such a tree. Further, for $i = 1,2,\ldots,6$,
let $\eta_i$ denote the number of vertices of degree $i$ in $T$.
First observe that $\eta_1 \le 3$. 
Indeed, if there is a leaf in $T$ that is not one of $v_1$, $v_2$, $v_3$,
then we could remove this leaf 
along with the single\pagebreak[3.99] edge incident upon it, to obtain 
a smaller spanning tree for $v_1,v_2,v_3$ in $\Ghex$.
Now, the order of $T$ is 
$|V| = \eta_1 + \eta_2 + \eta_3 + \eta_4 + \eta_5 + \eta_6$,
while its size is 
$|E| = (\eta_1 + 2\eta_2 + 3\eta_3 + 4\eta_4 + 5\eta_5 + 6\eta_6)/2$.
Since $T$ is a~tree, we have $|V| - |E| = 1$.
With $|V|$, $|E|$ expressed in terms of $\eta_1,\eta_2,\ldots,\eta_6$,
this condition is equivalent to\vspace{-.5ex}
$$
\eta_3 + 2\eta_4 + 3\eta_5 + 4\eta_6
\ = \
\eta_1 - 2
\ \le \
1
$$
\looseness=-1
It follows that $\eta_4 = \eta_5 = \eta_6 = 0$ and $\eta_3 \le 1$.
In other words, there are only two possible configurations for $T$:
either it is {\dfn star-like\/}, with a single vertex of degree $3$
and $v_1,v_2,v_3$~as its three leaves, or it is {\dfn snake-like\/}
with only some two of $v_1,v_2,v_3$ as leaves and all other vertices
of degree $2$. If $T$ is star-like, we take $v$ to be the unique vertex
of degree $3$~in~$T$. If $T$ is snake-like, we take $v$ to be one
of $v_1,v_2,v_3$, the one which is {not} a leaf 
in~$T$.~\qed\vspace{.5ex}
\end{proof}

\begin{theorem}
\label{hex-tristance}
Let  $v_1 = (x_1,y_1)$, $v_2 = (x_2,y_2)$,
$v_3 = (x_3 ,y_3 )$ be distinct 
points in $A_2$. Let $\xmid$ denote the
middle value among $x_1,x_2,x_3$ --- that is,
if $x'_1,x'_2,x'_3$ is a permutation of $x_1,x_2,x_3$
such that $x'_1 \le x'_2 \le x'_3$, then $\xmid = x'_2$.
Let $\ymid$ be similarly defined.
Then\vspace{-.75ex}
\Eq{
\label{d3-hex}
\dhex(v_1,v_2,v_3)
\, = \,
\sum_{i=1}^{3}
\left( \Strut{3.5ex}{0ex}
\hmax\,\bigl\{x_i \,{-}\, \xmid,\, y_i \,{-}\, \ymid\bigr\} 
\ - \ 
\hmin\,\bigl\{x_i \,{-}\, \xmid,\, y_i \,{-}\, \ymid\bigr\} 
\right)
\vspace{1.5ex}
}
\end{theorem}

\begin{proof}
Let $v = (x,y)$ be a point satisfying
$
\dhex(v_1,v_2,v_3)
 = 
\dtwohex(v_1,v) + \dtwohex(v_2,v) + \dtwohex(v_3,v)
$.
Such a point exists by \Lref{hextristance-lemma}. Then 
by~\eq{hex-distance} we have
\Eq{
\label{d3-hex-aux}
\dhex(v_1,v_2,v_3)
\, = \,
\sum_{i=1}^{3}
\left( \Strut{3.5ex}{0ex}
\hmax\,\bigl\{x_i \,{-}\, x,\, y_i \,{-}\, y\bigr\} 
\ - \ 
\hmin\,\bigl\{x_i \,{-}\, x,\, y_i \,{-}\, y\bigr\} 
\right)
}
Clearly, the expression in~\eq{d3-hex-aux} is an upper
bound on $\dhex(v_1,v_2,v_3)$ for all $x,y \in \Z$.
To establish the tristance, it remains to find $x,y \in \Z$
that minimize this expression.  
This is a tedious, but simple, optimization problem
(an optimal solution must
satisfy $x \,{\in}\, \{x_1,x_2,x_3\}$ and $y \,{\in}\, \{y_1,y_2,y_3\}$, 
so there are nine cases to consider). The reader~can~easily 
verify that $x = \xmid$ and $y = \ymid$ is indeed an optimal 
solution.~\qed\vspace{.75ex}
\end{proof}

{\bf Remark.} 
\looseness=-1
The expressions $\max\{\cdot,\cdot,0\}$ and $\min\{\cdot,\cdot,0\}$ in 
\eq{hex-distance}, \eq{d3-hex} arise from~the~asymmetry
in our coordinate system for $\Ghex$.
If, instead, we represent a generic
point of $A_2$ as $v = (x,y,z)$, with the understanding 
that $v = x + \omega y + \omega^2 z$, then \eq{hex-distance}
becomes
$$
\dtwohex(v_1,v_2)
\ = \
\max\,\{x_1\,{-}\,x_2,y_1\,{-}\,y_2,z_1\,{-}\,z_2\}
\ - \ 
\min\,\{x_1\,{-}\,x_2,y_1\,{-}\,y_2,z_1\,{-}\,z_2\}
$$
and \eq{d3-hex} should be modified accordingly.
In some sense, these expressions are more~natural,
since they reflect the three edge directions in $\Ghex$.
On the other hand, this coordinate system is redundant: 
$(x,y,z)$ and $(x-\delta,y-\delta,z-\delta)$ represent 
the same point of $A_2$ for all $\delta \,{\in}\, \Z$, since
$1 + \omega + \omega^2 = 0$. One can use this property
to zero-out any one of the three coordinates. Zeroing 
out the last coordinate by choosing $\delta = z$
(as we have done) is precisely the source for the
remnant zeros in $\max\{\cdot,\cdot,0\}$ and $\min\{\cdot,\cdot,0\}$ 
in \eq{hex-distance}, \eq{d3-hex}. \hfill$\Box$\hspace*{-1ex}
\vspace{+.5ex}

From here, we proceed along the lines of \S2.3. Let 
a {\dfn hexagon} $\sH(\al_1,\al_2,\ldots,\al_6)$ be~a~subset 
of $A_2$ defined by the inequalities
\Eq{
\label{hexagon}
\alpha_1 \,\le\, x \,\le\, \alpha_4\, ,
\hspace{4ex}
\alpha_2 \,\le\, y \,\le\, \alpha_5\, ,
\hspace{4ex}
\alpha_3 \,\le\, x-y \,\le\, \alpha_6
}
The next lemma and theorem show that the optimal centered anticode 
$\cA^{\hex}_d(v_1,v_2) \subset A_2$, centered about 
an arbitrary pair of points $v_1,v_2 \in A_2$, is a hexagon 
for all $d \ge \dtwohex(v_1,v_2)$.\vspace*{-1ex}\pagebreak[3.99]

\begin{lemma}
\label{hex-P}
Let $v_1 = (x_1,y_1)$, $v_2 = (x_2,y_2)$ be distinct 
points of $A_2$. Then $\dhex(v_1,v_2,v_3) = \dtwohex(v_1,v_2)$
if and only if $v_3$ belongs to the {\dfn bounding 
parallelepiped\/} of $v_1,v_2$,  which is a subset of $A_2$
defined by the inequalities
\begin{eqnarray}
\label{hex-P1}
\min\{x_1,x_2\}\ \le\hspace{-3.00ex} 
&x& 
\hspace{-3.00ex}\le\ \max\{x_1,x_2\}\\
\label{hex-P2}
\min\{y_1,y_2\}\ \le\hspace{-3.00ex} 
&y& 
\hspace{-3.00ex}\le\ \max\{y_1,y_2\}\\
\label{hex-P3}
\min\{x_1\,{-}\,y_1,x_2\,{-}\,y_2\}\ \le\hspace{-.75ex} 
&x-y& 
\hspace{-.75ex}\le\ \max\{x_1\,{-}\,y_1,x_2\,{-}\,y_2\}
\end{eqnarray}
\end{lemma}

\begin{proof}
Let $\cP(v_1,v_2)$ denote the bounding parallelepiped of $v_1,v_2$
(note that it is, indeed, a parallelepiped since one of 
\eq{hex-P1}--\eq{hex-P3} is always redundant).
Assume w.l.o.g.\ that $x_1 \le x_2$.
($\Leftarrow$) Suppose $v_3 \in \cP(v_1,v_2)$. Then \eq{hex-P1}, 
\eq{hex-P2} imply that $\xmid = x_3$ and $\ymid = y_3$ in~\eq{d3-hex}.
Thus \eq{d3-hex} reduces to
\Eq{
\label{lemma23-aux}
\dhex(v_1,v_2,v_3)
\ = \
\left\{
\begin{array}{l@{\hspace{5ex}}l}
(y_1 \,{-}\, y_3) - (x_1 \,{-}\, x_3) + (x_2 \,{-}\, x_3) - (y_2 \,{-}\, y_3) &
\text{if $y_1 \ge y_2$}\\[1ex]
\max\{x_2 {-}\, x_3,y_2 {-}\, y_3\} \ - \ \min\{x_1 {-}\, x_3,y_1 {-}\, y_3\} &
\text{if $y_1 \le y_2$}
\end{array}
\right.
}
If $y_1 \ge y_2$, then $\dhex(v_1,v_2,v_3) = \dtwohex(v_1,v_2)$ 
directly 
by \eq{lemma23-aux} and \eq{hex-distance}. 
If $y_1 \le y_2$, then~the~inequality
$\dhex(v_1,v_2,v_3) \le \dtwohex(v_1,v_2)$\linebreak[3.99]
follows by straightforward manipulation
from \eq{lemma23-aux}, \eq{hex-P3}.
($\Rightarrow$) Now suppose that $\dhex(v_1,v_2,v_3) = \dtwohex(v_1,v_2)$.
Then every minimal spanning tree for $v_1,v_2,v_3$ in $\Ghex$
must be snake-like, with $v_1,v_2$ as its leaves. Indeed, by 
\Lref{hextristance-lemma}~we~have 
$$
\dhex(v_1,v_2,v_3) 
\ = \
\dtwohex(v_1,v) + \dtwohex(v_2,v) + \dtwohex(v_3,v) 
\ = \
\dtwohex(v_1,v_2)
$$
for some $v \in A_2$. This is only possible if $\dtwohex(v_3,v) = 0$,
since $\dtwohex(v_1,v) + \dtwohex(v_2,v) \ge \dtwohex(v_1,v_2)$ by the
triangle inequality.
The fact that $\dtwohex(v_3,v) = 0$ implies that the third term  
in the summation of~\eq{d3-hex} is zero, which is only possible 
if $\xmid = x_3$ and $\ymid = y_3$. This~estab\-lishes 
\eq{hex-P1}~and~\eq{hex-P2}. Moreover, the expression for
$\dhex(v_1,v_2,v_3)$ in \eq{d3-hex} once again reduces 
to \eq{lemma23-aux}. If $y_1 \le y_2$,~then 
$\dhex(v_1,v_2,v_3) = \dtwohex(v_1,v_2)$ further
reduces to
$$
\max\{x_2 {-}\, x_3,y_1 {-}\, y_3\} \ - \ \min\{x_1 {-}\, x_3,y_1 {-}\, y_3\}
\ = \ 
\max\{x_2 {-}\, x_1,y_2 {-}\, y_1\}
$$
by \eq{lemma23-aux} and \eq{hex-distance}.
It is straightforward to show that this condition
is equivalent to \eq{hex-P3}. Otherwise, if $y_1 \ge y_2$, then 
\eq{hex-P1} and \eq{hex-P2} imply \eq{hex-P3}, and we are 
done.~\qed\vspace{.5ex}
\end{proof}

\begin{theorem}
\label{hex-centeroct}
Let $v_1 = (x_1,y_1)$, $v_2 = (x_2,y_2)$ be distinct 
points of $A_2$, and let $\cA^{\hex}_d(v_1,v_2)$
be the optimal tristance anticode in $\Ghex$
of diameter $d \ge \dtwohex(v_1,v_2)$ centered about $v_1$~and~$v_2$. 
Write $c = d - \dtwohex(v_1,v_2)$.
Then $\cA^{\hex}_d(v_1,v_2)$ consists of all $v = (x,y)$ in $A_2$ 
such that 
\begin{eqnarray}
\label{hex-centeroct1}
\min\{x_1,x_2\} - c \ \le\hspace{-3.00ex} 
&x& 
\hspace{-3.00ex}\le\ \max\{x_1,x_2\} + c\\
\min\{y_1,y_2\} - c \ \le\hspace{-3.00ex} 
&y& 
\hspace{-3.00ex}\le\ \max\{y_1,y_2\} + c\\
\label{hex-centeroct3}
\min\{x_1\,{-}\,y_1,x_2\,{-}\,y_2\} - c \ \le\hspace{-.75ex} 
&x-y& 
\hspace{-.75ex}\le\ \max\{x_1\,{-}\,y_1,x_2\,{-}\,y_2\} + c
\end{eqnarray}
\end{theorem}

\begin{proof}
When $c = 0$, the theorem follows immediately from \Lref{hex-P}.
Otherwise, it is easy to see that 
$v \in \cA^{\hex}_d(v_1,v_2)$
if and only if the distance in $\Ghex$ from $v$ to (the closest point of)
the {bounding parallelepiped} $\cP(v_1,v_2)$ is at most 
$c = d - \dtwohex(v_1,v_2)$.
This is precisely the property expressed 
by \eq{hex-centeroct1}--\eq{hex-centeroct3}.~~\qed\vspace{1.00ex}
\end{proof}

\begin{corollary}
\label{hex-shape}
Let $\cA^{\hex}_d$ be an optimal tristance anticode of diameter $d$ 
in $\Ghex$. 
Then $\cA^{\hex}_d$
is a hexagon $\sH(\alpha_1,\alpha_2,\ldots,\alpha_6)$
for some $\alpha_1,\alpha_2,\ldots,\alpha_6 \in \Z$.\pagebreak[3.99]
\end{corollary}

\begin{proof}
It is obvious from \Tref{hex-centeroct} and \eq{hexagon}
that the set $\cA^{\hex}_d(v_1,v_2)$ is a hexagon.
The fact that the intersection of any two hexagons is
a hexagon is also obvious (cf.~\Lref{octint}).
The corollary now follows in exactly the same way as 
\Lref{shape} and \Cref{star-shape}.~~\qed\vspace{1.00ex}
\end{proof}

Since $A_2$ is invariant under translation by a lattice point, we can 
again shift $\cA^{\hex}_d$ to the~ori\-gin, so that 
$\cA^{\hex}_d = \sH(0,0,\al_3,\al_4,\al_5,\al_6)$,
and then write it as $\sH(a,b,c_1,c_3)$, where~$a = \al_4$, $b = \al_5$,
$c_1 = \al_4 - \al_6$, and $c_3 = \al_3 + \al_5$
(cf.~Figure\,3). This is just
a special case of~\eq{new-oct1}--\eq{new-oct2}.

\begin{lemma}
\label{hex-diameter}
The diameter of an optimal tristance 
anticode $\cA^{\hex}_d = \sH(a,b,c_1,c_3)$ in the hexagonal 
graph $\Ghex$ is given by
$d = a + b - \min\{c_1,c_3\}$.\vspace{-.50ex}
\end{lemma}

\begin{proof}
\looseness=-1
Assume w.l.o.g.\ that $c_1 \le c_3$. Observe that we can further
assume w.l.o.g.\ that $0 \le c_1,c_3 \le \min\{a,b\}$; 
otherwise at least one of the inequalities in
\Eq{
\label{lemma26-aux}
0 \,\le\, x \,\le\, a\, ,
\hspace{4ex}
0 \,\le\, y \,\le\, b\, ,
\hspace{4ex}
c_3-b\,\le\, x-y \,\le\, a-c_1
}
is redundant, and the hexagon $\sH(a,b,c_1,c_3)$ can be translated
and/or re-parametrized so that $0 \le c_1,c_3 \le \min\{a,b\}$ holds.
Thus the points $v_1 = (0,0)$, $v_2 = (a,c_1)$, 
and $v_3 = (a,b)$\linebreak[3.99]
belong to $\sH(a,b,c_1,c_3)$ by~\eq{lemma26-aux}, and 
$\dhex(v_1,v_2,v_3) = a \kern1.5pt+\kern1.5pt b \kern1.5pt-\kern1.5pt c_1$ 
by~\eq{d3-hex}.
Now let $v_1,v_2,v_3$~be arbitrary points in $\sH(a,b,c_1,c_3)$.
We assume w.l.o.g.\ that $x_1 \le x_2 \le x_3$, and distinguish
between six cases. In each case, we compute $\dhex(v_1,v_2,v_3)$
using~\eq{d3-hex}.

\noindent{\bf Case\,1:} $y_3 \le y_2 \le y_1$. 
Then $\dhex(v_1,v_2,v_3) = (x_3-x_1) + (y_1-y_3)$.
\vspace{-.50ex}

\noindent{\bf Case\,2:} $y_3 \le y_1 \le y_2$. 
Then $\dhex(v_1,v_2,v_3) = (x_3-x_1) + (y_2-y_3)$.
\vspace{-.50ex}

\noindent{\bf Case\,3:} $y_2 \le y_3 \le y_1$. 
Then $\dhex(v_1,v_2,v_3) = (x_3-x_1) + (y_1-y_2)$.
\vspace{-.50ex}

\noindent{\bf Case\,4:} $y_2 \le y_1 \le y_3$. 
Then $\dhex(v_1,v_2,v_3) = (x_2-x_1) + (y_1-y_2) + \max\{x_3-x_2,y_3-y_1\}$.\vspace{-.50ex}

\noindent{\bf Case\,5:} $y_1 \le y_3 \le y_2$. 
Then $\dhex(v_1,v_2,v_3) = (x_3-x_2) + (y_2-y_3) + \max\{x_2-x_1,y_3-y_1\}$.\vspace{-.50ex}

\noindent{\bf Case\,6:} $y_1 \le y_2 \le y_3$. 
Then $\dhex(v_1,v_2,v_3) = \max\{x_3-x_2,y_3-y_2\} + \max\{x_2-x_1,y_2-y_1\}$.

In each of these cases, it is straightforward to show 
that 
$\dhex(v_1,v_2,v_3) \le a + b - c_3$ 
or
$\dhex(v_1,v_2,v_3) \le a + b - c_1$ by~\eq{lemma26-aux},
and the lemma follows.~~\qed
\end{proof}

\begin{theorem}
\label{optimal-hex-tristance}
Let $\cA^{\hex}_d$ be an optimal tristance anticode of diameter~$d$
in the hexagonal graph. Then\vspace{-.25ex}
$$
|\cA^{\hex}_d|
\ = \
\ceilenv{\frac{d^2 + 3d + 2}{3}}
\ = \
\ceilenv{\frac{(d+1)(d+2)}{3}}
\vspace{.5ex}
$$
Moreover, up to rotation by an angle of $\pi/3$ and translation, 
$\cA_d^{\hex} = \sH(a,b,c_1,c_3)$ where the parameters 
$a,b,c_1,c_3$ are given as a function of $d$ in 
Table\,3.\vspace{-.50ex}
\end{theorem}

\begin{proof}
The optimal tristance anticode is a hexagon $\sH(a,b,c_1,c_3)$ by
\Cref{hex-shape}. Its cardinality is  
$
|\cA^{\hex}_d| = (a+1)(b+1) - \shalf{1}c_1(c_1+1) - \shalf{1}c_3(c_3+1)
$
as in~\eq{A_d-size}, and its diameter is $a + b - \min\{c_1,c_3\}$
by \Lref{hex-diameter}. Clearly, the choice $c_1 = c_3 = c$
maximizes $|\cA^{\hex}_d|$ for a~given diameter. It remains 
to maximize $(a+1)(b+1) - c(c+1)$ subject to $a+b-c = d$. 
The solution to this optimization problem is given in 
Table\,3.~~\qed\vspace{1ex}
\end{proof}

Some examples of optimal tristance anticodes in $\Ghex$ are illustrated 
in Figure\,8. It can be seen from Table\,3
that such anticodes are regular hexagons if and only if 
$d \equiv 0\!\pmod{3}$.

$$
\begin{array}{|c|c|c|c|c|}
\hline \hline
& & & & \\[-1.50ex]
\,d \hspace{-1.50ex}\pmod{3} & a & b & c_1,c_3 & |\cA^{\hex}_d| 
\\[0.75ex]
\hline\hline
& & & & \\[-1.50ex]
0  & \frac{2d}{3} & \frac{2d}{3} & \frac{d}{3} & \frac{d^2+3d+3}{3} 
\\[0.55ex]
\hline
& & & & \\[-1.50ex]
1  & 
\begin{array}{c}
\frac{2d+1}{3}\\[.75ex]
\frac{2d+1}{3}
\end{array}
& 
\begin{array}{c}
\frac{2d+1}{3}\\[.75ex]
\frac{2d-2}{3}
\end{array}
&
\begin{array}{c}
\frac{d+2}{3}\\[.75ex]
\frac{d-1}{3}
\end{array}
& 
\frac{d^2+3d+2}{3}
\\[2.25ex]
\hline
& & & & \\[-1.50ex]
2  & 
\begin{array}{c}
\frac{2d-1}{3}\\[.75ex]
\frac{2d+2}{3}
\end{array}
& 
\begin{array}{c}
\frac{2d-1}{3}\\[.75ex]
\frac{2d-1}{3}
\end{array}
&
\begin{array}{c}
\frac{d-2}{3}\\[.75ex]
\frac{d+1}{3}
\end{array}
& 
\frac{d^2+3d+2}{3}
\\[3.15ex]
\hline\hline
\end{array}
$$

\begin{center}
{\bf Table\,3.}
{\sl Parameters of optimal tristance anticodes in the hexagonal graph $\Ghex$}
\vspace{5.0ex}
\end{center}

\centerline{%
\psfig{figure=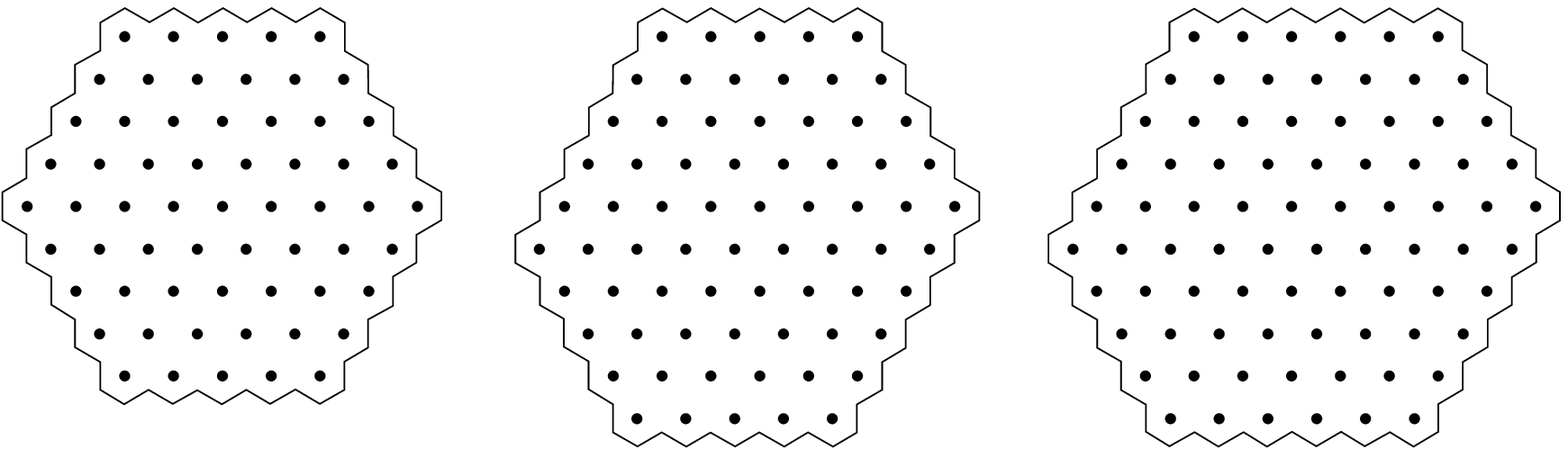,width=6.20in,silent=}}
\vspace{.50ex}

\begin{center}
{\bf Figure\,8.}
{\sl Optimal tristance anticodes in $\Ghex$ of diameter $d = 12,13,14$}
\vspace{.25ex}
\end{center}

{\bf Remark.}
Using similar methods, we can also characterize the optimal
{\sl distance\/} anticodes in the hexagonal graph. For even 
diameter $d$, such anticodes are regular hexagons (spheres 
in $\Ghex$) centered about a lattice point $v_0 = (x_0,y_0)$, namely
$$
\sS_{d/2}(v_0)
\ = \ 
\set{ (x,y) \in A_2 ~:~ 
|x - x_0| \le \frac{\raisebox{-1pt}{$d$}}{2}\,,\hspace{2ex}
|y - y_0| \le \frac{\raisebox{-1pt}{$d$}}{2}\,,\hspace{2ex}
|(x{-}y) - (x_0{-}y_0)| \le \frac{\raisebox{-1pt}{$d$}}{2}
\, }
$$
For odd $d$, optimal distance anticodes are again ``spheres'' 
of radius $(d+1)/2$, but no longer centered about a lattice
point. Specifically, an optimal distance anticode of an odd
diameter $d$ is given by
$$
\set{ (x,y) \in A_2 ~:~ 
|x - x'_0| \le \frac{\raisebox{-1pt}{$d+1$}}{2}\,,\hspace{2ex}
|y - y'_0| \le \frac{\raisebox{-1pt}{$d+1$}}{2}\,,\hspace{2ex}
|(x{-}y) - (x'_0{-}y'_0)| \le \frac{\raisebox{-1pt}{$d+1$}}{2}
\, }
$$
\looseness=-1
where $(x'_0,y'_0) = v_0 + \xi$, with $v_0 = (x_0,y_0)$ being 
an arbitrary point of $A_2$ and $\xi$ being one of
$(\sthird{1},-\sthird{1})$, 
$(\sthird{2},\sthird{1})$, 
$(\sthird{1},\sthird{2})$, 
$(-\sthird{1},\sthird{1})$, 
$(-\sthird{2},-\sthird{1})$, 
$(-\sthird{1},-\sthird{2})$.
Such anticodes can be construed as regular hexagons over $\R^2 \simeq \C$,
but they are not regular hexagons when viewed as subsets of $A_2$.
The cardinality of an optimal distance anticode of diameter $d$
is $1+ 3d(d+2)/4$
if $d$ is even, and $3(d+1)^2/4$ if $d$ is odd. \hfill$\Box$

\section{Higher dimensions and higher dispersions}
\vspace{-1.75ex} 
\label{sec4}

\looseness=-1
In general, extending our results for tristance anticodes in $\Z^2$ 
to higher dimensions and/or higher dispersions appears to be 
a~difficult problem. Nevertheless, we pursue such generalizations 
in~this~section, in part to illustrate the difficulties that arise
along the way.\vspace{-.25ex}

\looseness=-1
In \S4.1, we study tristance anticodes in $\Z^3$.
Here,~the~general approach developed in \S2 still works: we first characterize
the optimal centered anticodes $\cA_d(v_1,v_2) \subset \Z^3$ 
and~thus~determine, using \Lref{intersection}, the shape of an 
optimal unrestricted tristance anticode $\cA_d \subseteq \Z^3$.
The problem is that the expressions for the diameter and the cardinality 
of~$\cA_d$~are much more involved than their counterparts for 
$\Z^2$ in \eq{A_d-size} and~\eq{Cor13}. The resulting optimization task 
involves $23$ variables and does not appear~to~be tractable.
We~conjecture, however, that optimal tristance anticodes in $\Z^3$
satisfy a~certain symmetry condition. Subject to this conjecture,
we determine the parameters~of such anticodes and their cardinality.
\vspace{-.25ex}

\looseness=-1
In \S4.2, we consider quadristance anticodes in $\Z^2$. This serves
to illustrate a situation~whe\-re the approach of \S2 breaks down.
We can still characterize the optimal centered quad\-ristance
anticodes $\cA_d(z_1,z_2,z_3) \subset \Z^2$ and 
(the appropriate generalization of) \Lref{intersection} still
applies. However, such centered anticodes are no longer convex
and their general shape is \emph{not} preserved under intersection.
Thus \Lref{intersection} tells us nothing about~the~shape of 
unrestricted optimal quadristance anticodes in $\Z^2$.
We use the octagon shape to derive a~lower bound on the 
cardinality of such anticodes. We conjecture that
this bound~is, in fact, exact.
Observe, however, that shapes other than octagons occur among
optimal quadristance anticodes, at least for certain 
diameters (cf.~Figure\,12).

\vspace{-1.15ex}
\subsection{Optimal tristance anticodes in the grid graph of $\Z^3$}
\vspace{-1.15ex}
\label{sec4.1}

We first need an expression for tristance in $\Gplusthree$, the 
grid graph of $\Z^3$. Fortunately, the~tristance formula of 
\Tref{tristance-formula} 
easily generalizes to arbitrary dimensions.\vspace{-.50ex}

\begin{theorem}
\label{nD-tristance}
Let $v = (v_1,v_2,\ldots,v_n)$, 
$v' = (v'_1,v'_2,\ldots,v'_n)$, 
and
$v'' = (v''_1,v''_2,\ldots,v''_n)$~be
distinct points in $\Z^n$. Then\vspace{-1.75ex}
\Eq{
\label{d3-n}
d_3(v,v',v'')
\ = \ 
\sum_{i=1}^n
\left(
\max\{v_i,v'_i,v''_i\} \ - \ \min\{v_i,v'_i,v''_i\} 
\Strut{2.80ex}{0ex}
\right)
\vspace{1.25ex}
}
\end{theorem}

\begin{proof}
It is easy to see that, for all $i = 1,2,\ldots,n$, 
any spanning tree for $v,v',v''$~must 
contain at least 
$
\max\{v_i,v'_i,v''_i\} - \min\{v_i,v'_i,v''_i\}
$
edges that are parallel to the $i$-th coordinate
axis. Thus the sum on the right-hand side of~\eq{d3-n}
is a lower bound on $d_3(v,v',v'')$. To show that 
this bound holds with equality, we use induction on $n$,
with \Tref{tristance-formula} serving as the 
induction base. Assume w.l.o.g.\ that $v'_n \le v_n \le v''_n$
and let $u, w \in \Z^n$ be defined~by\hspace*{.6ex}
$$
u 
\ = \
(v'_1,v'_2,\ldots,v'_{n-1},v_n)
\hspace{4.5ex}\text{and}\hspace{5ex}
w 
\ = \
(v''_1,v''_2,\ldots,v''_{n-1},v_n)
$$
It takes $v_n - v'_n$ edges to connect $v'$ with $u$
and another $v''_n - v_n$ edges to connect $v''$ with~$w$, 
altogether 
$
v''_n - v'_n
= 
\max\{v_n,v'_n,v''_n\} \,-\, \min\{v_n,v'_n,v''_n\}
$ 
edges.
Since the points $u,v,w$~belong to the same coset 
of $\Z^{n-1}$ in $\Z^n$, the claim now follows
by induction hypothesis.~~\qed\vspace{.50ex}
\end{proof}

\looseness=-1
Next, we generalize to three dimensions the definition 
of a bounding rectangle in \S2.1. 
Let 
$v_1 = (x_1,y_1,z_1), v_2 = (x_2,y_2,z_2), \ldots, v_n = (x_n,y_n,z_n)$ 
be $n$ distinct points in~$\Z^3$. Then\pagebreak[3.99]
the {\dfn bounding cuboid\/}
of $v_1,v_2,\ldots,v_n$ is the smallest cuboid 
$\cC(v_1,v_2,\ldots,v_n)$ with edges parallel to the axes that 
contains all the $n$ points. Explicitly, define
$
\xmax = \max\{x_1,x_2,\ldots,x_n\}
$
and
$
\xmin = \min\{x_1,x_2,\ldots,x_n\}
$.
Let $\ymax$, $\ymin$, $\zmax$, and $\zmin$
be defined similarly. Then
$$
\cC(v_1,v_2,\ldots,v_n)
\ \ \deff \ 
\left\{ \Strut{2.5ex}{0ex}
(x,y,z) \in \Z^3 \,:\, 
\xmin \le x \le \xmax, \:
\ymin \le y \le \ymax, \: 
\zmin \le z \le \zmax 
\right\} 
$$
By \Tref{nD-tristance}, the tristance of \emph{any}
three points that lie in the bounding cuboid $\cC(v_1,v_2)$
of $v_1 = (x_1,y_1,z_1)$ and $v_2 = (x_2,y_2,z_2)$ is 
at most 
$
d(v_1,v_2)
= 
|x_1 - x_2| + |y_1 - y_2| + |z_1 - z_2|
$.
This immediately leads to the following 
characterization of optimal tristance 
anticodes in~$\Z^3$ that are {centered} 
about two given points $v_1,v_2 \in \Z^3$
(cf.~\Tref{centeroct}).

\begin{theorem}
\label{center-cuboid}
Let $v_1 = (x_1,y_1,z_1)$, $v_2 = (x_2,y_2,z_2)$ be distinct 
points~in~$\Z^3$. Let~$\cA_d(v_1,v_2)$
be the optimal tristance anticode in $\Gplusthree$
of diameter $d \ge d(v_1,v_2)$ centered about~$v_1$~and~$v_2$. 
Write $\delta = d - d(v_1,v_2)$.
Let $\xmax = \max\{x_1,x_2\}$, $\xmin = \min\{x_1,x_2\}$
with $\ymax$, $\ymin$~and $\zmax$, $\zmin$ defined similarly.
Then $\cA_d(v_1,v_2)$ consists of all $v = (x,y,z)$ in $\Z^3$ 
such that 
\begin{eqnarray}
\label{center-cuboid1}
\xmin -\delta \ \le\hspace{-4.75ex} 
& x &
\hspace{-4.75ex}\le\ \xmax + \delta
\\[-.30ex]
\ymin -\delta \ \le\hspace{-4.75ex}
& y &
\hspace{-4.75ex}\le\ \ymax + \delta
\\[-.30ex]
\zmin -\delta \ \le\hspace{-4.75ex} 
& z &
\hspace{-4.75ex}\le\ \zmax + \delta
\\[-.30ex] 
\xmin + \ymin - \delta \ \le\hspace{-2.75ex}
& x+y &
\hspace{-2.75ex}\le\ \xmax + \ymax + \delta
\\[-.30ex]
\xmin - \ymax - \delta \ \le\hspace{-2.75ex}
& x-y &
\hspace{-2.75ex}\le\ \xmax - \ymin + \delta
\\[-.30ex]
\xmin + \zmin - \delta \ \le\hspace{-2.75ex}
& x+z &
\hspace{-2.75ex}\le\ \xmax + \zmax + \delta
\\[-.30ex]
\xmin - \zmax - \delta \ \le\hspace{-2.75ex}
& x-z &
\hspace{-2.75ex}\le\ \xmax - \zmin + \delta
\\[-.30ex]
\ymin + \zmin - \delta \ \le\hspace{-2.75ex}
& y+z &
\hspace{-2.75ex}\le\ \ymax + \zmax + \delta
\\[-.30ex]
\ymin - \zmax - \delta \ \le\hspace{-2.75ex}
& y-z &
\hspace{-2.75ex}\le\ \ymax - \zmin + \delta
\\[-.30ex] 
\xmin + \ymin + \zmin - \delta \ \le\hspace{-1ex}
& x+y+z &
\hspace{-1ex}\le\ \xmax + \ymax + \zmax + \delta
\\[-.30ex]
\xmin - \ymax + \zmin - \delta \ \le\hspace{-1ex}
& x-y+z &
\hspace{-1ex}\le\ \xmax - \ymin + \zmax + \delta
\\[-.30ex]
\xmin + \ymin - \zmax - \delta \ \le\hspace{-1ex}
& x+y-z &
\hspace{-1ex}\le\ \xmax + \ymax - \zmin + \delta
\\[-.30ex]
\label{center-cuboid13}
\xmin - \ymax - \zmax - \delta \ \le\hspace{-1ex}
& x-y-z &
\hspace{-1ex}\le\ \xmax - \ymin - \zmin + \delta
\\[-2.85ex]
\nonumber
\end{eqnarray}
\end{theorem}

\begin{proof}
It follows from \Tref{nD-tristance} that
$\cA_d(v_1,v_2) = \cC(v_1,v_2)$ if $\delta = 0$.
Hence for $\delta > 0$,\linebreak[4.00] the set
$
\cA_d(v_1,v_2)
=
\{ \kern-.5pt v\kern-.5pt 
\in \Z^3 : d_3(v_1,v_2,v) \le d(v_1,v_2) \kern.75pt+\kern.75pt \delta\}
$
consists of all points \mbox{$(x,y,z)\kern-1pt \in \Z^3$} whose $L_1$-distance
from the bounding cuboid $\cC(v_1,v_2)$ is at most $\delta$.
It can be readily verified that this is precisely the set described 
by equations \eq{center-cuboid1} through 
\eq{center-cuboid13}.~~\qed\vspace{.5ex}
\end{proof}

The centered anticode $\cA_d(v_1,v_2)$ in 
\Tref{center-cuboid} is an example of a set
we call the \icos. In general, we define
an {\dfn \icos} $\sI(\al_1,\al_2,\ldots,\al_{26})$
as the set of all points of $\Z^3$ that lie within 
the convex polyhedron with $26$ faces, given by the inequalities
\begin{eqnarray}
\label{icosi-first}
&
\al_1 \,\le\, x \,\le\, \al_{14}\, ,
\hspace{8.0ex}
\al_2 \,\le\, y \,\le\, \al_{15}\, ,
\hspace{8.0ex}
\al_3 \,\le\, z \,\le\, \al_{16}
&
\\
&
\al_4 \,\le\, x+y \,\le\, \al_{17}\, ,
\hspace{4ex}
\al_5 \,\le\, x+z \,\le\, \al_{18}\, ,
\hspace{4ex}
\al_6 \,\le\, y+z \,\le\, \al_{19}
&
\\
&
\al_7 \,\le\, x-y \,\le\, \al_{20}\, ,
\hspace{4ex}
\al_8 \,\le\, x-z \,\le\, \al_{21}\, ,
\hspace{4ex}
\al_9 \,\le\, y-z \,\le\, \al_{22}
&
\\
&
\al_{10} \,\le\, x+y+z \,\le\, \al_{23}\, ,
\hspace{6ex}
\al_{11} \,\le\, x-y-z \,\le\, \al_{24}
&
\\
&
\al_{12} \,\le\, x-y+z \,\le\, \al_{25}\, ,
\hspace{6ex}
\al_{13} \,\le\, x+y-z \,\le\, \al_{26}
&
\label{icosi-last}
\end{eqnarray}
It is clear from \eq{icosi-first}\,--\,\eq{icosi-last} 
that an intersection of two \icoses\ is again an \icos.
Along with \Lref{intersection} 
and \Tref{center-cuboid}, this immediately implies the 
following.\pagebreak[3.99]

\begin{corollary}
\label{Z^3-shape}
Let $\cA_d$ be an optimal tristance anticode of diameter $d$ in $\Gplusthree$. 
Then $\cA_d$
is an \icos\ $\sI(\alpha_1,\alpha_2,\ldots,\alpha_{26})$
for some $\alpha_1,\alpha_2,\ldots,\alpha_{26} \in \Z$.\vspace{-1.0ex}
\end{corollary}

As in \S2.3, we can assume that  $\al_1 = \al_2 = \al_3 = 0$
in \eq{icosi-first}, up to a translation in $\Z^3$. We 
furthermore re-parametrize an \icos\
$\sI(0,0,0,\alpha_4,\al_5,\ldots,\alpha_{26})$ as follows
\begin{eqnarray}
\label{icos-first}
0 \le x \le a\, ,
\hspace{4.0ex}
0\hspace*{-3.75ex} 
& \,\le y \le\, & 
\hspace*{-3.75ex}b\, ,
\hspace{4.0ex}
0 \le z \le c
\\[.25ex]
\label{icosi-xy}
e_{\bx\by} 
\hspace*{-1.25ex} & \,\le\, x+y \,\le\, & \hspace*{-1.25ex} 
a + b - e_{\x\y}\hspace*{7ex}
\\[.25ex]
\label{icosi-xy1}
e_{\bx\y} - b 
\hspace*{-1.25ex} & \,\le\, x-y \,\le\, & \hspace*{-1.25ex} 
a - e_{\x\by}
\\[.25ex]
\label{icosi-xz}
e_{\bx\bz} 
\hspace*{-1.25ex} & \,\le\, x+z \,\le\, & \hspace*{-1.25ex} 
a + c - e_{\x\z}
\\[.25ex]
\label{icosi-xz1}
e_{\bx\z} - c 
\hspace*{-1.25ex} & \,\le\, x-z \,\le\, & \hspace*{-1.25ex} 
a - e_{\x\bz}
\\[.25ex]
\label{icosi-yz}
e_{\by\bz} 
\hspace*{-1.25ex} & \,\le\, y+z \,\le\, & \hspace*{-1.25ex} 
b + c - e_{\y\z}
\\[.25ex]
\label{icosi-yz1}
e_{\by\z} - c 
\hspace*{-1.25ex} & \,\le\, y-z \,\le\, & \hspace*{-1.25ex} 
b - e_{\y\bz}
\end{eqnarray}\vspace*{-5ex}
\begin{eqnarray}
\label{icosi-theta1}
\theta_{\bx\by\bz} 
\hspace*{-1.25ex} & \,\le\, x+y+z \,\le\, & \hspace*{-1.25ex} 
a + b + c - \theta_{\x\y\z}
\\
\theta_{\bx\y\bz} - b
\hspace*{-1.25ex} & \,\le\, x-y+z \,\le\, & \hspace*{-1.25ex} 
a + c - \theta_{\x\by\z}
\\
\theta_{\bx\by\z} - c
\hspace*{-1.25ex} & \,\le\, x+y-z \,\le\, & \hspace*{-1.25ex} 
a + b - \theta_{\x\y\bz}
\\
\theta_{\bx\y\z} - b - c 
\hspace*{-1.25ex} & \,\le\, x-y-z \,\le\, & \hspace*{-1.25ex} 
a - \theta_{\x\by\bz}
\label{icos-last}\vspace{-.5ex}
\end{eqnarray}
where 
$a = \al_{14}$, 
$b = \al_{15}$, 
$c = \al_{16}$
while
$e_{\bx\by} = \al_4$,
$e_{\x\y} = \al_{14}+\al_{15}-\al_{17}$,
$e_{\bx\y} = \al_7+\al_{15}$,
$e_{\x\by} = \al_{14}-\al_{20}$,
$e_{\bx\bz} = \al_5$,
$e_{\x\z} = \al_{14}+\al_{16}-\al_{18}$,
$e_{\bx\z} = \al_8+\al_{16}$,
$e_{\x\bz} = \al_{14}-\al_{21}$,
$e_{\by\bz} = \al_6$,
$e_{\y\z} = \al_{15}+\al_{16}-\al_{19}$,
$e_{\by\z} = \al_9+\al_{16}$,
$e_{\y\bz} = \al_{15}-\al_{22}$,
and
$\theta_{\bx\by\bz} = \al_{10}$,
$\theta_{\bx\by\z}  = \al_{13}+\al_{16}$,
$\theta_{\bx\y\bz}  = \al_{12}+\al_{15}$,
$\theta_{\bx\y\z}   = \al_{11}+\al_{15}+\al_{16}$,
$\theta_{\x\by\bz}  = \al_{14}-\al_{24}$,
$\theta_{\x\by\z}   = \al_{14}+\al_{16}-\al_{25}$,
$\theta_{\x\y\bz}   = \al_{14}+\al_{15}-\al_{26}$,
$\theta_{\x\y\z} = \al_{14}+\al_{15}+\al_{16}-\al_{23}$.
\vspace{3.75ex}

\centerline{%
\hspace*{-.3in}\raisebox{-.05in}%
{\smash{\psfig{figure=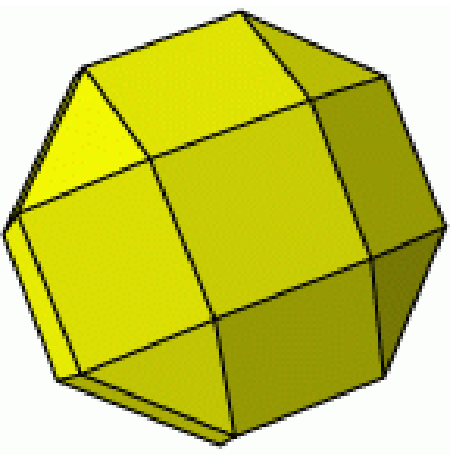,angle=-111,silent=}}} 
\hspace{.50in}
\psfig{figure=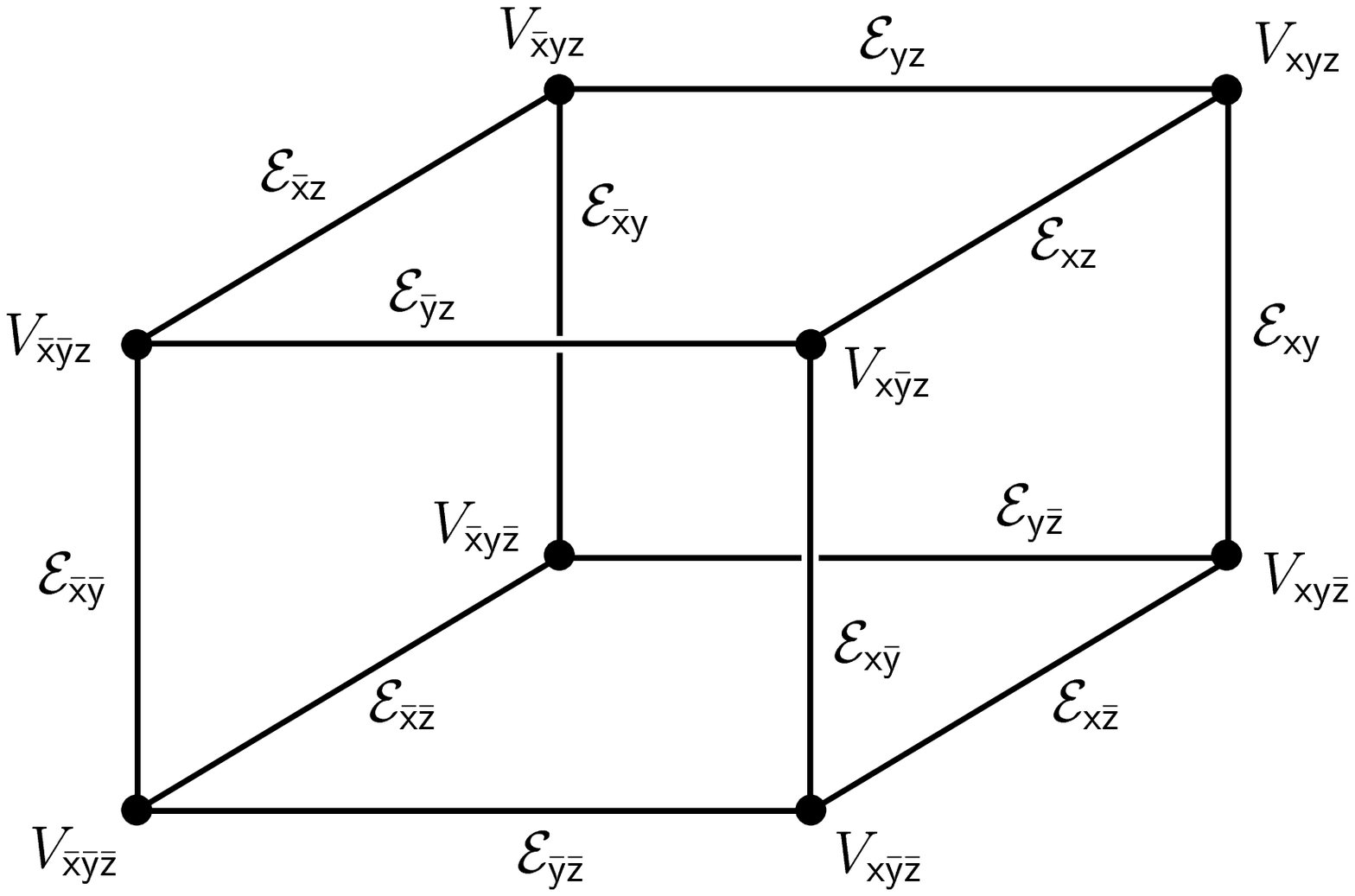,width=3.30in,silent=}} 
\vspace{.250ex}

\begin{center}
{\bf Figure\,9.}
{\sl A generic \icos\ and a labeling of its\/ $20$ truncations}
\vspace{.25ex}
\end{center}

Equations \eq{icos-first}\,--\,\eq{icos-last} make
it apparent that an \icos\ is just a truncated
cuboid: the eight values 
$\theta_{\bx\by\bz},\theta_{\bx\by\z},\ldots,\theta_{\x\y\z}$
give the amount of truncation at the vertices, while the 
twelve values
$e_{\bx\by},e_{\bx\y},\ldots,e_{\y\z}$
describe the amount of truncation along the edges.
Figure\,9 shows a generic \icos\ 
as a $26$-faceted three-dimensional solid 
along with our labeling of the edges and vertices
of the corresponding cuboid, as reflected in 
\eq{icos-first}\,--\,\eq{icos-last}.\pagebreak[3.99]

Observe that one can assume w.l.o.g.\ that each of the $26$ 
inequalities in~\eq{icosi-first}\,--\,\eq{icosi-last} holds
with equality for at least one point of 
$\sI(\alpha_1,\alpha_2,\ldots,\alpha_{26})$
---
otherwise, we can always increase the corresponding constant
$\al_i$ if $i \le 13$ or decrease it if $i \ge 14$. This implies
that each of the inequalities in~\eq{icos-first}\,--\,\eq{icos-last} 
must also hold with equality for at least one point of the \icos.
We will make use of this observation later on. 

\looseness=-1
We next determine the diameter of an \icos\
$\sI(a,b,c,\{e_{\p\q},e_{\p\r},e_{\q\r}\},\{\theta_{\p\q\r}\})$
para\-metrized as in~\eq{icos-first}\,--\,\eq{icos-last}.
To this end, we need to consider certain configurations~of
edges and vertices of a cuboid. Referring to Figure\,9, 
we define 
$\bar{\bx} = \x$,
$\bar{\by} = \y$,
and
\mbox{$\bar{\bz} = \z$},
so that the complement operation $\bar{\cdot}$ is an involution,
as expected.
Let $\p$ denote $\x$ or $\bx$, 
let $\q$~denote $\y$ or $\by$, 
and
let $\r$ denote $\z$ or $\bz$.
With this notation, we say
that a vertex $V_{\p\q\r}$~{\dfn lies opposite\/}
the edges $\cE_{\bp\bq}$, $\cE_{\bp\br}$, and $\cE_{\bq\br}$
(indeed, these are the three edges incident upon the diagonally
opposite vertex $V_{\bp\bq\br}$). We also say that the edges 
$\cE_{\p\q}$, $\cE_{\bp\r}$, $\cE_{\bq\br}$ 
{\dfn span the~cubo\-id\/} (these are the $8$ 
possible choices of three edges such that each face contains 
one of them).

\begin{lemma}
\label{icos-diameter}
Let $\sI(a,b,c,\{e_{\p\q},e_{\p\r},e_{\q\r}\},\{\theta_{\p\q\r}\})$
be an \icos, parametrized~as in~\eq{icos-first}\,--\,\eq{icos-last}.
Define\vspace{-1.00ex}
\begin{eqnarray}
\label{s-def}
s 
& \!{\deff}\! &
\min_{\p,\q,\r} 
\left\{ e_{\p\q} + e_{\bp\r} + e_{\bq\br} \right\}
\\
\label{t-def}
t
& \!{\deff}\! &
\min_{\p,\q,\r} 
\left\{ 
\min\{e_{\bp\bq},e_{\bp\br},e_{\bq\br}\} + \theta_{\p\q\r}
\right\}
\end{eqnarray}
where the minimum in~\eq{s-def} and\/ \eq{t-def} is taken 
over the eight possible assignments of values to $\p$, $\q$, and $\r$.
Then the diameter of 
$\sI(a,b,c,\{e_{\p\q},e_{\p\r},e_{\q\r}\},\{\theta_{\p\q\r}\})$
is given by 
\Eq{
\label{3D-diameter}
d \ = \ a + b + c - \min\{s,t\}
\vspace{-0.75ex}
}
\end{lemma}

\begin{proof}
Let 
$v_1 = (x_1,y_1,z_1)$,
$v_2 = (x_2,y_2,z_2)$,
and
$v_3 = (x_3,y_3,z_3)$
be a triple~of~distinct points of 
$\sI(a,b,c,\{e_{\p\q},e_{\p\r},e_{\q\r}\},\{\theta_{\p\q\r}\})$
and consider their bounding cuboid $\cC(v_1,v_2,v_3)$.
The key observation is that each of the six faces of 
$\cC(v_1,v_2,v_3)$ must contain at least one of the
three points. This leads to the following three cases.\vspace{-1.5ex}
\begin{description}
\item[\hspace*{5.5ex}Case\,1.]
Suppose that none of the points $v_1,v_2,v_3$ is
a vertex of $\cC(v_1,v_2,v_3)$.
Then each of the points $v_1,v_2,v_3$ must belong to an edge 
of $\cC(v_1,v_2,v_3)$~and, moreover,~the three edges
must span the cuboid. Thus if the cuboid $\cC(v_1,v_2,v_3)$
is labeled~as~in~Fi\-gure\,9, we can assume w.l.o.g.\ 
that $v_1 \in \cE_{\p\q}$, $v_2 \in \cE_{\bp\r}$, 
and~$v_3 \in \cE_{\bq\br}$ for some~$\p,\q,\r$.
Referring to Figure\,9, it 
follows that in the definition of $\cC(v_1,v_2,v_3)$, we must have
$
\xmax - \xmin  = |x_1 - x_2|
$,
$
\ymax - \ymin = |y_1 - y_3|
$,
and
$
\zmax - \zmin = |z_2 - z_3|
$.
Thus 
$$
d_3(v_1,v_2,v_3) 
\,=\,
|x_1 - x_2| + |y_1 - y_3| + |z_2 - z_3|
\ = \
\pm (x_1 \pm y_1) \mp (x_2 \pm z_2) \pm (y_3 \pm z_3)
$$
There are eight cases, depending on the $8$ possible 
values of $\p,\q,\r$; but in each case 
$x_1 \pm y_1$ is bounded by \eq{icosi-xy} or \eq{icosi-xy1},
$x_2 \pm z_2$ is bounded by \eq{icosi-xz} or \eq{icosi-xz1},
and
$y_3 \pm z_3$ is bounded~by \eq{icosi-yz} or \eq{icosi-yz1}.
In all cases, these bounds produce the same result, namely
\Eq{
\label{d-edges}
d_3(v_1,v_2,v_3) 
\ \le \
a + b + c \ - \ \left(e_{\p\q} + e_{\bp\r} + e_{\bq\br}\right)
\ \le \
a + b + c \ - \ s
}
Note that we can always achieve the second inequality in~\eq{d-edges}
with equality~by~choosing $v_1 \in \cE_{\p\q}$, $v_2 \in \cE_{\bp\r}$, 
$v_3 \in \cE_{\bq\br}$ so that  $\p,\q,\r$ attains the
minimum~in~\eq{s-def}. The 
\looseness=-1
first inequality in~\eq{d-edges} can  
be also achieved with equality, because of the assumption
that each of \eq{icosi-xy}\,--\,\eq{icosi-yz1} holds with
equality for some point of the \icos.

\item[\hspace*{5.5ex}Case\,2.]
Suppose that one of the three points $v_1,v_2,v_3$ is a vertex 
of $\cC(v_1,v_2,v_3)$,~say $v_1 = V_{\p\q\r}$, but none
of the other two points is in the diagonally opposite
vertex~$V_{\bp\bq\br}$. Then one of $v_2,v_3$ must belong
to an edge that lies opposite $V_{\p\q\r}$, 
say $v_2 \in \cE_{\bp\bq}$,~while 
the other point must belong to the remaining face of $\cC(v_1,v_2,v_3)$.
Referring once again to Figure\,9, we see that
$
\xmax - \xmin  = |x_1 - x_2|
$,
$
\ymax - \ymin = |y_1 - y_2|
$,
and
$
\zmax - \zmin = |z_1 - z_3|
$.
It follows that
$$
d_3(v_1,v_2,v_3) 
\,=\,
|x_1 - x_2| + |y_1 - y_2| + |z_1 - z_3|
\ = \
\pm (x_1 \pm y_1 \pm z_1) \mp (x_2 \pm y_2) \pm z_3
$$
As before, there are eight cases depending on the values
of $\p,\q,\r$, but in each case $x_1 \pm y_1 \pm z_1$
is bounded by one of \eq{icosi-theta1}\,--\,\eq{icos-last}, 
$x_2 \pm y_2$ is bounded by \eq{icosi-xy} or \eq{icosi-xy1},
and
$\pm z_3$ is bounded by~\eq{icos-first}.
In all the eight cases, we get the same result, namely
\Eq{
\label{d-edge/vertex}
d_3(v_1,v_2,v_3) 
\ \le \
a + b + c \ - \ \left(e_{\bp\bq} + \theta_{\p\q\r}\right)
\ \le \
a + b + c \ - \ t
}
Once again, we can attain the first inequality in~\eq{d-edge/vertex}
with equality by choosing~suitable points $v_1,v_2,v_3$ in 
the \icos.
The other two cases where $v_2 \in \cE_{\bp\br}$ or 
\mbox{$v_2 \in \cE_{\bq\br}$}
are similar, leading to the minimization among 
$e_{\bp\bq}, e_{\bp\br}, e_{\bq\br}$ in~\eq{t-def}.

\item[\hspace*{5.5ex}Case\,3.]
Now suppose that one of the points $v_1,v_2,v_3$ 
is a vertex of $\cC(v_1,v_2,v_3)$~and another
of the points is the diagonally opposite vertex, 
say $v_1 = V_{\p\q\r}$ and $v_2 = V_{\bp\bq\br}$. 
Then 
$
\xmax - \xmin = |x_1 - x_2|
$,
$
\ymax - \ymin = |y_1 - y_2|
$,
and
$
\zmax - \zmin = |z_1 - z_2|
$,
so
$$
d_3(v_1,v_2,v_3) 
\,=\,
|x_1 - x_2| + |y_1 - y_2| + |z_1 - z_2|
\ = \
\pm (x_1 \pm y_1 \pm z_1) \mp (x_2 \pm y_2 \pm z_2)
$$
We again get eight cases depending on the values of 
$\p,\q,\r$, with $x_1 \pm y_1 \pm z_1$ and $x_2 \pm y_2 \pm z_2$
both bounded by the same equation 
-- one of \eq{icosi-theta1}\,--\,\eq{icos-last} --
one from above and the other from below. This produces
\Eq{
\label{d-vertex/vertex}
d_3(v_1,v_2,v_3) 
\ \le \
a + b + c \ - \ \left(\theta_{\p\q\r} + \theta_{\bp\bq\br}\right)
}
We can again achieve the bound in ~\eq{d-vertex/vertex} with 
equality but, as we shall see, this case does not produce
a diametric triple of points in 
$\sI(a,b,c,\{e_{\p\q},e_{\p\r},e_{\q\r}\},\{\theta_{\p\q\r}\})$.
\vspace{-1ex}
\end{description}
Since the three cases above are exhaustive, in order to 
complete the proof of the lemma, it would suffice to show
that
\Eq{
\label{theta-constraint}
\theta_{\bp\bq\br}
\ \ge \
\frac{e_{\bp\bq} + e_{\bp\br} + e_{\bq\br}}{2}
\ \ge \
\min\left\{e_{\bp\bq}, e_{\bp\br}, e_{\bq\br}\right\}
}
which, in conjunction with~\eq{d-vertex/vertex}, would imply 
that $d_3(v_1,v_2,v_3) \le a + b + c - t$ in Case\,3.
This follows from the fact that each of the $26$ inequalities
in~\eq{icos-first}\,--\,\eq{icos-last} must hold with equality 
at some point of the \icos.
For example, let $\p = \x$, $\q = \y$, and $\r = \z$. 
Adding the first inequalities of \eq{icosi-xy}, \eq{icosi-xz},
\eq{icosi-yz} yields 
$
2(x+y+z) \ge e_{\bx\by} + e_{\bx\bz} + e_{\by\bz}
$.
Thus if the first inequality in~\eq{icosi-theta1}
is to hold with equality, we must have 
$
\theta_{\bx\by\bz}
\ge 
\shalf{1} \kern-2pt \left(e_{\bx\by} + e_{\bx\bz} + e_{\by\bz}\right)
$.
The other seven ways to assign values to $\p,\q,\r$ can be treated
similarly.~~\qed\vspace{1ex}
\end{proof}

The next task is to determine the volume of 
$\sI(a,b,c,\{e_{\p\q},e_{\p\r},e_{\q\r}\},\{\theta_{\p\q\r}\})$
in terms of its parameters. This  innocuous task is surprisingly
arduous: a complete expression~for
$|\sI(a,b,c,\{e_{\p\q},e_{\p\r},e_{\q\r}\},\{\theta_{\p\q\r}\})|$
would entail hundreds of cases depending upon the rela\-tionships
between various parameters. Moreover, given such an expression,\pagebreak[3.99]
we would need to solve a nonlinear integer optimization problem
involving $23$ variables --- the parameters 
$a,b,c,\{e_{\p\q},e_{\p\r},e_{\q\r}\},\{\theta_{\p\q\r}\}$
in~\eq{icos-first}\,--\,\eq{icos-last}. This problem 
does not appear to be tractable. The situation simplifies 
considerably, however, with the help of the following.

\begin{conjecture}
For each diameter $d \ge 2$, there exists an optimal tristance 
anticode~in~$\Gplusthree$
$\cA_d = \sI(a,b,c,\{e_{\p\q},e_{\p\r},e_{\q\r}\},\{\theta_{\p\q\r}\})$
with equally truncated edges; that is, such that
\Eq{
\label{e-def}
e_{\bx\by} \,=\, 
e_{\bx\y}  \,=\, 
e_{\bx\bz} \,=\, 
e_{\bx\z}  \,=\, 
e_{\by\bz} \,=\, 
e_{\by\bz} \,=\, 
e_{\y\bz}  \,=\, 
e_{\y\z}   \,=\, 
e_{\x\bz}  \,=\, 
e_{\x\z}   \,=\, 
e_{\x\by}  \,=\, 
e_{\x\y}   
\,\ \ \deff\ \
e
}
\end{conjecture}

It is easy to see from~\eq{t-def} that if 
$\cA_d = \sI(a,b,c,\{e_{\p\q},e_{\p\r},e_{\q\r}\},\{\theta_{\p\q\r}\})$
is an optimal~anti\-code with equally truncated edges, then its
vertices must also be equally truncated,~that~is
\Eq{
\label{theta-def}
\theta_{\bx\by\bz} \,=\,
\theta_{\bx\by\z} \,=\,
\theta_{\bx\y\bz} \,=\,
\theta_{\bx\y\z} \,=\,
\theta_{\x\by\bz} \,=\,
\theta_{\x\by\z} \,=\,
\theta_{\x\y\bz} \,=\,
\theta_{\x\y\z} 
\,\ \ \deff\ \
\theta
}
We will denote an \icos\ satisfying~\eq{e-def} and \eq{theta-def}
as $\sI(a,b,c,e,\theta)$. It~now~follows from 
\eq{s-def}\,--\,\eq{3D-diameter} and \eq{theta-constraint}
that if $\cA_d = \sI(a,b,c,e,\theta)$ then
\Eq{
\label{theta-limits}
\frac{3e}{2} 
\ \le \
\theta 
\ \le \
2e
}
The condition that each of the inequalities 
in \eq{icos-first}\,--\,\eq{icos-last} holds
with equality at some point of 
$\sI(a,b,c,e,\theta)$ further implies that 
$
2e \le \min\{a,b,c\}
$.
We present the next lemma~with\-out proof; 
while its proof is not conceptually 
difficult, it is rather tedious.

\begin{lemma}
\label{icos-volume}
Subject to the condition 
$\,\shalf{3} \kern1pt e \le \theta \le 2e \le \min\{a,b,c\}$,
the volume of an \icos\ $\sI(a,b,c,e,\theta)$ 
is given by
\begin{eqnarray*}
|\sI(a,b,c,e,\theta)|
&\hspace{-.5ex}=\hspace{-.5ex}&
(a+1)(b+1)(c+1) \ - \ 2e(e+1)(a+b+c+3)
\nonumber\\[1ex]
& &
+ \ 24e^3 \, + \ \frac{4}{3}\, \theta
\Bigl(
\kern1pt 3\theta(6e{-}1) - 9e(3e{-}1) - (2\theta {+} 1)(2\theta {-} 1)\,
\Bigr)
\nonumber
\end{eqnarray*}
\end{lemma}

Using the expression for $|\sI(a,b,c,e,\theta)|$ in \Lref{icos-volume}
along with~\eq{3D-diameter}, 
it can be furthermore shown that 
if $A_d = \sI(a,b,c,e,\theta)$ is an optimal tristance 
anticode in $\Gplusthree$, then $\theta = 2e$.
With this, the expression for the volume of the \icos\
further simplifies~to
\Eq{
\label{final-volume}
|\sI(a,b,c,e,2e)|
\ = \
(a\,{+}\,1)(b\,{+}\,1)(c\,{+}\,1) \ - \ 
2e(e\,{+}\,1)
\left( 
a\,{+}\,b\,{+}\,c\,{+}\,3 - \frac{4}{3}\,(2e\,{+}\,1)
\right)
}
It remains to maximize the cubic on the right-hand side
of~\eq{final-volume} subject to the constraints
$a+b+c-3e = d$ and $a \ge b \ge c \ge 2e$.
Note that for each fixed $e$, we have
$$
|\sI(a,b,c,e,2e)|
\ = \
(a\,{+}\,1)(b\,{+}\,1)(c\,{+}\,1) \ - \ 
\text{const}
$$
since $a+b+c = d+3e$. This immediately shows that the optimal values
of $a,b,c$~are~given by $c = \floorenv{d/3} + e$ with 
$a,b$ being equal to either $c$ or $c+1$. The complete
solution to the optimization problem is given in Table\,4,
where $[\mu]$ denotes the integer that is closest to the 
real number $\mu$ (rounding). 
We have verified by exhaustive computer search that~the 
anticodes in Table\,4 are, in fact, the 
\emph{unique optimal tristance anticodes in}
$\Gplusthree$~up~to~dia\-meter $d = 11$. 
Figure\,10 shows some of these
anticodes, for diameters $d = 9,10,11$.\pagebreak[3.99]

$$
\begin{array}{|c|c|c|c|c|}
\hline \hline
& & & & \\[-1.50ex]
\,d \hspace{-1.50ex}\pmod{3} & a & b & c & e 
\\[0.75ex]
\hline\hline
& & & & \\[-1.250ex]
0  & 
\frac{d}{3} + e & \frac{d}{3} + e & \frac{d}{3} + e & 
{\scriptstyle (d\,{+}\,1) \ - \ \left[\sqrt{\sthird{2}(d{+}1)(d{+}2)}\,\right]}
\\[1.50ex]
\hline
& & & & \\[-1.250ex]
1  & 
\frac{d+2}{3} + e & \frac{d-1}{3} + e & \frac{d-1}{3} + e & 
{\scriptstyle (d\,{+}\,1) \ - \ 
\left[\sqrt{\sthird{2}(d{+}1)(d{+}2) \:+\, \sthird{1}}\; \right]}
\\[1.50ex]
\hline
& & & & \\[-1.250ex]
2  & 
\frac{d+1}{3} + e & \frac{d+1}{3} + e & \frac{d-2}{3} + e & 
{\scriptstyle (d\,{+}\,1) \ - \ 
\left[\sqrt{\sthird{2}(d{+}1)(d{+}2) \:+\, \sthird{1}}\; \right]}
\\[2.0ex]
\hline\hline
\end{array}
$$

\begin{center}
{\bf Table\,4.}
{\sl Parameters of 
(conjecturally) optimal tristance anticodes in $\Gplusthree$}
\vspace{5.0ex}
\end{center}

\centerline{
\psfig{figure=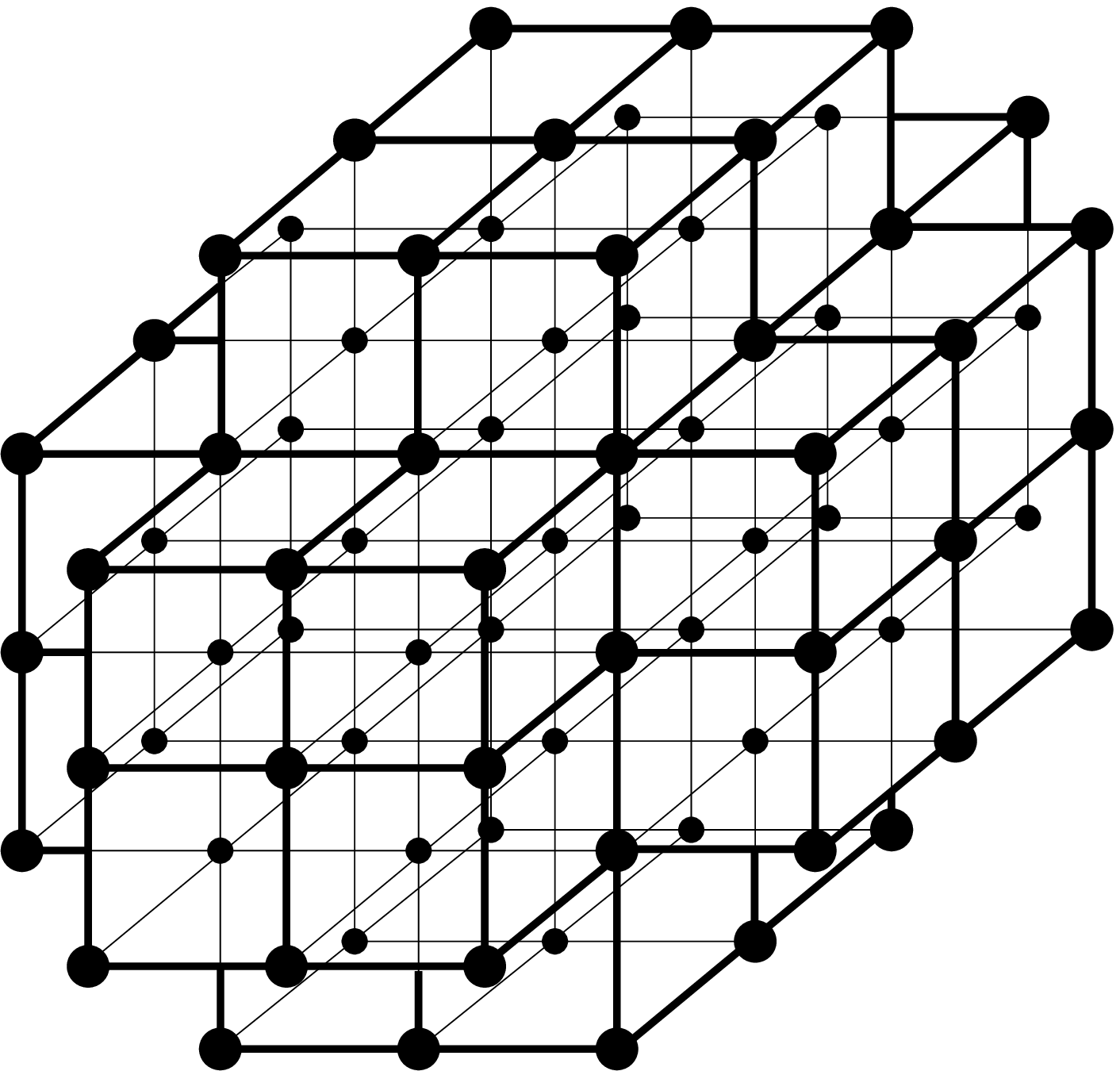,width=2.60in,silent=}%
\hspace{1.0in}%
\psfig{figure=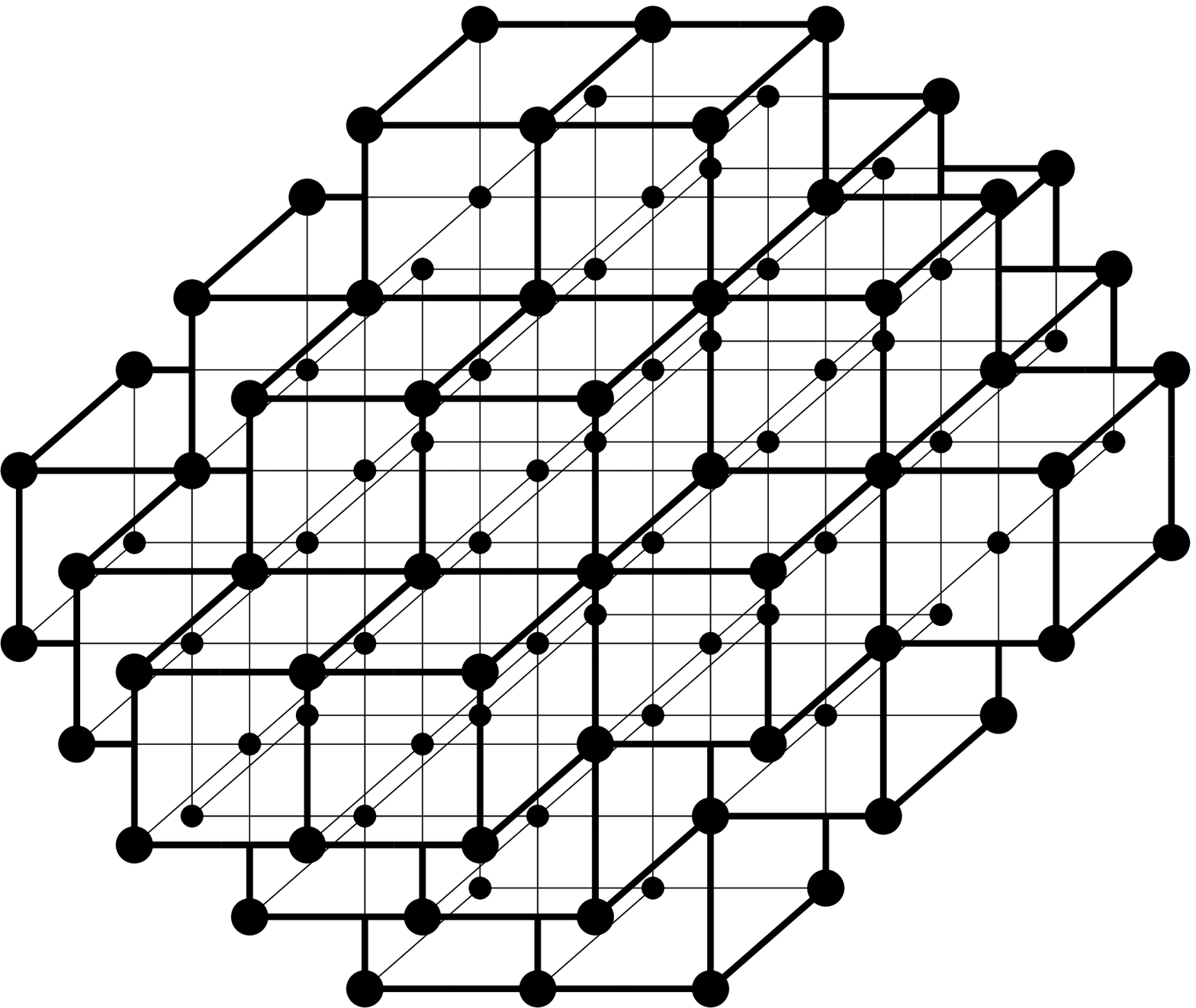,width=3.15in,silent=}}

\hspace*{.6in}\psfig{figure=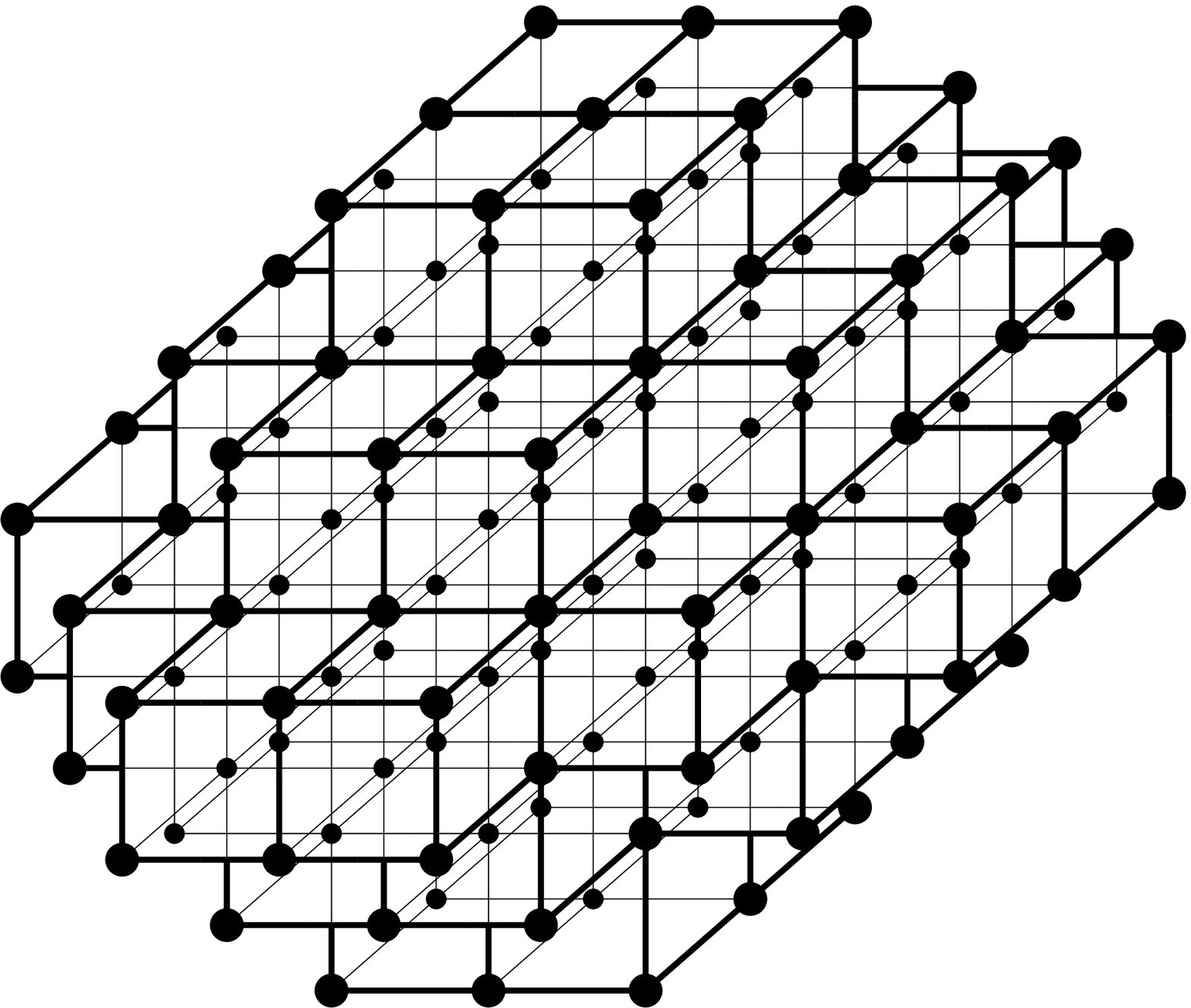,width=3.65in,silent=}
\vspace{2ex}

\begin{center}
{\bf Figure\,10.}
{\sl Optimal tristance anticodes in $\Gplusthree$ of diameter $d = 9,10,11$}
\end{center}

\subsection{Optimal quadristance anticodes in the grid graph of $\Z^2$}
\vspace{-1.00ex}
\label{sec4.2}

Recall that, given four distinct points $z_1,z_2,z_3,z_4$ in $\Z^2\!$,
the {\dfn quadristance\/} $d_4(z_1,z_2,z_3,z_4)$ is defined as 
the number of edges in a minimal spanning tree for $z_1,z_2,z_3,z_4$ 
in the grid graph $\Gplus$ of $\Z^2$. The following expression for
quadristance is implicit in~\cite{EV}.

\begin{theorem}
\label{quadristance-formula}
Let  $z_1 = (x_1,y_1)$, $z_2 = (x_2,y_2)$,
$z_3 = (x_3 ,y_3 )$, $z_4 = (x_4,y_4)$~be~\mbox{distinct~points}
in $\Z^2$. Let $\sigma$ and $\tau$ be permutations of\/ $\{1,2,3,4\}$
such that 
$x_{\sigma(1)} \le x_{\sigma(2)} \le x_{\sigma(3)} \le x_{\sigma(4)}$
and
$y_{\tau(1)} \le y_{\tau(2)} \le y_{\tau(3)} \le y_{\tau(4)}$.
Then
\Eq{
\label{quad-simple}
d_4(z_1,z_2,z_3,z_4)
\ = \
\Bigl(x_{\sigma(4)} - x_{\sigma(1)}\Bigr)
\ + \ 
\Bigl(y_{\tau(4)} - y_{\tau(1)}\Bigr)
}
provided $\tau \sigma^{-1}\! \in \Gamma$, where $\Gamma$ 
is the subgroup of the symmetric group 
generated by~the~permutations $(1,2)$,~$(3,4)$, 
and $(1,3)(2,4)$. If $\tau \sigma^{-1}\! \not\inn \Gamma$ then
\Eq{
\label{quad-plus}
\hspace*{-1ex}d_4(z_1,z_2,z_3,z_4)
\: = \:
\Bigl(x_{\sigma(4)} \,-\, x_{\sigma(1)}\Bigr)
\ + \ 
\Bigl(y_{\tau(4)} \,-\, y_{\tau(1)}\Bigr)
\ + \
\min\!\left\{x_{\sigma(3)} - x_{\sigma(2)},y_{\tau(3)} - y_{\tau(2)}\right\}
}
\end{theorem}

Note that one can assume w.l.o.g.\ that $\sigma$ is the
identity permutation. Then $d_4(z_1,z_2,z_3,z_4)$ is given 
by~\eq{quad-simple} precisely in the eight cases where
\begin{eqnarray*}
&
y_1 \le y_2 \le y_3 \le y_4 \,,\hspace{3.00ex}   
y_1 \le y_2 \le y_4 \le y_3 \,,\hspace{3.00ex}     
y_2 \le y_1 \le y_3 \le y_4 \,,\hspace{3.00ex}     
y_2 \le y_1 \le y_4 \le y_3   
&
\\[.5ex]
&
y_3 \le y_4 \le y_2 \le y_1 \,,\hspace{3.00ex}   
y_3 \le y_4 \le y_1 \le y_2 \,,\hspace{3.00ex}     
y_4 \le y_3 \le y_1 \le y_2 \,,\hspace{3.00ex}     
y_4 \le y_3 \le y_2 \le y_1   
&
\end{eqnarray*}

We next determine the optimal quadristance anticode 
$\cA_d(z_1,z_2,z_3) \subset \Z^2$ 
centered about three given 
points $z_1,z_2,z_3$, namely the set
\Eq{
\cA_d(z_1,z_2,z_3)
\ \ \deff\,\ 
\left\{ z \in \Z^2 ~:~ d_4(z_1,z_2,z_3,z) \le d \right\}
}
Clearly $\cA_d(z_1,z_2,z_3) = \varnothing$ for $d < d_3(z_1,z_2,z_3)$.
As before, we first consider the case where $d = d_3(z_1,z_2,z_3)$.
Recall that $\cR(z_1,z_2)$ denotes the bounding rectangle of 
$z_1,z_2 \in \Z^2$.

\begin{lemma}
\label{A-Delta}
Let  $z_1 = (x_1,y_1)$, $z_2 = (x_2,y_2)$, $z_3 = (x_3,y_3)$ 
be distinct points in $\Z^2$,
and~let $\Delta = d_3(z_1,z_2,z_3)$. Write 
$\cR_1 = \cR(z_1,z_2)$, 
$\cR_2 = \cR(z_1,z_3)$, 
and
$\cR_3 = \cR(z_2,z_3)$.
Then
\Eq{
\label{A_Delta}
\cA_{\Delta}(z_1,z_2,z_3)
\: = \:
\Bigl( \cR_1 \cap \cR_2 \Bigr)
\,\cup\,
\Bigl( \cR_1 \cap \cR_3 \Bigr)
\,\cup\,
\Bigl( \cR_2 \cap \cR_3 \Bigr)\vspace{1ex}
}
\end{lemma}

\begin{proof}
By definition, $z \in \cA_{\Delta}(z_1,z_2,z_3)$ iff 
$d_4(z_1,z_2,z_3,z) = d_3(z_1,z_2,z_3)$.~This
happens if and only if $z$ belongs to the vertex set of 
a minimal spanning tree for $z_1,z_2,z_3$. Hence
$$
\cA_{\Delta}(z_1,z_2,z_3)
\ = \hspace{-1ex}
\bigcup_{T(z_1,z_2,z_3)} \hspace{-2.00ex}
\bigl\{\text{vertex set of $T(z_1,z_2,z_3)$}\bigr\}
$$
where the union is over all the minimal spanning trees for $z_1,z_2,z_3$. 
Observe that given any $u,v \in \Z^2\!$, the union of all the shortest paths
(minimal spanning trees) between~$u$~and~$v$ in~$\Gplus$ is precisely
the bounding rectangle $\cR(u,v)$. 
Now consider $\cR(z_1,z_2,z_3)$, the bounding rectangle of $z_1,z_2,z_3$.
Since each of the four edges of $\cR(z_1,z_2,z_3)$ must 
contain 
at least one of the three points, at least one of $z_1,z_2,z_3$ must be 
a vertex of $\cR(z_1,z_2,z_3)$. Thus we can assume w.l.o.g.\ that
$z_1$ is a vertex of $\cR(z_1,z_2,z_3)$. 
We distinguish between 
two cases.\vspace{-1.5ex}
\begin{description}
\item[\hspace*{5.5ex}Case\,1.]
Suppose that one of the other two points, say $z_3$, is the 
opposite vertex of~$\cR(z_1,z_2,z_3)$, as 
illustrated on the right side of Figure\,11. It is easy to see 
that,~in this case, any minimal spanning tree for $z_1,z_2,z_3$ is
a union of a shortest path from $z_1$ to $z_2$ with 
a shortest path from $z_2$ to $z_3$.
Hence
\Eq{
\label{A_Delta-aux1}
\cA_{\Delta}(z_1,z_2,z_3)
\ = \ 
\cR(z_1,z_2) \cup \cR(z_2,z_3) 
}
Since $\cR_2 = \cR(z_1,z_3) = \cR(z_1,z_2,z_3)$ in this case,
we have $\cR_1 \cap \cR_2 = \cR(z_1,z_2)$ and 
$\cR_2 \cap \cR_3 = \cR(z_2,z_3)$. It follows that~\eq{A_Delta-aux1}
coincides with \eq{A_Delta}.

\item[\hspace*{5.5ex}Case\,2.]
Suppose that none of the other two points is a vertex 
of $\cR(z_1,z_2,z_3)$ opposite to $z_1$. Then $z_2$ and $z_3$
must belong to the two edges of $\cR(z_1,z_2,z_3)$ that lie
opposite $z_1$. This case is illustrated on the left side of Figure\,11.
As in \eq{d3-hex}, let~$\xmid$ and $\ymid$ denote the middle values 
among $x_1,x_2,x_3$ and $y_1,y_2,y_3$, respectively. Write
$v = (\xmid,\ymid)$. Then any minimal spanning tree for $z_1,z_2,z_3$ 
consists of a shortest path from $z_1$ to $v$ along with the 
unique shortest path from $v$ to $z_2$ and the 
unique shortest path from $v$ to $z_3$. Thus
\Eq{
\label{A_Delta-aux2}
\cA_{\Delta}(z_1,z_2,z_3)
\ = \ 
\cR(z_1,v) \cup \cR(z_2,v) \cup \cR(z_3,v) 
}
Notice that 
$\cR(z_1,v) = \cR_1 \cap \cR_2$,
$\cR(z_2,v) = \cR_1 \cap \cR_3$,
and
$\cR(z_3,v) = \cR_2 \cap \cR_3$~(even 
though the two rectangles $\cR(z_2,v), \cR(z_3,v)$ are degenerate).
Hence, 
\eq{A_Delta-aux2} again~co\-incides 
with \eq{A_Delta}, and we are done.~~\qed\vspace{1ex}
\end{description}
\end{proof}

\begin{theorem}
\label{quad-A_d}
Let  $z_1 = (x_1,y_1)$, $z_2 = (x_2,y_2)$, $z_3 = (x_3,y_3)$ 
be distinct points in $\Z^2$.
Let~$d$ be an integer such that $d \ge d_3(z_1,z_2,z_3)$,
and write $c = d - d_3(z_1,z_2,z_3)$.
Further define
$$
\begin{array}{r@{\hspace{1.00ex}}c@{\hspace{1ex}}l@{\hspace{5ex}}%
                                 r@{\hspace{1.00ex}}c@{\hspace{1ex}}l}
\al_1 
& = &
\max\bigl\{ \min\{x_1,x_2\}, \min\{x_1,x_3\} \bigr\}\,,
&
\beta_1 
& = &
\min\bigl\{ \max\{x_1,x_2\}, \max\{x_1,x_3\} \bigr\}
\\[.75ex]
\gamma_1 
& = &
\max\bigl\{ \min\{y_1,y_2\}, \min\{y_1,y_3\} \bigr\}\,,
&
\delta_1 
& = &
\min\bigl\{ \max\{y_1,y_2\}, \max\{y_1,y_3\} \bigr\}
\\[1.5ex]
\al_2
& = &
\max\bigl\{ \min\{x_1,x_2\}, \min\{x_2,x_3\} \bigr\}\,,
&
\beta_2
& = &
\min\bigl\{ \max\{x_1,x_2\}, \max\{x_2,x_3\} \bigr\}
\\[.75ex]
\gamma_2 
& = &
\max\bigl\{ \min\{y_1,y_2\}, \min\{y_2,y_3\} \bigr\}\,,
&
\delta_2
& = &
\min\bigl\{ \max\{y_1,y_2\}, \max\{y_2,y_3\} \bigr\}
\\[1.5ex]
\al_3
& = &
\max\bigl\{ \min\{x_1,x_3\}, \min\{x_2,x_3\} \bigr\}\,,
&
\beta_3
& = &
\min\bigl\{ \max\{x_1,x_3\}, \max\{x_2,x_3\} \bigr\}
\\[.75ex]
\gamma_3 
& = &
\max\bigl\{ \min\{y_1,y_3\}, \min\{y_2,y_3\} \bigr\}\,,
&
\delta_3
& = &
\min\bigl\{ \max\{y_1,y_3\}, \max\{y_2,y_3\} \bigr\}
\end{array}
$$
Then the centered
quadristance anticode $\cA_d(z_1,z_2,z_3)$ is a union 
of three octagons $\sO_1$,~$\sO_2$, and $\sO_3$, where for 
$i = 1,2,3$, the octagon $\sO_i$ consists of all $(x,y) \in \Z^2$ 
such that \vspace{-.50ex}
\begin{eqnarray*}
\label{quad-centeroct1}
\alpha_i - c \ \le\ x \ \le\ \beta_i + c\, ,
&\hspace{1ex}& 
\alpha_i + \gamma_i - c \ \le\ x+y \ \le\ \beta_i + \delta_i + c
\\[1.25ex]
\label{quad-centeroct2}
\gamma_i - c \ \le\ y \ \le\ \delta_i + c\, ,
&\hspace{1ex}& 
\alpha_i - \delta_i - c \ \le\ x-y \ \le\ \beta_i - \gamma_i + c
\end{eqnarray*} 
\end{theorem}

\begin{proof}
As in \Lref{A-Delta}, let $\Delta = d_3(z_1,z_2,z_3)$.
It is not difficult to show that,~once~again,
$z \in \cA_d(z_1,z_2,z_3)$
if and only if the $L_1$-distance from $z$ to (the closest point of)
$\cA_{\Delta}(z_1,z_2,z_3)$ is at most $c = d - \Delta$.
\Lref{A-Delta} proves that $\cA_{\Delta}(z_1,z_2,z_3)$ is 
a union of three rectangles.
For each rectangle, the set of
all $z \in \Z^2$ that are at $L_1$-distance at most $c$ from
it~is~an~octagon. Indeed, 
$\sO_1$ is precisely the set of all points that are 
at $L_1$-distance at most $c$ from $\cR_1 \cap \cR_2$, 
while $\sO_2$ and $\sO_3$ are constructed similarly
with respect to $\cR_1 \cap \cR_3$ and $\cR_2 \cap \cR_3$,
where $\cR_1, \cR_2, \cR_3$ are as defined in \Lref{A-Delta}.
Hence $\cA_d(z_1,z_2,z_3) = \sO_1 \cup \sO_2 \cup \sO_3$.~~\qed\pagebreak[3.99]
\end{proof}

\centerline{%
\psfig{figure=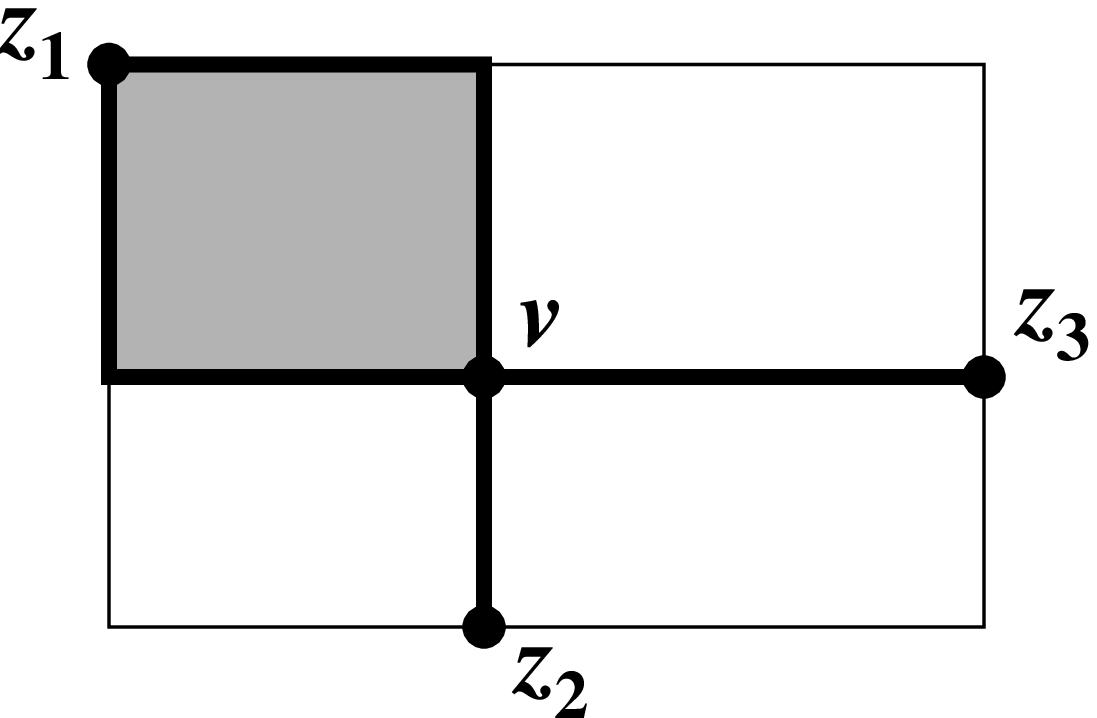,width=2.250in,silent=}
\hspace{.750in}
\raisebox{0.10in}%
{\psfig{figure=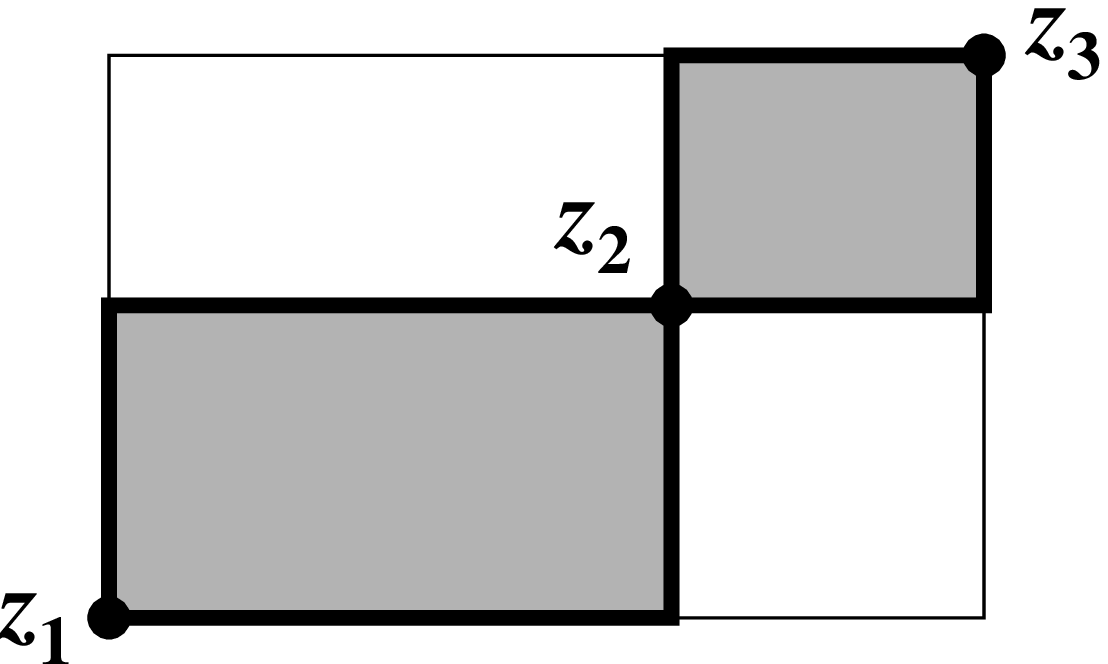,width=2.250in,silent=}}}
\vspace{-.5ex}

\begin{center}
{\bf Figure\,11.}
{\sl 
Centered quadristance anticodes $\cA_d(z_1,z_2,z_3)$
for $d = d_3(z_1,z_2,z_3)$}
\vspace{.25ex}
\end{center}

Now let $\cA_d$ denote an optimal (unrestricted) quadristance
anticode of diameter $d$ in $\Gplus$. Arguing as in \Lref{intersection},
it is easy to show that
\Eq{
\label{quad-intersection}
\cA_d
\ = \hspace{-1ex}
\bigcap_{z_1,z_2,z_3 \in \cA_d} \hspace{-2ex} \cA_d(z_1,z_2,z_3)
}
However, as can be seen from Figure\,11, the sets 
$\cA_d(z_1,z_2,z_3)$ are no longer convex and their
general shape (union of three octagons) is not preserved
under intersection. Thus \eq{quad-intersection} and \Tref{quad-A_d}
do not suffice to determine the shape of $\cA_d$. 
In fact, as we shall see in Figure\,12, $\cA_d$ may come in 
several different shapes, at least for certain diameters.

Nevertheless, we can use an \emph{arbitrary} shape in order 
to derive a \emph{lower bound}~on~the~cardinality of $\cA_d$.
Based on the available numerical evidence (cf.~Figure\,12), 
we will use an octagon with equally truncated corners, namely
the set of all $(x,y) \in \Z^2$ such that
\begin{eqnarray}
\label{quad-oct1}
0 \, \le\, x \, \le\: a\, ,
&\hspace{1ex}&  
c \ \le\ x+y \ \le\  a + b - c\\[0.750ex]
\label{quad-oct2}
0 \, \le\, y \, \le\: b\, ,\hspace{.5ex}
&& 
c - b \ \le\ x-y \ \le \ a - c
\end{eqnarray}
We will denote such a set by $\sO(a,b,c)$. 
We assume that each of the eight inequalities~in
\eq{quad-oct1} and~\eq{quad-oct2} holds with equality 
for some point of $\sO(a,b,c)$; otherwise, we can
always re-parametrize accordingly. 
This, in particular, implies that 
\Eq{
\label{c-cond}
2c \ \le\ \min\{a,b\}
}
Note that the cardinality 
of $\sO(a,b,c)$ is given by 
\eq{A_d-size} with $c_0 = c_1 = c_2 = c_3 = c$.
Thus the~next step is to determine the quadristance diameter of $\sO(a,b,c)$.

\begin{lemma}
\label{quadristance-d}
Let $\sO(a,b,c)$ be the octagon 
in\/ {\rm \eq{quad-oct1}} and\/ {\rm \eq{quad-oct2}} and assume w.l.o.g.\
that $a \ge b$. Then 
the quadristance diameter of $\sO(a,b,c)$ is given by
\Eq{
d \ = \ a + 2b - 2c
}
\end{lemma}

\begin{proof}
\looseness=-1
Let 
$z_1\kern-0.6pt =\kern-0.4pt (x_1,y_1)$, 
$z_2\kern-0.6pt =\kern-0.4pt (x_2,y_2)$, 
$z_3\kern-0.6pt =\kern-0.4pt (x_3,y_3)$, 
$z_4\kern-0.6pt =\kern-0.4pt (x_4,y_4)$
be four arbitrary~points
of $\sO(a,b,c)$,
and assume w.l.o.g.\ that
$x_1 \le x_2 \le x_3 \le x_4$. 
Let $\tau$ be a permutation such that 
$y_{\tau(1)} \le y_{\tau(2)} \le y_{\tau(3)} \le y_{\tau(4)}$.
If $\tau \in \Gamma$ so that $d_4(z_1,z_2,z_3,z_4)$ is given
by~\eq{quad-simple},~then
$$
d_4(z_1,z_2,z_3,z_4)
\ = \
\bigl(x_4 - x_1\bigr)
\, + \, 
\bigl(y_{\tau(4)} - y_{\tau(1)}\bigr)
\ \le \
a + b
\ \le \
a + 2b - 2c
$$
where the first inequality follows from \eq{quad-oct1} and~\eq{quad-oct2}
while the second inequality follows from \eq{c-cond}. Otherwise, 
$d_4(z_1,z_2,z_3,z_4)$ is given by~\eq{quad-plus} so that
\begin{eqnarray}
\hspace*{-2ex}d_4(z_1,z_2,z_3,z_4)\hspace{-.75ex}
& = &\hspace{-.75ex}
\bigl(x_4 - x_1\bigr)
\ + \ 
\bigl(y_{\tau(4)} \,-\, y_{\tau(1)}\bigr)
\ + \
\min\!\left\{x_3 - x_2,y_{\tau(3)} - y_{\tau(2)}\right\}
\\
\label{quad-aux2}
& \le &\hspace{-.75ex}
\bigl(x_4 - x_1\bigr)
\ + \ 
\Bigl(y_{\tau(4)} + y_{\tau(3)} - y_{\tau(2)} - y_{\tau(1)}\Bigr)
\\
\label{quad-aux3}
& = &\hspace{-.75ex}
\bigl(x_4 \pm y_4\bigr) 
\ - \
\bigl(x_1 \pm y_1\bigr) 
\ + \
\bigl(\pm y_2 \pm y_3\bigr) 
\end{eqnarray}
There are four simple cases depending on the signs of $y_1$ 
and $y_4$ in~\eq{quad-aux3}. Observe that, in view of~\eq{quad-aux2},
exactly two of $y_1,y_2,y_3,y_4$ contribute to $d_4(z_1,z_2,z_3,z_4)$
with a positive sign and two with a negative sign.
This immediately implies the following.

\noindent\hspace*{3ex}{\bf Case\,1:} 
$
d_4(z_1,z_2,z_3,z_4)
\: = \:
(x_4 - y_4) \,-\, (x_1+y_1) \,+\, (y_2 + y_3)
\: \le \:
(a-c) \,-\, c \,+\, (b+b)
$\vspace{-.50ex}

\noindent\hspace*{3ex}{\bf Case\,2:} 
$
d_4(z_1,z_2,z_3,z_4)
\: = \:
(x_4 + y_4) \,-\, (x_1-y_1) \,-\, (y_2 + y_3)
\: \le \:
(a+b-c) \,-\, (c-b)
$

\noindent 
In the other two cases, $y_1, y_4$ contribute to 
$d_4(z_1,z_2,z_3,z_4)$ with opposite signs. Hence~$y_2, y_3$ 
also have opposite signs,
so that the last term in~\eq{quad-aux3} is at most $|y_2 - y_3| \le b$.
Thus 

\noindent\hspace*{3ex}{\bf Case\,3:} 
$
d_4(z_1,z_2,z_3,z_4)
\: \le \:
(x_4 - y_4) \,-\, (x_1-y_1) \,+\, b
\: \le \:
(a-c) \,-\, (c - b) \,+\,b
$\vspace{-.50ex}

\noindent\hspace*{3ex}{\bf Case\,4:} 
$
d_4(z_1,z_2,z_3,z_4)
\: = \:
(x_4 + y_4) \,-\, (x_1+y_1) \,+\, b
\: \le \:
(a+b-c) \,-\, c \,+\,b
$

The above shows that 
$
d_4(z_1,z_2,z_3,z_4) \le a + 2b - 2c
$
for any 
$z_1,z_2,z_3,z_4 \in \sO(a,b,c)$.
To see that this bound holds with equality,
consider the points
$z_1 = (0,b\,{-}\,c)$,
$z_2 = (c,0)$,
$z_3 = (a\,{-}\,c,b)$,
$z_4 = (a,c)$.
Then 
$
d_4(z_1,z_2,z_3,z_4) = a + 2b - 2c
$
by~\eq{c-cond} and \Tref{quadristance-formula}.~~\qed\vspace{1ex}
\end{proof}

\begin{theorem}
\label{optimal-quadristance}
Let $\cA_d$ be an optimal quadristance anticode of diameter~$d$
in 
$\Gplus$. Then
$$
|\cA_d|
\ \ge \
\ceilenv{\frac{d^2 + 4d + 3}{6}}
\ = \
\ceilenv{\frac{(d+1)(d+3)}{6}}
\vspace{.75ex}
$$
\end{theorem}

\begin{proof}
The lower bound follows by considering quadristance 
anticodes of type $\sO(a,b,c)$. In view of~\eq{A_d-size}, \eq{c-cond},
and \Lref{quadristance-d}, the optimal parameters $a$, $b$, $c$
are obtained~by~maximizing
$
|\sO(a,b,c)|
 = 
(a+1)(b+1) - 2c(c+1)
$
subject to the constraints $a + 2b - 2c = d$ and $a \ge b \ge 2c$.
The solution to this optimization problem is compiled in 
Table\,5.~~\qed\vspace{1ex}
\end{proof}

\looseness=-1
We have also used exhaustive computer search to find optimal
quadristance anticodes~in~$\Gplus$ of diameters up to $d=9$.
The results of this search are presented in Figure\,12, which 
shows \underline{all} the optimal anticodes (up to obvious isomorphisms)
for $d = 4,5,\ldots,9$. For $d=3$, the\linebreak[3.99]
the optimal quadristance anticodes are simply the five 
tetromino shapes of~\cite{Golomb}.\vspace{1ex}

{\bf Remark.}
We observe that --- for the first time in this paper ---
optimal anticodes~of~a~given diameter do not have a unique shape.
Moreover, for $d = 3,5,9$ non-convex shapes
occur among optimal quadristance anticodes, again for the
first time in this~paper. 

Nevertheless, an octagon with equally truncated corners 
is always among the optimal shapes in Figure\,12. Hence 
the lower bound on $|\cA_d|$ 
in \Tref{optimal-quadristance} is exact at least up to~$d = 9$. 
We conjecture that this bound is, in fact, exact for all 
diameters.

$$
\begin{array}{|c|c|c|c|c|}
\hline \hline
& & & & \\[-1.50ex]
\,d \hspace{-1.50ex}\pmod{6} & a & b & c & |\sO(a,b,c)| 
\\[0.75ex]
\hline\hline
& & & & \\[-1.50ex]
0  & \frac{2d+3}{3} & \frac{d+3}{3} & \frac{d}{6} & \frac{d^2+4d+6}{6} 
\\[0.55ex]
\hline
& & & & \\[-1.95ex]
1  & \frac{2d+4}{3} & \frac{d+2}{3} & \frac{d-1}{6} & \frac{d^2+4d+7}{6} 
\\[0.55ex]
\hline
& & & & \\[-1.95ex]
2  & \frac{2d+5}{3} & \frac{d+1}{3} & \frac{d-2}{6} & \frac{d^2+4d+6}{6} 
\\[0.55ex]
\hline
& & & & \\[-1.50ex]
3  & 
\begin{array}{c}
\frac{2d+6}{3}\\[.75ex]
\frac{2d}{3}
\end{array}
& 
\begin{array}{c}
\frac{d}{3}\\[.75ex]
\frac{d+3}{3}
\end{array}
&
\begin{array}{c}
\frac{d-3}{6}\\[.75ex]
\frac{d-3}{6}
\end{array}
& 
\frac{d^2+4d+3}{6}
\\[2.25ex]
\hline
& & & & \\[-1.95ex]
4  & \frac{2d+1}{3} & \frac{d+2}{3} & \frac{d-4}{6} & \frac{d^2+4d+4}{6} 
\\[0.55ex]
\hline
& & & & \\[-1.50ex]
5  & 
\begin{array}{c}
\frac{2d+2}{3}\\[.75ex]
\frac{2d+2}{3}
\end{array}
& 
\begin{array}{c}
\frac{d+1}{3}\\[.75ex]
\frac{d+4}{3}
\end{array}
&
\begin{array}{c}
\frac{d-5}{6}\\[.75ex]
\frac{d+1}{6}
\end{array}
& 
\frac{d^2+4d+3}{6}
\\[3.15ex]
\hline\hline
\end{array}
$$

\begin{center}
{\bf Table\,5.}
{\sl Parameters of (conjecturally) optimal quadristance anticodes in $\Gplus$}
\vspace{.5ex}
\end{center}

\vspace{0.20in}

\newlength{\mylen}
\setlength{\mylen}{1.0em}

\hspace*{0.20in}\begin{minipage}[t]{2.50in}
\centerline{
\psfig{figure=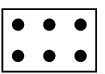,height=2.3\mylen,silent=}
\hspace{.20in}
\psfig{figure=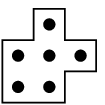,height=3.3\mylen,silent=}
}
\vspace{0.12in}

\centerline{diameter $d = 4$}
\vspace{0.35in}

\centerline{
\psfig{figure=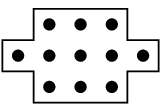,height=3.3\mylen,silent=}
}
\vspace{0.10in}

\centerline{diameter $d = 6$}
\vspace{0.35in}

\centerline{
\psfig{figure=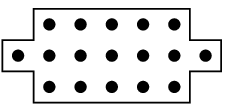,height=3.3\mylen,silent=}
\hspace{.20in}
\psfig{figure=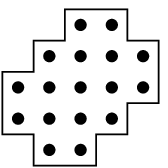,height=5.3\mylen,silent=}
}
\vspace{0.10in}

\centerline{diameter $d = 8$}

\end{minipage}
\hspace{-0.30in}
\begin{minipage}[t]{3.50in}
\centerline{
\psfig{figure=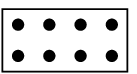,height=2.3\mylen,silent=}
\hspace{.20in}
\psfig{figure=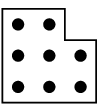,height=3.3\mylen,silent=}}
\vspace{0.15in}

\centerline{
\psfig{figure=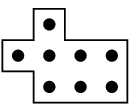,height=3.3\mylen,silent=}
\hspace{.20in}
\psfig{figure=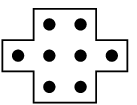,height=3.3\mylen,silent=}
}
\vspace{0.12in}

\centerline{diameter $d = 5$}
\vspace{0.50in}

\centerline{
\psfig{figure=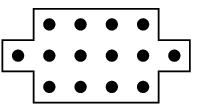,height=3.3\mylen,silent=}
}
\vspace{0.10in}

\centerline{diameter $d = 7$}
\vspace{0.60in}

\hspace*{-.50in}%
\psfig{figure=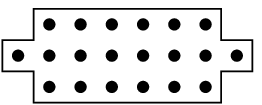,height=3.3\mylen,silent=}
\hspace{.20in}
\psfig{figure=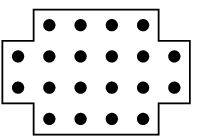,height=4.3\mylen,silent=}
\hspace{.20in}
\psfig{figure=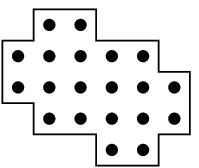,height=5.3\mylen,silent=}\vspace{0.10in}

\centerline{diameter $d = 9$}
\vspace{0.10in}

\end{minipage}

\begin{center}
{\bf Figure\,12.}
{\sl Optimal quadristance anticodes in $\Gplus$ of diameter $d = 4,5,\ldots,9$}
\vspace*{-2.5ex}
\end{center}

\newpage

\section{\hspace*{-5pt}Applications of tristance and quadristance anticodes}
\vspace{-1.25ex} 
\label{sec5}

\looseness=-1
Our study of tristance anticodes was originally motivated by
applications to multi-dimen\-sional interleaving~\cite{BBF,BBV,EV,ME}.
A two-dimensional {\dfn interleaving scheme\/}
$\cI(t,r)$ of strength~$t$ with $r$ repetitions
is a labeling $\cI(t,r): \Z^2 \to \{1,2,\ldots,\chi\}$ such that 
no integer~in~the~range $\{1,2,\ldots,\chi\}$ of $\cI(t,r)$ appears 
more than $r$ times among the labels of any 
connected subgraph of $\Gplus$ with $\le t$ vertices.
The integer $\chi$ is called the {\dfn interleaving degree\/} 
of $\cI(t,r)$~and denoted $\deg\cI(t,r)$.
Interleaving schemes of this kind may be used for error-control~in 
optical, holo\-graphic, and magnetic recording~\cite{BBF,BBV}.
The graphs $\Gstar$, $\Ghex$, and $\Gplusthree$ are also relevant
for these applications. In each case, the goal is to minimize
the interleaving~degree for a given strength $t$ and a given number of 
repetitions $r$. Usually, the values~of~interest are when $r$ is small
(say $r = 1,2,3$) and $t$ is large. 
We note that 
for $r = 1$, the problem has been completely solved in~\cite{BBV}.
For $r = 2,3$, {upper bounds} on the interleaving degree~are 
given in~\cite{EV}.
In particular, it is shown in~\cite{EV} that there
exist interleaving schemes $\cI(t,2)$ with
$\deg \cI(t,2) = (3/16)t^2 + O(t)$ and 
$\cI(t,3)$ with
$\deg \cI(t,3) = (8/81)t^2 + O(t)$.

Our results on tristance anticodes provide \emph{lower bounds}
on the minimum possible interleaving degree of\, $\cI(t,2)$
as follows.
If $\cA_d$ is a tristance anticode of diameter $d = t-1$, then 
any three points in $\cA_d$ belong to a connected subgraph with 
$\le t$ vertices, by definition. 
Therefore, no integer in the range of $\cI(t,2)$
can be used more than twice in labeling the points of $\cA_{t-1}$.
Hence $\deg\cI(t,2) \ge |\cA_{t-1}|/2$. 
In conjunction with Theorems \ref{optimal-tristance}, 
\ref{optimal-star-tristance}, \ref{optimal-hex-tristance}, 
this immediately implies that
\Eq{
\label{degree-bounds}
\deg\cI(t,2) 
\ \ge \
\left\{
\begin{array}{l@{\hspace{5ex}}l}
\displaystyle \ceilenv{\frac{t(t+1)}{7}} & 
\text{in the grid graph $\Gplus$}\\[4ex]
\displaystyle \ceilenv{\frac{t(t+1)}{6}} & 
\text{in the hexagonal graph $\Ghex$}\\[4ex]
\displaystyle \ceilenv{\frac{2t^2-1}{7}} & 
\text{in the infinity graph $\Gstar$}
\end{array}
\right.
}
For the grid graph $\Gplus$, a better bound was recently given 
in~\cite{ME} in the case where~$t$~is~even. In fact, it is shown
in~\cite{ME} that for even $t$ the minimum possible interleaving
degree of~$\cI(t,2)$ in $\Gplus$ is exactly
$\floorenv{(3t^2{+}\,4)/16}$ (the problem is still open 
for odd $t$). 
For $\Ghex$ and $\Gstar$, the bounds in~\eq{degree-bounds}
are the best known. Using similar reasoning,
\Tref{optimal-quadristance} implies that 
$$
\deg\cI(t,3) 
\ \ge \
\ceilenv{\frac{t(t+2)}{18}}
$$
in $\Gplus$. This bound is also the best known (cf.\,\cite{BBF}).
Finally, the results of \S4.1 herein imply a lower bound
on $\deg\cI(t,2)$ in the three-dimensional grid graph $\Gplusthree$.
For this graph, no upper bounds are yet known and even the 
problem of determining $\deg\cI(t,1)$ is still open.

\looseness=-1
Our results for the grid graph $\Gplus$ also have applications to 
the game of Go. Indeed, the game is played on a $19 \,{\times}\, 19$ 
square subgraph of $\Gplus$, called the goban. 
Two players~--\,Black 
and White -- alternate moves, each move consisting 
of one stone of the player's color being 
placed on one of the $361$ vertices of the goban. 
A set of stones of the same color~is~consid\-ered a connected group if the induced 
subgraph of $\Gplus$ is connected. Thus our results in \S2.2, \S2.3, 
and \S4.2 answer the following questions:\vspace{-6ex}
\begin{list}{}
{
\addtolength{\leftmargin}{-3.5ex}
\addtolength{\rightmargin}{4.5ex}
}
\item\noindent
\begin{itemize}
\item
How should three stones be played on an empty goban, so that 
they can then be all connected with at most $k$ moves?  
\vspace{.50ex}

\item
Given a stone, how should two stones be played on an empty
goban, so that all three stones can then be connected with 
at most $k$ moves?
\vspace{.50ex}

\item
Given two stones, where could one play a third stone so that
all three can then be connected with at most $k$ moves?
\vspace{.50ex}

\item
Given three stones, where could one play a fourth stone so that
all four can then be connected with at most $k$ moves?
\vspace{-1.50ex}
\end{itemize}
\end{list}
\looseness=-1
The answers to these questions are, respectively, the 
tristance anticode $\cA_{k+2}$ given~in~Theo\-rem\,\ref{optimal-tristance}, 
the centered tristance anticode $\cA_{k+2}(z_0)$ in \Tref{one-center}, 
the centered tristance anticode $\cA_{k+2}(z_1,z_2)$ in \Tref{centeroct},
and the centered quadristance anticode $\cA_{k+3}(z_1,z_2,z_3)$ 
in \Tref{quad-A_d}. It is interesting that the answers to the third
and fourth questions above 
are drastically different (compare Figures 2 and 11),
even though the questions themselves appear to be 
similar. Of course, all these results assume an empty goban and no
active opposition to the desired connection. Nevertheless, they 
could be of interest~for~computer
Go applications~\cite{Muller}.
We also have an algorithmic solution (to be presented elsewhere)
for the case where the goban already has black and white stones 
in arbitrary positions. 


\vfill
{\bf Acknowledgment.} We would like to thank 
Yael Merksamer for stimulating discussions.

\clearpage
\vspace{5ex}
%
%

\end{document}